\documentclass[10pt,a4paper,oneside]{amsart}
\usepackage{graphicx}
\usepackage{amsmath}
\usepackage{color}
\usepackage{enumerate}
\usepackage{amsthm}
\usepackage{latexsym}
\usepackage[active]{srcltx}
\newtheorem{theorem}{Theorem}[section]

\newtheorem{definition}[theorem]{Definition}

\newtheorem{lemma}[theorem]{Lemma}
\newtheorem{remark}[theorem]{Remark}
\newtheorem{proposition}[theorem]{Proposition}

\newtheorem{corollary}[theorem]{Corollary}

\numberwithin{equation}{section}

\newcommand{\tr}{\text{tr}}

\renewcommand{\H}{{\mathcal H}}
\def\C{\mathbb C}
\def\R{\mathbb R}

\def\C{\mathbb C}
\def\R{\mathbb R}
\def\I{\mathbb I}
\def\N{\mathbb N}
\def\al{\alpha}
\def\HH{\mathbb H}
\def\be{\beta}

\def\rh{\rho}
\def\et{\eta}
\def\th{\theta}
\def\ga{\gamma}

\def\la{\lambda}

\def\om{\omega}

\def\ta{\tau}

\def\g{\mathfrak{g}}

\def\h{\mathfrak{h}}

\def\ll{\mathfrak{l}}

\def\t{\mathfrak{t}}

\def\u{\mathfrak{u}}

\def\g{\mathfrak g}
\def\h{\mathfrak h}

\def\ga{\gamma}
\def\la{\lambda}

\def\si{\sigma}
\def\om{\omega}

\def\ph{\phi}
\def\ch{\chi}
\def\ta{\tau}
\def\ps{\psi}
\def\N{\mathbb{N}}
\def\Z{\mathbb{Z}}
\def\R{\mathbb{R}}
\def\C{\mathbb{C}}

\def\ol#1{\overline{#1}}
\def\nn{\nonumber}

\def\R{{\mathbb R}}
\def\C{{\mathbb C}}
\def\N{{\mathbb N}}

\def\Z{{\mathbb Z}}
\def\T{{\mathbb T}}
\def\I{{\mathbb I}}

  \def\Id{{\mathbb I}}

\def\B{{\mathcal B}}
\def\D{{\mathcal D}}
\def\F{{\mathcal F}}

\def\H{{\mathcal H}}
\def\P{{\mathcal P}}
\def\K{{\mathcal K}}
\def\LL{{\mathcal L}}

\def\RR{{\mathcal R}}

\def\W{{\mathcal W}}

\def\O{{\mathcal O}}

\def\Ad{{\text Ad}}
\def\tr{{\text tr}}

\def\iy{\infty}

\def\ol#1{\overline{#1}}

\def\hb#1{\hbox{#1}}
\def\val#1{\vert #1\vert}

\def\no#1#2{\Vert #1\Vert_{#2} }

\def\wh{\widehat}
\def\exp{\hb{exp}}
\def\ind{\hb{ind}}

\def\res#1{_{\vert #1}}
\def\inv{^{-1}}

\def\hb #1{\hbox{#1}}

\def\limk{\lim_{k\to\infty}}

\def\ti{\tilde}

\def\hb#1{\hbox{#1}}
\def\val#1{\vert #1\vert}

\def\ot{\otimes}

\def\cal{\mathcal}

\def\L1#1{L^1(#1)}

\def\L#1#2{L^{#1}(#2)}
\def\l#1#2{L^{#1}(#2)}

\def\Re{\mathrm{\, Re\, }}

\def\lef({\left(}
\def\rig){\right)}
\def\lan{\langle}
\def\ran{\rangle}

\usepackage{amssymb}

\begin{document}

\title[ On the dual topology of the groups $\mathbf{U(n)\ltimes\HH_n}$]{On the dual topology of the groups $\mathbf{U(n)\ltimes\HH_n}$}

\author{Mounir Elloumi, Janne-Kathrin G\"unther and Jean Ludwig}%
\address{Mathematics Department, College of Science, King Faisal University, P.O. 380, Ahsaa 31982, Kingdom of Saudi Arabia}
\email{melloumi@kfu.edu.sa}

\address{Universit\'e de Lorraine, Institut Elie Cartan de Lorraine, UMR 7502, Metz, F-57045, France \linebreak 
University of Luxembourg, Mathematical Research Unit, 6 Rue Richard Coudenhove-Kalergi, Luxembourg, L-1359, Luxembourg}
\email{janne.guenther@gmx.de}
\address{Universit\'e de Lorraine, Institut Elie Cartan de Lorraine, UMR 7502, Metz, F-57045, France}
\email{jean.ludwig@univ-lorraine.fr}
\thanks{2000 {\it Mathematics Subject Classification.} Primary 43A40; Secondary 22E45}
\keywords{Unitary group, semi-direct product, dual topology,
admissible coadjoint orbit space.} \vskip 0.2cm

\dedicatory{}

\begin{abstract}
Let $G_n=U(n)\ltimes \HH_n $ be the semi-direct product of the unitary group acting by automorphisms on the Heisenberg group $\HH_n$. According to Lipsman, the unitary dual  $\widehat {G_n} $
of $G_n $ is in one to one correspondence with the space of admissible coadjoint orbits $\mathfrak g_n^\ddagger /G_n $ of $G_n $. In this paper, we determine the topology of the space  $\mathfrak g_n^\ddagger /G_n $ and we show that the  correspondence with $\widehat {G_n} $ is a
homeomorphism.
\end{abstract}
\maketitle
\tableofcontents
\section{Introduction}
Let $G$ be a locally compact group and $\widehat{G}$ the unitary dual of $G$, i.e., the set of equivalence classes of irreducible unitary representations
 of $G$, endowed with the pullback of the hull-kernel topology on the primitive ideal space of $C^*(G)$, the $C^*$-algebra of $G$.
 Besides the fundamental problem of determining $\widehat{G}$ as a set, there is a genuine interest in a precise and neat description of the topology
 on $\widehat{G}$. For several classes of Lie groups, such as simply connected nilpotent Lie groups or, more generally, exponential solvable Lie groups,
 the Euclidean motion groups and also the extension groups $U(n)\ltimes\HH_n$ considered in this paper, there is a nice geometric object parameterizing
  $\widehat{G}$, namely the space of admissible coadjoint orbits in the dual $\g^*$ of the Lie algebra $\g$ of $G$.

In such a situation, the natural and important question arises  of whether the bijection between the orbit space, equipped with the quotient topology,
and $\widehat{G}$ is a homeomorphism. In [\ref{leplud}], H. Leptin and J. Ludwig have proved that for an exponential solvable Lie group $G=\exp\g$,
the dual space $\widehat{G}$ is homeomorphic to the space of coadjoint orbits $\g^*/G$ through the Kirillov mapping. On the other hand, it had been
 shown in [\ref{El-Lu}] that the dual topology of the classical motion groups $SO(n)\ltimes\R^n$ for $n\geq2$ can be linked to the topology of the quotient
 space of admissible coadjoint orbits.

In this paper, we consider the semi-direct product
$G_n=U(n)\ltimes\HH_n$ for $n\geq1$ and we identify its dual space
$\widehat{G_n}$ with the lattice of admissible coadjoint orbits. Lipsman
showed in [\ref{Lipsman}] that each irreducible unitary
representation of $G_n$ can be constructed by holomorphic induction
from an admissible linear functional $\ell$ of the Lie algebra
$\g_n$ of $G_n$. Furthermore, two irreducible representations in
$\widehat{G_n}$ are equi\-va\-lent if and only if their respective linear
functionals are in the same $G_n$-orbit. We prove  that this
identification is a homeomorphism.

This paper is structured in the following way. Section 2 contains
preliminary material and summarizes results from previous work
concerning the dual space of $G_n$ which is identified with the space of its
admissible coadjoint orbits. The representations attached to
an admissible linear functional are obtained via  Mackey's
little-group method and the dual space $\widehat{G_n} $ of $G_n$ is given by the
parameter space \mbox{$ \K_n=\Big\{\al\in\R^*, r>0,\rh_\mu\in\widehat{U(n-1)},
\ta_\la\in\widehat{U(n)}\Big\}$}. In Section 3, we shall link the
convergence of sequences of admissible coadjoint orbits to the
convergence in $ \K_n$. Section 4 describes the dual topology of a
second countable locally compact group. In the last two sections,  we
discuss the topology of the dual space of our groups $G_n$.

\section{The space of admissible coadjoint orbits }\label{prel-un}
Let $\C^n$ be the $n$-dimensional complex vector space equipped with the
standard scalar product $\lan .,.\ran_{\C^n}$ given by
$$\lan x,y \ran_{\C^n}=\sum \limits_{j=1}^{n} x_j \ol{y_j} \:\:\: \forall x=(x_1,...,x_n),y=(y_1,...,y_n) \in \C^n.$$
Moreover, let $(.,.)_{\C^n}$ and $\omega(.,.)_{\C^n}$ denote the real and the imaginary part of $\lan .,.\ran_{\C^n}$, respectively, i.e.
\begin{equation*}
\lan .,.\ran_{\C^n}=(.,.)_{\C^n}+i \omega(.,.)_{\C^n}.
\end{equation*}
The bilinear forms $(.,.)_{\C^n}$ and $\omega(.,.)_{\C^n}$ define a
positive definite inner product and a symplectic structure on the
underlying real vector space $\R^{2n}$ of $\C^n$, respectively. The associated
Heisenberg group \mbox{$\HH_n=\C^n\times\R$} of dimension $2n+1$ over $\R$
is given by the group multiplication
\begin{equation*}
(z,t)(z',t'):= \Big(z+z',t+t'-\frac{1}{2}\omega(z,z')_{\C^n} \Big) \:\:\: \forall (z,t),(z',t') \in \HH_n.
\end{equation*}
Furthermore, consider the unitary group $U(n)$ of automorphisms of $\HH_n$
preserving $\lan .,.\ran_{\C^n}$ on $\C^n$ which embeds into $\text{Aut}(\HH_n)$
via
\begin{equation*}
A.(z,t):=(Az,t) \:\:\: \forall A \in U(n)~ \forall (z,t) \in \HH_n.
\end{equation*}
Then, $U(n)$ yields a maximal compact connected subgroup of
$\text{Aut}(\HH_n)$ (see [\ref{Folland}], \mbox{Theorem 1.22} and [\ref{knapp}], Chapter I.1).
Moreover, $G_n=U(n)\ltimes\HH_n$
denotes the semi-direct product of $U(n)$ with the Heisenberg group
$\HH_n$ equipped with the group law
\begin{equation*}
(A,z,t)(B,z',t'):= \Big(AB,z+Az',t+t'-\frac{1}{2}\omega(z, Az')_{\C^n} \Big) \:\:\: \forall (A,z,t),(B,z',t') \in G_n.
\end{equation*}
The Lie algebra $\h_n$ of $\HH_n$ will be identified with $\HH_n$ itself via the
exponential map. The Lie bracket of $\h_n$ is defined as
\begin{eqnarray*}
\big[(z,t),(w,s) \big]:=\big(0,-\om(z,w)_{\C^n}\big) \:\:\: \forall (z,t),(w,s) \in \h_n
\end{eqnarray*}
and the derived action of the Lie algebra  $\u(n)$ of $U(n)$ on
$\h_n$ is
\begin{equation*}
A.(z,t):=(Az,0) \:\:\: \forall A \in \u(n) ~\forall (z,t) \in \h_n.
\end{equation*}
Denoting by $\g_n=\u(n)\ltimes\h_n$ the Lie algebra of $G_n$, for all $(A,z,t)\in G_n$ and all $(B,w,s)\in \g_n$, one gets
\begin{eqnarray} \label{unAd}
\nn  &  &
\Ad(A,z,t)(B,w,s)\\
&=&\frac{d}{dy}\Big|_{y=0}\Ad(A,z,t) \big(e^{yB},yw,ys \big) \\
\nn &=& \Big(ABA^*,-ABA^*z+Aw,s-\om(z,Aw)_{\C^n}+\frac{1}{2}\om(A^*z,BA^*z)_{\C^n}\Big),
\end{eqnarray}
where $A^*$ is the adjoint matrix of $A$. In particular
\begin{equation}\label{unea}
\Ad(A)(B,w,s)=(ABA^*,Aw,s).
\end{equation}
From Identity (\ref{unAd}), one can deduce the Lie bracket
\begin{eqnarray*}
\big[(A,z,t),( B,w,s)\big]&=&\frac{d}{dy}\Big|_{y=0}\Ad\big(e^{yA},yz,yt\big)(B,w,s)\\
&=&\big(AB-BA,Aw-Bz,-\om(z,w)_{\C^n}\big)
\end{eqnarray*}
for all $(A,z,t), (B,w,s)\in \g_n$. \\

\subsection{The coadjoint orbits of $G_n$}
~\\

In this subsection, the coadjoint orbit space of $G_n$ will be described according to [\ref{BGLR}], \mbox{Section 2.5.}\\

In the following, $\u(n)$ will be identified with its vector dual space $\u^*(n)$ with the help of the
\mbox{$U(n)$-invariant} inner product
\begin{eqnarray*}
\lan A,B\ran_{\u(n)}:=\tr(AB) \:\:\: \forall A,B \in \u(n).
\end{eqnarray*}
For $z\in\C^n$, define the linear form $z^\vee$ in $(\C^n)^*$ by
\begin{eqnarray*}
z^\vee(w):=\om(z,w)_{\C^n} \:\:\: \forall w \in \C^n.
\end{eqnarray*}
Furthermore, one defines a map $\times: \C^n\times\C^n\longrightarrow\u^*(n)$,
$(z,w)\mapsto z\times w$ by
\begin{eqnarray*}
z\times w(B):=w^\vee(Bz)=\om(w, Bz)_{\C^n} \:\:\: \forall B\in\u(n).
\end{eqnarray*}
One can verify that for $A\in U(n)$, $B\in \u(n)$ and $z,w\in
\C^n$,
\begin{eqnarray}
Az^\vee:=z^\vee\circ A\inv &=&(Az)^\vee,\label{uneb}\\
\nn z^\vee\circ B&=&-(Bz)^\vee,\\
\nn z\times w&=&w \times z \:\:\: \text{and} \\
\nn A(z\times w)A^*&=&(Az)\times (Aw).
\end{eqnarray}
Hence, the dual $\g^*_n= \big(\u(n)\ltimes\h_n \big)^*$ will be identified with
$\u(n)\oplus\h_n$, i.e. each element $\ell\in\g^*_n$ can be
identified with an element $(U,u,x)\in \u(n)\times\C^n\times\R$ such
that
\begin{eqnarray*}
\big\lan (U,u,x), (B,w,s)\big\ran_{\g_n}=\lan U, B\ran_{\u(n)} + u^\vee(w) + xs \:\:\:
\forall (B,w,s)\in\g_n.
\end{eqnarray*}
From (\ref{unea}) and (\ref{uneb}), one obtains for all $A \in U(n)$,
\begin{equation}\label{unorbeq}
\Ad^*(A)(U,u,x)=(AUA^*,Au,x) \:\:\: \forall (U,u,x)\in \u(n)\times\C^n\times\R
\end{equation}
and for all $(A,z,t) \in G_n$ and all $(U,u,x)\in \u(n)\times\C^n\times\R$,
\begin{equation}
\label{unorbitlem} \Ad^*(A,z,t)(U,u,x)= \Big(AUA^*+z\times
(Au)+\frac{x}{2}z\times z,Au+xz,x \Big).
\end{equation}
Letting $A$ and $z$ vary over $U(n)$ and $\C^n$, respectively, the
coadjoint orbit ${\cal O}_{(U,u,x)}$ of the linear form $(U,u,x)$ can then
be written as
\begin{equation*} 
{\cal O}_{(U,u,x)}=\bigg\{ \Big(AUA^*+z\times (Au)+\frac{x}{2}z\times
z,Au+xz,x \Big) \Big|\: A\in U(n), z\in\C^n \bigg\}
\end{equation*}
or equivalently, by replacing $z$ by $Az$ and using Identity
(\ref{unorbeq}),
\begin{equation*} 
{\cal O}_{(U,u,x)}=\bigg\{\Ad^*(A) \Big(U+z\times u+\frac{x}{2}z\times z,u+xz,x \Big) \Big|\:
A\in U(n), z\in\C^n \bigg\}.
\end{equation*}

Here, $z$ is regarded as a column vector $z=(z_1,\ldots,z_n)^{T}$ and
$z^*:=\ol{z}^t$. \\
One can show as follows that $z\times u\in\u^*(n)\cong\u(n)$ is the $n \times n$ skew-Hermitian matrix $\frac{i}{2}(uz^*+zu^*)$: \\
For all $B\in\u(n)$,
\begin{eqnarray*}
\big\lan uz^*+zu^*, B \big\ran_{\u(n)}=\tr \big((uz^*+zu^*)B\big)=\sum_{1\leq i,j\leq
n}B_{ji}z_i\ol{u}_j-\sum_{1\leq i,j\leq n}u_i\ol{B}_{ij}\ol{z}_j=-2i
z\times u(B).
\end{eqnarray*}
In particular, $z\times z$ is the skew-Hermitian matrix $i zz^*$
whose entries are determined by $(i zz^*)_{lj}=iz_l\ol{z}_j$.

\section{The spectrum of $G_n$}\label{duun}

\subsection{The spectrum and the admissible coadjoint orbits of $G_n$}
~\\

The description of the spectrum of $G_n$ is based on a method by Mackey (see [\ref{Ma}], \mbox{Chapter 10}), which states that one has to determine the irreducible unitary representations of the subgroup $\HH_n$ in order to construct representations of $G_n$. \\
~\\
First, regard the infinite-dimensional irreducible representations of the Heisenberg group $\HH_n$, which are parameterized by $\R^*$: \\
For each element $\al\in\R^*$, the coadjoint orbit ${\cal O}_\alpha$ of the
irreducible representation $\sigma_\al$ is the hyperplane
${\cal O}_\alpha=\big\{(z,\al)|\:z\in\C^n \big\}$. Since for every
$\alpha$, this orbit is invariant under the
action of $U(n)$, the unitary group $U(n)$ preserves the
equivalence class of $\sigma_\alpha$.

The representation $\sigma_\alpha$ can be realized for $\al>0 $ in the Fock space
\begin{equation*}
\F_\alpha(n)=\Bigg\{f:\C^n\longrightarrow\C \textrm{ entire } \Big|~
\int \limits_{\C^n}\vert f(w)\vert^2 e^{-\frac{\val{\alpha}}{2}\vert
w\vert^2} dw<\iy \Bigg\}
\end{equation*}
as
\begin{eqnarray*}
\sigma_\alpha(z,t)f(w):=e^{i\alpha
t-\frac{\alpha}{4}|z|^2-\frac{\alpha}{2}\langle
 w,z\rangle_{\C^n}}f(w+z)
\end{eqnarray*}
and for $\al<0 $ on the space $\F_\al(n) =\ol{\F_{-\al}(n)}$ as
\begin{eqnarray*}
 \sigma_\alpha(z,t) f(\ol{w}):=e^{i\alpha
t+\frac{\alpha}{4}|z|^2+\frac{\alpha}{2}\langle
 \ol{w},\overline{z}\rangle_{\C^n}}f(\ol{w}+\ol{z}).
 \end{eqnarray*}
See [\ref{Folland}], Chapter 1.6 or [\ref{Howe}], Section 1.7 for a discussion of the Fock space. \\
~\\
For each $A\in U(n)$, the operator
$W_\alpha(A)$ defined by
\begin{equation*}
W_\alpha(A):\F_{\alpha}(n)\rightarrow \F_{\alpha}(n),~W_\alpha(A)f(z):=f \big(A^{-1}z \big)\:\:\: \forall f \in \F_{\alpha}(n)~\forall z \in \C^n
\end{equation*}
intertwines $\sigma_\alpha$ and $(\sigma_\alpha)_A$ given by
$(\sigma_\alpha)_A(z,t):=\sigma_\alpha(Az,t)$. $W_\al$ is called the
projective intertwining representation of $U(n)$ on the Fock space.
Then, by [\ref{Ma}], Chapter 10, for each $\alpha\in\R^*$ and each element
$\ta_\la$ in $\widehat{U(n)}$,
\begin{equation*}
\pi_{(\la,\alpha)}(A,z,t):=\ta_\la(A)\otimes \big(\sigma_\alpha(z,t)\circ
W_\alpha(A) \big) \:\:\: \forall(A,z,t)\in G_n
\end{equation*}
is an irreducible unitary representation of $G_n$ realized in
$\H_{(\la,\al)}:=\H_\la\otimes\F_\alpha(n)$, where $\H_\la$ is the Hilbert space of
$\ta_\la$.

Associate to $\pi_{(\la,\alpha)}$ the linear functional
$\ell_{\la,\alpha}:=(J_\la,0,\alpha) \in \g_n^*$ given by
\begin{equation*}
J_\la:=\left( \begin{array}{ccc}
i\la_1 & \ldots & 0\\
\vdots & \ddots & \vdots\\
0 & \ldots & i\la_n\\
\end{array}\right).
\end{equation*}
Denote by $G_n[\ell_{\la,\alpha}]$, $U(n)[\ell_{\la,\alpha}]$ and
$\HH_n[\ell_{\la,\alpha}]$ the stabilizers of $\ell_{\la,\alpha}$ in $G_n$, $U(n)$ and $\HH_n$, respectively. By Formula
(\ref{unorbitlem}),
\begin{eqnarray*}
G_n[\ell_{\la,\alpha}]&=& \bigg\{(A,z,t)\in G_n\Big|~ \Big(A J_\la A^* +\frac{i}{2}\alpha zz^*,\alpha z,\alpha \Big)=(J_\la,0,\alpha) \bigg\} \\
&=& \big\{(A,0,t)\in G_n |~ A J_\la A^*=J_\la \big\}, \\
&&~\\
U(n)[\ell_{\la,\alpha}]&=& \big\{A\in U(n) |~ (A J_\la A^*,0,\alpha)=(J_\la,0,\alpha) \big\} \\
&=& \big\{A\in U(n)|~ A J_\la A^*=J_\la \big\} \:\:\:  \text{and} \\
\HH_n[\ell_{\la,\alpha}]&=& \bigg\{(z,t)\in \HH_n \Big|~ \Big(J_\la
+\frac{i}{2}\alpha zz^*,\alpha z,\alpha \Big)=(J_\la,0,\alpha) \bigg\}
~=~\{0\}\times\R.
\end{eqnarray*}
It follows that
$G_n[\ell_{\la,\alpha}]=U(n)[\ell_{\la,\alpha}]\ltimes\HH_n[\ell_{\la,\alpha}]$.
Hence, $\ell_{\la,\alpha}$ is aligned in the sense of Lipsman (see [\ref{Lipsman}], Lemma 4.2). \\
~\\
The finite-dimensional irreducible representations of $\HH_n$ are
the characters $\chi_v$ for $v\in\C^n$, defined by
\begin{eqnarray*}
\chi_{v}(z,t):=e^{-i(v,z)_{\C^n}} \:\:\: \forall (z,t) \in \HH_n.
\end{eqnarray*}
Denote by $U(n)_v$ the stabilizer of the character $\chi_v$, or
equivalently of the vector $v$, under the action of $U(n)$. Then, for every
irreducible unitary representation $\rho$ of $U(n)_v$, the tensor product
$\rho\otimes\chi_v$ is an irreducible representation of
$U(n)_v\ltimes\HH_n$ whose restriction to $\HH_n$ is a multiple of
$\chi_v$, and the induced representation
$$\pi_{(\rho,v)}:=\ind_{U(n)_v\ltimes\HH_n}^{U(n)\ltimes\HH_n}\rho\otimes\chi_v$$
is an irreducible representation of $G_n$. Furthermore, the restriction of
$\pi_{(\rho,v)}$ to $U(n)$ is equi\-va\-lent to $\ind_{U(n)_v}^{U(n)}\rho$. \\
For any $v'=Av$ for $A\in U(n)$ (i.e. $v$ and $v'$ belong to the same sphere
centered at $0$ and of radius $r=\Vert v\Vert_{\C^n}$), one has
$U(n)_{v'}=AU(n)_vA^*$ and thus, the representation $\pi_{(\rho,v)}$ is equivalent with $\pi_{(\rho',v')}$ for any $\rho' \in \widehat{U(n)_{v'}}$ such that $\rho'(B)=\rho(A^*BA)$ for each $B\in
U(n)_{v'}$. Hence, one can regard the character $\chi_r$ associated to
the linear form $v_r$ which is identified with the vector
$(0,\ldots,0,r)^T$ in $\C^n$. Throughout this text, denote by
$\rho_\mu$ the representation of the subgroup $U(n-1)=U(n)_{v_r}$
with highest weight $\mu$ and by $\pi_{(\mu,r)}$ the representation
$\pi_{(\rho_\mu,v_r)}=\ind_{U(n-1)\ltimes\HH_n}^{G_n} \rho_\mu\otimes\chi_r$. Its Hilbert space $\H_{(\mu,r)}$ is given by
$$\H_{(\mu,r)}=L^2 \Big(G_n/\big(U(n-1)\ltimes \HH_n \big),
\rho_\mu\otimes\chi_r \Big).$$

Again, $\pi_{(\mu,r)}$ is linked to the linear functional
$\ell_{\mu,r}:=(J_\mu,v_r,0) \in \g_n^*$ for
\begin{equation*}
J_\mu:= \left( \begin{array}{cccc}
i\mu_1 & \ldots & 0&0\\
\vdots & \ddots & \vdots&\vdots\\
0 & \ldots & i\mu_{n-1}&0\\
0 & \ldots & 0&0\\
\end{array}\right).
\end{equation*}
By (\ref{unorbitlem}), one can check that
\begin{eqnarray*}
G_n[\ell_{\mu,r}]&=& \Big\{(A,z,t)\in G_n \big|~ \big(A J_\mu A^*+z\times (Av_r),Av_r,0\big)=\big(J_\mu,v_r,0 \big) \Big\} \\
&=& \bigg\{(A,z,t)\in G_n \Big|~ A\in U(n-1),~ A J_\mu A^*+\frac{i}{2} \big(v_r z^*+z
v_r^* \big)=J_\mu \bigg\} \\
&=&  \big\{(A,z,t)\in G_n |~ z\in i\R v_r, A\in U(n-1),~ A J_\mu A^* =J_\mu
\big\},
\end{eqnarray*}
since $ A J_\mu A^*\in \u(n-1) $ and
\begin{equation*} 
v_r z^*+zv_r^*=\left(\begin{array}{cccc}
0 & \ldots & 0 & rz_1\\
\vdots & \ddots & \vdots & \vdots\\
0 & \ldots & 0 & rz_{n-1}\\
r\ol{z}_1 & \ldots & r \ol{z}_{n-1} & 2r\Re(z_n)\\
\end{array}
 \right).
\end{equation*}
In addition,
\begin{eqnarray*}
U(n)[\ell_{\mu,r}]&=& \big\{A\in U(n-1)|~A J_\mu A^*=J_\mu \big\} \:\:\:\:\:\: \text{and} \\
\HH_n[\ell_{\mu,r}]&=&i\R v_r\times\R.
\end{eqnarray*}
Thus, similarly to $\ell_{\la,\al}$, the linear functional $\ell_{\mu,r}$ is aligned. \\
~\\
One obtains in this way all the  irreducible unitary
representations of  $G_n$ which are not trivial on $\HH_n$. \\
~\\

The trivial extension of each element $\ta_\la$ of
$\widehat{U(n)}$ to the entire group $G_n$ is an irreducible
representation which will also be denoted by $\ta_\la$.

First of all, as $U(n)$ is a compact group, one knows that its spectrum is discrete and that every representation of $U(n)$ is finite-dimensional. \\

Now, let
\begin{eqnarray*}
\T_n=\left\{T=
\left( \begin{array}{ccc}
e^{i \th_1} & \ldots & 0\\
\vdots & \ddots & \vdots\\
0 & \ldots & e^{i \th_n}\\
\end{array}\right)
\Bigg|~ \th_j\in\R~\forall j \in \{1,...,n\} \right\}
\end{eqnarray*}
be a maximal torus of the unitary group $U(n)$ and let $\t_n$ be its
Lie algebra. By complexification of $\u(n)$ and $\t_n$, one gets the complex Lie algebras $\u^\C(n)=\g\ll(n,\C)=M(n,\C)$
and
\begin{eqnarray*}
\t_n^\C= \left\{H=\left( \begin{array}{ccc}
h_1 & \ldots & 0\\
\vdots & \ddots & \vdots\\
0 & \ldots & h_n\\
\end{array}\right)
\Bigg|~ h_j\in\C~ \forall j \in \{1,...,n \} \right\},
\end{eqnarray*}
respectively, which is a Cartan subalgebra of $\u^\C(n)$. For $j \in \{1,...,n \}$, define a linear functional $e_j$ by
\begin{equation*}
e_j\left( \begin{array}{ccc}
h_1 & \ldots & 0\\
\vdots & \ddots & \vdots\\
0 & \ldots & h_n\\
\end{array}\right):=h_j.
\end{equation*}
Let $P_n$ be the set of all dominant integral forms $\la$ for $U(n)$
which may be written in the form $\sum \limits_{j=1}^n i\la_j e_j$, or
simply as $\la=(\la_1,\cdots,\la_n)$, where $\la_j$ are integers for every $j \in \{1,...,n\}$ such that
$\la_1\geq\cdots\geq\la_n$. $P_n$ is a lattice in the
vector dual space $\t_n^*$ of $\t_n$.\\
Since each irreducible unitary representation $(\ta_\la, \H_{\la})$ of $U(n)$ is determined
by its highest weight $\la\in P_n$, the spectrum $\widehat{U(n)}$ of $U(n)$ is in bijection with the set $P_n$. \\
~\\
For each $\la$ in $P_n$, the highest vector $\phi^\la$ in the Hilbert space
$\H_\la$ of $\ta_\la$ verifies $\ta_\la(T)\phi^\la=\chi_\la(T)\phi^\la$, where $\chi_\la$ is the
character of $\T_n$ associated to the linear functional $\la$ and
defined by
\begin{equation*}
\chi_\la \left(T=\left( \begin{array}{ccc}
e^{i \th_1} & \ldots & 0\\
\vdots & \ddots & \vdots\\
0 & \ldots & e^{i \th_n}\\
\end{array}\right)\right):=e^{-i\la_1\th_1}\cdots
e^{-i\la_n\th_n}.
\end{equation*}
Moreover, for two irreducible unitary representations
$(\ta_\la,\mathcal{H}_{\la})$ and
 $(\ta_{\la'},\mathcal{H}_{\la'})$, the Schur orthogonality relation states that for all
$\xi,\eta\in\mathcal{H}_{\la}$ and all $\xi',\eta'\in\mathcal{H}_{\la'}$,
\begin{eqnarray}\label{e1}
\int \limits_{U(n)} \big\lan\ta_{\la}(g)\xi,\et \big\ran_{\H_{\la}}\ol{\big\lan\ta_{\la'}(g)\xi',\et' \big\ran_{\H_{\la'}}}~dg=\begin{cases}
0 & \textrm{if } \la\ne\la',\\
\frac{\langle \xi,\xi'\rangle_{\H_{\la}} \langle \eta',\eta\rangle_{\H_{\la}}}{d_\la} &
 \textrm{if }\la=\la',
\end{cases}
\end{eqnarray}
where $d_\la$ denotes the dimension of the representation $\ta_\la$. \\

The linear functional corresponding to the irreducible $G_n$-representation $\tau_\la $ for $\la \in P_n $ is given by  $\ell_\la:=(J_\la,0,0)$.
~\\

According to [\ref{hol-bie}], Chapter 1,
if $\rho_\mu$ is an irreducible representation of
$U(n-1)$ with highest weight $\mu=(\mu_1,...,\mu_{n-1})$, the
induced representation $\pi_\mu:=\ind_{U(n-1)}^{U(n)}\rho_\mu$ of
$U(n)$ decomposes with multiplicity one, and the representations of
$U(n)$ that appear in this decomposition are exactly those with the highest weight
$\la=(\lambda_1,...,\lambda_n)$ such that
\begin{eqnarray}\label{lamu}
\lambda_1\geq \mu_1\geq \lambda_2\geq \mu_2\geq...\geq
\lambda_{n-1}\geq \mu_{n-1}\geq \lambda_n.
\end{eqnarray}

Therefore, by Mackey's theory, the spectrum $\wh{G_n}$ consists of the following families of representations:
\begin{enumerate}[(i)]
\item $\pi_{(\la,\al)}$ for $\la \in P_n$ and $\al \in \R^*=\R \setminus \{0\}$,
\item $\pi_{(\mu,r)}$ for $\mu \in P_{n-1}$ and $r \in \R_{>0}~~$ and
\item $\tau_{\la}$ for $\la \in P_n$.
\end{enumerate}
Hence, $\wh{G_n}$ is in bijection with the set
$$(P_n\times\R^*) \cup  \big(P_{n-1}\times\R_{>0}\big) \cup P_n.$$
~\\
A linear functional $\ell$ in $\g_n^*$ is defined to be
admissible if there exists a unitary character $\chi$ of the
connected component of $G_n[\ell]$ such that
$d\chi=i\ell\res{\g_n[\ell]}$. A calculation shows that all the linear
functionals $\ell_{\la,\al}$, $\ell_{\mu,r}$ and $\ell_\la$ are
admissible. Then, according to [\ref{Lipsman}], the re\-pre\-sen\-tations
$\pi_{(\la,\alpha)}$, $\pi_{(\mu,r)}$ and $\ta_\la$ described above
are equivalent to the representations of $G_n$ obtained by
holomorphic induction from their respective linear functionals
$\ell_{\la,\al}$, $\ell_{\mu,r}$ and $\ell_\la$.

Denote by ${\cal O}_{(\la,\al)}$, ${\cal O}_{(\mu,r)}$ and  $
{\cal O}_\la $ the coadjoint orbits associated to the linear forms
$\ell_{\la,\al}$, $\ell_{\mu,r}$ and $\ell_\la$, respectively. Let $ \g_n^\ddagger
\subset \g_n^*$ be the union of all the elements in ${\cal O}_{(\la,\al)}$, ${\cal O}_{(\mu,r)} $ and ${\cal O}_{\la} $ and denote by $\g_n^\ddagger/G_n  $ the corresponding set in the orbit space. Now, from [\ref{Lipsman}] follows that $ \g_n^\ddagger $ is the set of all admissible linear functionals of $ \g_n $.

We obtain in this way the Kirillov-Lipsman bijection
\begin{eqnarray*}
 \K: \g_n^\ddagger/G_n &\to& \widehat{G_n},\\
\O_\ell &\mapsto& [\pi_\ell]
 \end{eqnarray*}
between the space of admissible coadjoint orbits and the space of equivalence classes of irreducible unitary representations of $G_n $.

\section{Convergence in the quotient space $\g_n^\ddagger/G_n$} \label{conv-quo}
According to the last subsection, the spectrum of $G_n$ is
parameterized by the dominant integral forms $\la$ for $U(n)$ and
$\mu$ for $U(n-1)$, the non-zero $\alpha\in\R$ attached to the
generic orbits ${\cal O}_\al$ in $\h_n^*$ and the positive real $r$ derived
from the natural action of the unitary group $U(n)$ on the
characters of the Heisenberg group $\HH_n$. \\
Moreover, it has been elaborated that the quotient space $\g_n^\ddagger/G_n$ of admissible coadjoint
orbits is in bijection with $\widehat{G_n}$. \\
Now, the convergence of the admissible coadjoint orbits will be linked to the convergence in the pa\-ra\-me\-ter space
$\Big\{\al\in\R^*, r>0,\rh_\mu\in\widehat{U(n-1)}, \ta_\la\in\widehat{U(n)} \Big\}$. \\
~\\
Letting $\W$ be the subspace of $\u(n)$ generated by the matrices
$z\times v_r=\frac{i}{2}(v_rz^*+zv_r^*)$ for $z\in \C^n$, the
space $\g_n^\ddagger/G_n$ is the set of all orbits
\begin{eqnarray*}
\mathcal{O}_{(\la,\alpha)}&=&\bigg\{ \Big(AJ_\la A^*+\frac{i\alpha}{2}
zz^*,\alpha z,\alpha \Big) \Big|~ z\in\C^n,~ A\in U(n) \bigg\},\\
\mathcal{O}_{(\mu,r)}&=&\Big\{\big(A(J_{\mu}+  \W)A^*,Av_r,0 \big) \big|~A\in U(n) \Big\} \:\:\: \text{and} \\
\mathcal{O}_\la&=& \big\{(AJ_{\lambda}A^*,0,0) \big|~ A\in U(n) \big\}
\end{eqnarray*}
for $\al \in \R^*$, $r \in \R_{>0}$, $\mu \in P_{n-1}$ and $\la \in P_n$. \\
~\\
Before beginning the discussion on the convergence of the admissible coadjoint orbits, the following preliminary lemmas are needed:

\begin{lemma}\label{fondegal2}
~\\
For $n \in \N^*$ and for any scalars
 $X_1,...,X_n$ and $Y_1,...,Y_{n-1}$ fulfilling $Y_i\ne Y_j$ for $i\ne j$, one has
\begin{eqnarray*}
\sum_{j=1}^{n-1}~\frac{\prod \limits^n_{\substack{i=1\\ i\ne k}}(X_i-Y_j)}{\prod \limits^{n-1}_{\substack{i=1\\ i\ne j}}(Y_i-Y_j)}~=~\sum_{\substack{j=1\\ j\ne k}}^n X_j-\sum_{j=1}^{n-1}Y_j
\end{eqnarray*}
for each $k \in \{1,...,n \}$.
\end{lemma}

Proof: \\
For $n=1$, the formula is trivial. \\
So, let $n>1$ and assume that the assertion is true for this $n$. Consider the relation at $n + 1 $.
\\
For $k=n+1$, a simple calculation gives the result. If $k\ne n+1$,
one gets

\begin{eqnarray*}
\sum \limits_{j=1}^{n}\frac{\prod \limits^{n+1}_{\substack{i=1\\ i\ne k}}(X_i-Y_j)}{\prod \limits^{n}_{\substack{i=1\\ i\ne j}}(Y_i-Y_j)}&=&\frac{\prod \limits^{n+1}_{\substack{i=1\\ i\ne k}}(X_i-Y_n)}{\prod \limits^{n-1}_{i=1}(Y_i-Y_n)}~+~\sum \limits_{j=1}^{n-1}\frac{\prod \limits^{n+1}_{\substack{i=1\\ i\ne k}}(X_i-Y_j)}{\prod \limits^{n}_{\substack{i=1\\ i\ne j}}
(Y_i-Y_j)}\\[2pt]
&=&(X_{n+1}-Y_n)~\frac{\prod \limits \limits^{n}_{\substack{i=1\\ i\ne k}}(X_i-Y_n)}{\prod \limits^{n-1}_{i=1}(Y_i-Y_n)}~+~\sum \limits_{j=1}^{n-1}\frac{\prod \limits^{n}_{\substack{i=1\\ i\ne k}}
(X_i-Y_j)}{\prod \limits^{n-1}_{\substack{i=1\\ i\ne j}}(Y_i-Y_j)} \cdot \frac{(X_{n+1}-Y_j)}{Y_n-Y_j}\\[2pt]
&=&(X_{n+1}-Y_n)~\frac{\prod \limits^{n}_{\substack{i=1\\ i\ne k}}(X_i-Y_n)}{\prod \limits^{n-1}_{i=1}(Y_i-Y_n)}~+~\sum \limits_{j=1}^{n-1}\frac{\prod \limits^{n}_{\substack{i=1\\ i\ne k}}
(X_i-Y_j)}{\prod \limits^{n-1}_{\substack{i=1\\ i\ne j}}(Y_i-Y_j)} \cdot \frac{(X_{n+1}-Y_n)}{Y_n-Y_j}
~+~\underbrace{\sum \limits_{j=1}^{n-1}\frac{\prod \limits^{n}_{\substack{i=1\\ i\ne k}}(X_i-Y_j)}{\prod \limits^{n-1}_{\substack{i=1\\ i\ne j}}(Y_i-Y_j)}}_{=\sum \limits_{\substack{j=1\\ j\ne k}}^n X_j-\sum \limits_{j=1}^{n-1}Y_j}\\
&=&(X_{n+1}-Y_n)~\underbrace{\sum \limits_{j=1}^{n}\frac{\prod \limits^n_{\substack{i=1\\ i\ne k}}(X_i-Y_j)}{\prod \limits^{n}_{\substack{i=1\\ i\ne j}}(Y_i-Y_j)}}_{=1 \textrm{ by [\ref{El-Lu}], Lemma 5.3}}~+~\sum \limits_{\substack{j=1\\ j\ne k}}^n X_j~-~\sum \limits_{j=1}^{n-1}Y_j~=~\sum \limits_{\substack{j=1\\ j\ne k}}^{n+1} X_j~-~\sum \limits_{j=1}^{n}Y_j
\end{eqnarray*}
and the claim is shown.
~\\
\qed

\begin{lemma}\label{unlemmatrix}
~\\
Let $\mu\in P_{n-1}$ and $\la\in P_n$. Then,
$\la_1\geq\mu_1\geq\lambda_2\geq...\geq\mu_{n-1}\geq\lambda_{n}$ if
and only if there is a skew-Hermitian matrix
\begin{equation*} 
B=\left( \begin{array}{cccccccc}
0&0& \ldots &0&-z_1\\
0&0& \ldots &0&-z_2\\
\vdots &\vdots&\ddots&\vdots&\vdots\\
0&0&\ldots&0&-z_{n-1}\\
\ol{z}_1&\ol{z}_2&\ldots &\ol{z}_{n-1}&ix\\
\end{array}\right)
\end{equation*}
in $\W$ such that $A(J_{\mu}+B)A^*=J_\la$ for an element $A\in U(n)$.
\end{lemma}

Proof: \\
For $y\in\R$, a computation shows that $\det(J_{\mu}+B-iy\I)=(-i)^nP(y)$, where
\begin{equation*}
P(y):=(y-x)\prod_{i=1}^{n-1}(y-\mu_i)~-~\sum_{j=1}^{n-1}
\bigg(\val{z_j}^2\prod_{\substack{i=1\\ i\ne j}}^{n-1}(y-\mu_i)\bigg).
\end{equation*}
Furthermore, one can observe that
$P(y)\overset{y \to \iy}{\longrightarrow}\infty$ and that $P(\mu_j)\leq0$
if $j$ is odd and $P(\mu_j)\geq0$ if $j$ is even. \\
Now, if
$A(J_{\mu}+B)A^*=J_\la$ for an element $A\in U(n)$, by the spectral theorem,
$i\la_1,i\la_2,\cdots,i\la_n$ are all the elements of the spectrum
of $J_\mu+B$ fulfilling
$\la_1\geq\mu_1\geq\lambda_2\geq...\geq\mu_{n-1}\geq\lambda_{n}$. \\
Conversely, suppose first that all the $\mu_j$ for $j \in \{1,...,n-1 \}$ are
pairwise distinct. In this case, let $B$ be the skew-Hermitian
matrix with the entries $z_1,\cdots,z_{n-1},x$ satisfying
\begin{eqnarray*}
\val{z_j}^2&=&-\frac{\prod \limits_{i=1}^{n}(\lambda_i-\mu_j)}{\prod \limits^{n-1}_{\substack{i=1\\ i\ne j}}(\mu_i-\mu_j)} \:\:\:\:\:\:\text{for every}~j \in \{1,...,n-1 \}\:\:\:\text{and} \\[1pt]
x&=&\sum_{j=1}^n\la_j-\sum_{j=1}^{n-1}\mu_j.
\end{eqnarray*}
From Lemma \ref{fondegal2},
\begin{eqnarray*}
P(\la_k)&=&\Bigg(\sum_{j=1}^{n-1}\mu_j-\sum_{\substack{j=1\\ j\ne k}}^n\la_j\Bigg)\prod_{i=1}^{n-1}(\lambda_k-\mu_i)~+~\sum_{j=1}^{n-1}\left(\frac{
\prod \limits_{\substack{i=1\\ i\ne k}}^{n}(\lambda_i-\mu_j)}{\prod \limits_{\substack{i=1\\ i\ne j}}^{n-1}(\mu_i-\mu_j)}\prod_{i=1}^{n-1}(\lambda_k-\mu_i)\right)\\[5pt]
&=&\prod_{i=1}^{n-1}(\lambda_k-\mu_i)\left(\sum_{j=1}^{n-1}\mu_j~-~\sum_{\substack{j=1\\ j\ne k}}^n\la_j~+~\sum_{j=1}^{n-1}\frac{\prod \limits_{\substack{i=1\\ i\ne k}}^{n}(\lambda_i-\mu_j)}{\prod \limits_{\substack{i=1\\ i\ne j}}^{n-1}(\mu_i-\mu_j)}
\right)=0.
\end{eqnarray*}
Hence, the spectrum of the matrix $J_\mu+B$ is the set
$\{i\la_1,i\la_2,\cdots,i\la_n\}$ and thus, the spectral theorem implies that $A(J_{\mu}+B)A^*=J_\la$ for an element $A\in U(n)$. \\
Now, if the $\mu_j$ for $j \in \{1,...,n-1 \}$ are not pairwise distinct, there exist two
families of integers $\{p_l|~1\leq l\leq s\}$ and $\{q_l|~1\leq l\leq s\}$ such that $1\leq p_1<
 q_1<p_2<q_2<\cdots<p_s<q_s\leq n-1$ and $\mu_{p_l}=\mu_{p_l+1}=\cdots=\mu_{q_l-1}=\mu_{q_l}$,
 $\mu_{q_l}\not=\mu_{q_l+1}$
and
 $\mu_{p_l-1}\not=\mu_{p_l}$ for all $l \in \{ 1,..., s\}$. Let
  \begin{eqnarray*}
Q(y):=\prod_{i=1}^{p_1}~\prod_{i=q_1+1}^{p_2}\cdots\prod_{i=q_s+1}^{n-1}
(y-\mu_i),~~~\tilde{Q}_l(y):=\prod_{\substack{i=1\\ i\ne p_l}}^{p_1}~\prod_{\substack{i=q_1+1 \\ i\ne p_l}}^{
p_2}\cdots\prod_{\substack{i=q_{s-1}+1 \\ i\ne p_l}}^{p_s}~\prod_{i=q_s+1}^{n-1}(y-\mu_i)
\end{eqnarray*}
\begin{eqnarray*}
\textrm{and}\:\:\:\:\:\:Q_j(y):=\prod_{\substack{i=1 \\ i\not=j}}^{p_1}~\prod_{\substack{i=q_1+1 \\ i\not=j}}^{p_2}
\cdots\prod_{\substack{i=q_s+1 \\ i\not=j}}^{n-1}
(y-\mu_i).
\end{eqnarray*}
Then, for
\begin{eqnarray*}
P(y):=(y-x)Q(y)~-~\sum_{l=1}^{s}\bigg(\sum_{j=p_l}^{q_l}\val{z_j}^2\bigg)\tilde{Q}_l(y)~-~\sum_{j=1}^{p_1-1}\sum_{j=q_1+1}^{p_2-1}...
\sum_{j=q_s+1}^{n-1}\Big(\val{z_j}^2Q_j(y)\Big),
\end{eqnarray*}
one gets $\det(J_{\mu}+B-iy\I)=(-i)^n\prod \limits_{l=1}^s (y-\mu_{p_l} )^{q_l-p_l}P(y)$. \\
Now, the entries $z_j$ of the skew-Hermitian matrix $B$ can be chosen as follows:
 \begin{eqnarray*}
\val{z_j}^2:=-\frac{\prod \limits_{i=1}^{n}(\lambda_i-\mu_j)}{\prod \limits_{\substack{i=1 \\ i\not=j}}^{n-1}(\mu_i-\mu_j)}~=~-\frac{\prod \limits_{i=1}^{p_1}\prod \limits_{i=q_1+1}^{p_2}
\cdots \prod \limits_{i=q_s+1}^{n}(\lambda_i-\mu_j)}{\prod \limits_{\substack{i=1 \\ i\not=j}}^{p_1} \prod \limits_{\substack{i=q_1+1 \\ i\not=j}}^{p_2} \cdots \prod \limits_{\substack{i=q_s+1 \\ i\not=j}}^{n-1}(\mu_i-\mu_j)}
\end{eqnarray*}
for each $j \in \{1,...,p_1-1,q_1+1,...,p_s-1,q_s+1,...,n-1\}$ and
\begin{eqnarray*}
\val{z_{p_l}}^2+...+\val{z_{q_l-1}}^2+\val{z_{q_l}}^2:=-\frac{\prod \limits_{i=1}^{p_1}\prod \limits_{i=q_1+1}^{p_2}\cdots\prod \limits_{i=q_s+1}^{n}(\lambda_i-\mu_{p_l})}
{\prod \limits_{\substack{i=1 \\ i\ne p_l}}^{p_1}
\prod \limits_{\substack{i=q_1+1 \\ i\ne p_l}}^{p_2}\cdots
\prod \limits_{\substack{i=q_{s-1}+1 \\ i\ne p_l}}^{p_s}~\prod \limits_{i=q_s+1}^{n-1}(\mu_i-\mu_{p_l})}
\end{eqnarray*}
for each $l \in \{1,...,s \}$. The entry $x$ can be defined as
\begin{equation*}
x:=\sum_{j=1}^n\la_j-\sum_{j=1}^{n-1}\mu_j=\sum_{j=1}^{p_1}\sum_{j=q_1+1}^{p_2}...\sum_{j=q_s+1}^{n}\la_j~-~\sum_{j=1}^{p_1}\sum_{j=q_1+1}^{p_2}
...\sum_{j=q_s+1}^{n-1}\mu_j.
\end{equation*}
Then, if $\lambda_k=\mu_{p_l}$, one obviously has $P(\lambda_k)=Q(\lambda_k)=0$ and the multiplicity of the root $\la_k = \mu_{p_l} $ of $P $ is
$q_l - p_l $.
Otherwise, one gets
\begin{eqnarray*}
P(\lambda_k)&=&\bigg(\la_k-\sum_{j=1}^{p_1}\sum_{j=q_1+1}^{p_2}\cdots\sum_{j=q_s+1}^{n}\la_j+\sum_{j=1}^{p_1}\sum_{j=q_1+1}^{p_2}\cdots
\sum_{j=q_s+1}^{n-1}\mu_j\bigg)Q(\la_k)\\[4pt]
&&+~\sum_{j=1}^{p_1-1}\sum_{j=q_1+1}^{p_2-1}\cdots\sum_{j=q_s+1}^{n-1}\Bigg(~\frac{\prod \limits_{i=1}^{p_1}\prod \limits_{i=q_1+1}^{p_2}
\cdots\prod \limits_{i=q_s+1}^{n}(\lambda_i-\mu_j)}{\prod \limits_{\substack{i=1 \\ i\not=j}}^{
p_1}\prod \limits_{\substack{i=q_1+1 \\ i\not=j}}^{p_2}\cdots\prod \limits_{\substack{i=q_s+1 \\ i\not=j}}^{n-1}(\mu_i-\mu_j)}~Q_j(\la_k)\Bigg)\\[3pt]
&&+~\sum_{l=1}^{s}\Bigg(~\frac{\prod \limits_{i=1}^{p_1}\prod \limits_{i=q_1+1}^{p_2}\cdots\prod \limits_{i=q_s+1}^{n}(\lambda_i-\mu_{p_l})}{\prod \limits_{\substack{i=1 \\ i\ne p_l}}^{p_1}
\prod \limits_{\substack{i=q_1+1 \\ i\ne p_l}}^{p_2}\cdots\prod \limits_{\substack{i=q_{s-1}+1 \\ i\ne p_l}}^{p_s}~\prod \limits_{i=q_s+1}^{n-1}(\mu_i-\mu_{p_l})}~\tilde{Q}_l(\la_k)
\Bigg)\\
\end{eqnarray*}
\begin{eqnarray*}
&=&Q(\la_k)\Bigg(\sum_{j=1}^{p_1}\sum_{j=q_1+1}^{p_2}\cdots\sum_{j=q_s+1}^{n-1}\mu_j-\sum_{\substack{j=1 \\ j\ne k}}^{p_1}\sum_{\substack{j=q_1+1 \\ j\ne k}}^{p_2}\cdots\sum_{\substack{j=q_s+1 \\ j\ne k}}^{n}\la_j\\[4pt]
&&+~\sum_{j=1}^{p_1-1}\sum_{j=q_1+1}^{p_2-1}\cdots\sum_{j=q_s+1}^{n-1}\frac{\prod \limits_{\substack{i=1 \\ i\ne k}}^{p_1}
\prod \limits_{\substack{i=q_1+1 \\ i\ne k}}^{p_2}\cdots\prod \limits_{\substack{i=q_s+1 \\ i\ne k}}^{n}(\lambda_i-\mu_j)}{\prod \limits_{\substack{i=1 \\ i\ne j}}^{p_1}\prod \limits_{\substack{i=q_1+1 \\ i\not=j}}^{p_2}\cdots\prod \limits_{\substack{i=q_s+1 \\ i\not=j}}^{n-1}(\mu_i-\mu_j)}\\[3pt]
&&+~\sum_{l=1}^{s}\frac{\prod \limits_{\substack{i=1 \\ i\ne k}}^{p_1}\prod \limits_{\substack{i=q_1+1 \\ i\ne k}}^{p_2}\cdots\prod \limits_{\substack{i=q_s+1 \\ i\ne k}}^{n} (\lambda_i-\mu_{p_l})}{\prod \limits_{\substack{i=1 \\ i\ne p_l}}^{p_1}\prod \limits_{\substack{i=q_1+1 \\ i\ne p_l}}^{p_2}\cdots\prod \limits_{\substack{i= q_{s-1}+1 \\ i\ne p_l}}^{p_s}~\prod \limits_{i=q_s+1}^{n-1}(\mu_i-\mu_{p_l})}~\Bigg)\\[5pt]
&=&Q(\la_k)\Bigg(\sum_{j=1}^{p_1}\sum_{j=q_1+1}^{p_2}\cdots\sum_{j=q_s+1}^{n-1}\mu_j-\sum_{\substack{j=1 \\ j\ne k}}^{p_1}
\sum_{\substack{j=q_1+1 \\ j\ne k}}^{p_2}\cdots\sum_{\substack{j=q_s+1 \\ j\ne k}}^{n}\la_j\\[4pt]
&&+~\sum_{j=1}^{p_1}\sum_{j=q_1+1}^{p_2}\cdots\sum_{j=q_s+1}^{n-1}\frac{\prod \limits_{\substack{i=1 \\ i\ne k}}^{p_1}
\prod \limits_{\substack{i=q_1+1 \\ i\ne k}}^{p_2}\cdots\prod \limits_{\substack{i=q_s+1 \\ i\ne k}}^{n}(\lambda_i-\mu_j)}{\prod \limits_{\substack{i=1 \\ i\not=j}}^{p_1}\prod \limits_{\substack{i=q_1+1 \\ i\not=j}}^{p_2}\cdots\prod \limits_{\substack{i=q_s+1 \\ i\not=j}}^{n-1}(\mu_i-\mu_j)}~\Bigg)=0.
\end{eqnarray*}
Hence, the spectrum of the matrix $J_\mu+B$ equals the set
$\{i\la_1,i\la_2,\cdots,i\la_n\}$. As above, this completes the
proof.
~\\
\qed

\begin{lemma}\label{unspctrela}
~\\
\begin{enumerate}
\item Let $\la\in P_n$, $\alpha\in\R^*$ and $z\in \C^n$. Then, the matrix $J_\la+\frac{i}{\alpha}zz^*$ admits $n$
eigenvalues $i\be_1,i\be_2,\ldots,i\be_n$ such that $\be_1\geq\la_1\geq\be_2\geq\la_2\geq\cdots\geq\be_n\geq\la_n$ if $\alpha>0$ and
$\la_1\geq\be_1\geq\la_2\geq\be_2\geq\cdots\geq\la_n\geq\be_n$ if $ \alpha<0$.
\item Let $\la, \be \in P_n$. If $\be_1 \geq \la_1 \geq \be_2 \geq \la_2\geq ...\geq \be_n \geq \la_n$, there exists a number \mbox{$z \in \C^n$} such that the matrix $J_{\la}+izz^*$ admits the $n$ eigenvalues $i \be_1,...,i \be_n$. If \mbox{$\la_1 \geq \be_1 \geq \la_2 \geq \be_2\geq ...\geq \la_n \geq \be_n$}, there exists $z \in \C^n$ such that the matrix $J_{\la}-izz^*$ admits the $n$ eigenvalues $i \be_1,...,i \be_n$.
\end{enumerate}
\end{lemma}

Proof: \\
1) One can prove by induction that the characteristic polynomial of the matrix $\frac{1}{i}J_\la+\frac{zz^*}{\alpha}$ is equal to
$Q_n^{\la,z,\alpha}$ defined by
\begin{equation*}
Q_n^{\la,z,\alpha}(x):=\prod_{i=1}^n(x-\la_i)~-~\sum_{j=1}^n\frac{\val{z_j}^2}{\alpha}\prod_{\substack{i=1 \\ i\ne j}}^n(x-\la_i).
\end{equation*}
Assume that $\alpha$ is negative. Then, $Q_n^{\la,z,\alpha}(x)\overset{x \to \iy}{\longrightarrow}\infty$ and $Q_n^{\la,z,\alpha}(\la_j)\geq0$ if $j$ is odd and
$Q_n^{\la,z,\alpha}(\la_j)\leq0$ if $j$ is even. Furthermore, $Q_n^{\la,z,\alpha}(x)\overset{x \to \iy}{\longrightarrow}-\iy$ if
$n$ is odd and
$Q_n^{\la,z,\alpha}(x)\overset{x \to -\iy}{\longrightarrow}\iy$ if
$n$ is even and therefore, one can deduce that the matrix $\frac{1}{i}J_\la+\frac{zz^*}{\alpha}$ admits $n$
eigenvalues $\be_1,\be_2,\ldots,\be_n$ verifying
$\la_1\geq\be_1\geq\la_2\geq\be_2\geq\cdots\geq\la_n\geq\be_n$. Hence, $J_\la+\frac{i}{\alpha}zz^*$ admits the $n$ eigenvalues $i\be_1,i\be_2,\ldots,i\be_n$ fulfilling $\la_1\geq\be_1\geq\la_2\geq\be_2\geq\cdots\geq\la_n\geq\be_n$. \\
The same reasoning applies when $\alpha$ is positive. \\

2) Let $\be_1 \geq \la_1 \geq \be_2 \geq \la_2\geq ...\geq \be_n \geq \la_n$. \\
For any $z \in \C^n$, the characteristic polynomial of $\frac{1}{i}J_{\la} +zz^*$ is equal to $Q_n^{\la,z,1}$ with $Q_n^{\la,z,1}=:Q_n^{\la,z}$ like above. \\
First, assume that $\be_1 > \la_1 >  ...> \be_n > \la_n$. \\
Let
$$|z_j|^2:=- \frac{\prod \limits_{i=1}^n(\la_j-\be_i)}{\prod \limits_{\substack{i=1\\ i\not=j}}^n (\la_j-\la_i)}.$$
Then, as $\la_j<\be_i$ for all $i \in \{1,...,j\}$, as $\la_j>\be_i$ for all $i \in \{j+1,...,n\}$, as $\la_j<\la_i$ for all $i \in \{1,...,j-1\}$ and as $\la_j>\la_i$ for all $i \in \{j+1,...,n\}$, one gets $\text{sgn} \big(|z_j|^2 \big)=(-1) \frac{(-1)^j}{(-1)^{j-1}}=1$ and thus, this definition is meaningful. \\
One now has to show that $Q_n^{\la,z}(\be_{\ell})=0$ for all $\ell \in \{1,...,n\}$.
\begin{eqnarray*}
Q_n^{\la,z}(\be_{\ell})&=&\prod \limits_{i=1}^n (\be_{\ell}-\la_i)~+~ \sum \limits_{j=1}^n \frac{\prod \limits_{i=1}^n(\la_j-\be_i)}{\prod \limits_{\substack{i=1\\ i\not=j}}^n (\la_j-\la_i)} ~\prod \limits_{\substack{i=1\\ i\not=j}}^n (\be_{\ell}- \la_i) \\[1pt]
&=& \prod \limits_{i=1}^n (\be_{\ell}-\la_i)~ \left(1+\sum \limits_{j=1}^n \frac{\prod \limits_{i=1}^n(\la_j-\be_i)}{\prod \limits_{\substack{i=1\\ i\not=j}}^n (\la_j-\la_i)(\be_{\ell}-\la_j)} \right)\\[3pt]
&=& \prod \limits_{i=1}^n (\be_{\ell}-\la_i) ~\left(1-\sum \limits_{j=1}^n \frac{\prod \limits_{\substack{i=1\\ i\not=\ell}}^n(\la_j-\be_i)}{\prod \limits_{\substack{i=1\\ i\not=j}}^n (\la_j-\la_i)} \right)=0,
\end{eqnarray*}
as by [\ref{El-Lu}], Lemma 5.3, one obtains $\sum \limits_{j=1}^n \frac{\prod \limits_{\substack{i=1\\ i\not=\ell}}^n(\la_j-\be_i)}{\prod \limits_{\substack{i=1\\ i\not=j}}^n (\la_j-\la_i)}=1$. \\
~\\
Now, regard arbitrary $\be_1 \geq \la_1 \geq \be_2 \geq \la_2\geq ...\geq \be_n \geq \la_n$. \\
For $n=1$, one can choose $|z_1|^2:=(\be_1-\la_1)\geq 0$ and the claim is shown. \\
Let $n>1$ and assume that the assertion is true for $n-1$. \\
If $\la_{\ell-1} \not=\be_{\ell} \not=\la_{\ell}$ for all $\ell \in \{1,...,n\}$, the claim is already shown above. So, without restriction let $\ell \in \{1,...,n\}$ with $\be_{\ell} =\la_{\ell}$. The case $\la_{\ell-1}=\be_{\ell}$ is very similar. \\
Hence, for $\la^{\ell}:=(\la_1,...,\la_{\ell-1},\la_{\ell+1},...,\la_n)$ and $\be^{\ell}:=(\be_1,...,\be_{\ell-1},\be_{\ell+1},...,\be_n)$,
$$\be_1 \geq \la_1 \geq...\geq \be_{\ell-1} \geq \la_{\ell-1} \geq \be_{\ell+1} \geq \la_{\ell+1} \geq ...\geq \be_n \geq \la_n$$
and thus, by the induction hypothesis, there exists $\C^{n-1} \ni z^{\ell}:=(z_1,...,z_{\ell-1},z_{\ell+1},...,z_n)$ such that $Q_{n,\ell}^{\la^{\ell},z^{\ell}}(\be_i)=0$ for all
$i \in \{1,...,\ell-1,\ell+1,...,n\}$, where
$$Q_{n,\ell}^{\la^{\ell},z^{\ell}}(x):=\prod \limits_{\substack{i=1\\ i\not=\ell}}^n (x-\la_i)~-~ \sum \limits_{\substack{j=1\\ j\not=\ell}}^n |z_j|^2 \prod \limits_{\substack{i=1\\ i\not=j,i\not=\ell}}^n (x- \la_i).$$
Now, let $z_{\ell}:=0$, i.e. $z:=(z_1,...,z_{\ell-1},0,z_{\ell+1},...,z_n)$. Then,
\begin{eqnarray*}
Q_{n}^{\la,z}(x)&=&(x-\la_{\ell})\prod \limits_{\substack{i=1\\ i\not=\ell}}^n (x-\la_i)~-~ \sum \limits_{\substack{j=1\\ j\not=\ell}}^n |z_j|^2 (x- \la_{\ell})
\prod \limits_{\substack{i=1\\ i\not=j,i\not=\ell}}^n (x- \la_i)~-~|z_{\ell}|^2\prod \limits_{\substack{i=1\\ i\not=\ell}}^n (x- \la_i) \\[1pt]
&=&(x-\la_{\ell})Q_{n,\ell}^{\la^{\ell},z^{\ell}}(x)~-~|z_{\ell}|^2\prod \limits_{\substack{i=1\\ i\not=\ell}}^n (x- \la_i) \\
&=&(x-\la_{\ell})Q_{n,\ell}^{\la^{\ell},z^{\ell}}(x).
\end{eqnarray*}
If $i \in \{1,...,\ell-1,\ell+1,...,n\}$, then $Q_{n,\ell}^{\la^{\ell},z^{\ell}}(\be_i)=0$ and thus, $Q_{n}^{\la,z}(\be_i)=0$. Furthermore, \mbox{$Q_{n}^{\la,z}(\be_{\ell})=0$}, as $\be_{\ell}=\la_{\ell}$. \\
Therefore, $Q_{n}^{\la,z}(\be_i)=0$ for all $i \in \{1,...,n\}$ and the claim is shown. \\
~\\
Next, let $\la_1 \geq \be_1 \geq \la_2 \geq \be_2\geq ...\geq \la_n \geq \be_n$. \\
Then, for any $z \in \C^n$, the characteristic polynomial of $\frac{1}{i}J_{\la} -zz^*$ is equal to $Q_n^{\la,z,-1}$. \\
If $\la_1 > \be_1 >  ...> \la_n > \be_n$, let
$$|z_j|^2:=\frac{\prod \limits_{i=1}^n(\la_j-\be_i)}{\prod \limits_{\substack{i=1\\ i\not=j}}^n (\la_j-\la_i)}.$$
Here, $\text{sgn} \big(|z_j|^2 \big)= \frac{(-1)^{j-1}}{(-1)^{j-1}}=1$ and hence, this definition is meaningful. \\
The rest of the proof is the same as in the first part of (2).
~\\
\qed
~\\
~\\
With these lemmas, one can now prove the following theorem which describes the topology of the space of admissible coadjoint orbits of $G_n$.

\begin{theorem}\label{unconorballam}
~\\
Let $\al\in\R^*$, $r>0$, $\mu\in P_{n-1}$ and $\la\in P_n$. Then, the following holds:
\begin{enumerate}
\item A sequence of coadjoint orbits
$\big(\mathcal{O}_{(\la^k,\alpha_k)}\big)_{k \in \N}$ converges to the orbit
$\mathcal{O}_{(\lambda,\alpha)}$  in $\g_n^\ddagger/G_n$ if and only
if $\underset{k\to\infty}{\lim}\alpha_k=\alpha$ and $\la^k=\la$ for
large $k$.
\item A sequence of coadjoint orbits
$\big(\mathcal{O}_{(\la^k,\alpha_k)} \big)_{k\in\mathbb{N}}$ converges to the
orbit $\mathcal{O}_{(\mu,r)}$  in $\g_n^\ddagger/G_n$ if and only if
$\underset{k\to\infty}{\lim}\alpha_k=0$ and the sequence
$\big(\mathcal{O}_{(\la^k,\alpha_k)}\big)_{k\in\mathbb{N}}$ satisfies one of
the following conditions:
\begin{enumerate}[(i)]
\item For $k$ large enough, $\alpha_k>0$, $\la_j^k=\mu_j$ for all
$j \in \{1,..., n-1\}$ and
$\underset{k\to\infty}{\lim}\alpha_k\la^k_{n}=-\frac{r^2}{2}$.
\item For $k$ large enough, $\alpha_k<0$, $\la_j^k=\mu_{j-1}$ for all
$j \in \{2,..., n\}$ and
$\underset{k\to\infty}{\lim}\alpha_k\la^k_{1}=-\frac{r^2}{2}$.
\end{enumerate}
\item A sequence of coadjoint orbits
$\big(\mathcal{O}_{(\la^k,\alpha_k)}\big)_{k\in\mathbb{N}}$ converges to the
orbit $\mathcal{O}_\la$  in $\g_n^\ddagger/G_n$ if and only if
$\underset{k\to\infty}{\lim}\alpha_k=0$ and the sequence
$\big(\mathcal{O}_{(\la^k,\alpha_k)}\big)_{k\in\mathbb{N}}$ satisfies one of
the following conditions:
\begin{enumerate}[(i)]
\item For $k$ large enough, $\alpha_k>0$, $\la_1\geq\la_1^k\geq...\geq\la_{n-1}\geq\la_{n-1}^k\geq\la_n\geq\la_n^k$ and $\underset{k\to\infty}{\lim}\alpha_k\la^k_{n}=0$.
\item For $k$ large enough, $\alpha_k<0$, $\la_1^k\geq\la_1\geq\la_2^k\geq\la_2\geq...\geq\la_{n-1}\geq\la_{n}^k\geq\la_n$ and \mbox{$\underset{k\to\infty}{\lim}\alpha_k\la^k_{1}=0$}.
\end{enumerate}
\item A sequence of coadjoint orbits $\big(\mathcal{O}_{(\mu^k,r_k)}\big)_{k \in \N}$
converges to the orbit $\mathcal{O}_{(\mu,r)}$ in $\g^\ddagger_n/G_n$ if and
only if $\underset{k\to\infty}{\lim}r_k=r$ and $\mu^k=\mu$ for large
$k$.
\item A sequence of coadjoint orbits $\big(\mathcal{O}_{(\mu^k,r_k)}\big)_{k \in \N}$
converges to $\mathcal{O}_{\la}$ in $\g^\ddagger_n/G_n$ if and only
if $(r_k)_{k \in \N}$ tends to $0$ and
$\la_1\geq\mu^k_1\geq\la_2\geq\mu^k_2\geq...\geq\la_{n-1}\geq\mu^k_{n-1}\geq
\la_n $ for $k$ large enough.

\item A sequence of coadjoint orbits $\big(\mathcal{O}_{\la^k}\big)_{k \in \N}$
converges to the orbit $\mathcal{O}_{\lambda}$  in
$\g_n^\ddagger/G_n$ if and only if $\la^k=\la$ for large $k$.
\end{enumerate}
\end{theorem}

Proof: \\
Examining the shape of the coadjoint orbits listed at the beginning of this subsection, 1) and 6) are clear and \mbox{Assertion 5)} follows immediately from
Lemma \ref{unlemmatrix}. Furthermore, the proof of 4) is similar to that of [\ref{El-Lu}], Theorem 4.2. \\

2) Assume that
$\big(\mathcal{O}_{(\la^k,\alpha_k)}\big)_{k\in\mathbb{N}}$ converges to the
orbit $\mathcal{O}_{(\mu,r)}$. Then, there exist a sequence
$(A_k)_{k\in\N}$ in $U(n)$ and a sequence of vectors
$\big(z(k)\big)_{k\in\N}$ in $\C^n$ such that
\begin{equation*}
\underset{k\to\infty}{\lim} \bigg(A_k \Big(J_{\la^k} +\frac{i}{\alpha_k}
z(k)z(k)^* \Big)A_k^*,\sqrt{2}A_k z(k),\alpha_k \bigg)=(J_\mu,v_r,0).
\end{equation*}
Let $A=(a_{mj})_{1\leq m,j\leq n}$ be the limit of a subsequence
$(A_s)_{s\in I}$ for $I \subset \N$. Then,
$$\underset{s\to\infty}{\lim}J_{\la^s} +\frac{i}{\alpha_s}z(s)z(s)^*=A^*J_\mu A, \:\:\: \underset{s\to\infty}{\lim}z_j(s)=\frac{r}{\sqrt{2}}~\ol{a}_{nj}~
\text{for} \: j \in \{1,..., n\} \:\:\: \text{and} \:\:\: \lim \limits_{s \to \iy} \al_s=0.$$
On the other hand, one has $(A^*J_\mu
A)_{mj}=i\sum \limits_{l=1}^{n-1}\mu_l\ol{a}_{lm}a_{lj}$ and
\begin{eqnarray*}
J_{\la^s} +\frac{i}{\alpha_s} z(s)z(s)^*=\left(
\begin{array}{cccccccc}
i\la^s_1+i\frac{\val{z_1(s)}^2}{\alpha_s}&i\frac{z_1(s)\ol{z}_2(s)}{\alpha_s}& \ldots &i\frac{z_1(s)\ol{z}_n(s)}{\alpha_s}\\
i\frac{z_2(s)\ol{z}_1(s)}{\alpha_s}&i\la^s_2+i\frac{\val{z_2(s)}^2}{\alpha_s}& \ldots &i\frac{z_2(s)\ol{z}_n(s)}{\alpha_s}\\
\vdots &\vdots&\ddots&\vdots\\
i\frac{z_n(s)\ol{z}_1(s)}{\alpha_s}&i\frac{z_n(s)\ol{z}_2(s)}{\alpha_s}&\ldots&i\la^s_n+i\frac{\val{z_n(s)}^2}{\alpha_s}\\
\end{array}\right).
\end{eqnarray*}
Hence, for $m\ne j$,
$\underset{s\to\infty}{\lim}\big|\frac{z_m(s)\ol{z}_j(s)}{\alpha_s}\big|=\Big|\sum \limits_{l=1}^{n-1}\mu_l\ol{a}_{lm}a_{lj}\Big|<\iy$,
and since $\underset{s\to\infty}{\lim}\Vert
z(s)\Vert_{\C^n}=\frac{r}{\sqrt{2}}\ne0$, there is a unique $i_0 \in \{1,...,n\}$ such that
$\underset{s\to\infty}{\lim}z_{i_0}(s)=\frac{r}{\sqrt{2}}e^{i\theta}$
for a certain number $\theta\in\R$ and $\underset{s\to\infty}{\lim}z_j(s)=0$ for $j\ne
i_0$. One obtains $a_{ni_0}=e^{-i\theta}$ and $a_{nj}=0$ for $j\ne
i_0$, i.e. the matrices $A$ and $A^*J_\mu A$ can be written in the
following way:
\begin{eqnarray*}
A=\underset{\underbrace{}_{i_0-{\rm{th}}\rm{\:position}}}{\left(
\begin{array}{cccccccc}
*&\cdots&*&0&*&\cdots&*\\
\vdots& &\vdots&0&\vdots& &\vdots\\
\vdots& &\vdots&\vdots&\vdots& &\vdots\\
\vdots& &\vdots&0&\vdots& &\vdots\\
*&\cdots&*&0&*&\cdots&*\\
0&\cdots&0&e^{-i\theta}&0&\cdots&0\\
\end{array}\right)}
~~~\text{and}
\end{eqnarray*}
\begin{eqnarray*}
A^*J_\mu
A=\underset{\underbrace{}_{i_0-{\rm{th}}\rm{\:position}}}{\left(
\begin{array}{cccccccc}
*&\cdots&*&0&*&\cdots&*\\
\vdots& &\vdots&\vdots&\vdots& &\vdots\\
*&\cdots&*&0&*&\cdots&*\\
0&\cdots&0&0&0&\cdots&0\\
*&\cdots&*&0&*&\cdots&*\\
\vdots& &\vdots&\vdots&\vdots& &\vdots\\
*&\cdots&*&0&*&\cdots&*\\
\end{array}\right)}\}i_0-{\rm{th}}\rm{\:position}
\end{eqnarray*}
since $\ol{(A^*J_\mu A)}_{i_0 j}=-(A^*J_\mu A)_{j
i_0}=-i\sum \limits_{l=1}^{n-1}\mu_l\ol{a}_{l j}a_{l i_0}=0$ for
$j \in \{1,...,n\}$. As for each  $j\ne i_0$,
\begin{eqnarray*}
\underset{s\to\infty}{\lim}\frac{z_j(s)\ol{z}_{i_0}(s)}{\alpha_s}=0, \:\:\: \underset{s\to\iy}{\lim}\frac{z_j(s)}{\alpha_s}=0 \:\:\: \textrm{and
} \:\:\: \underset{s\to\infty}{\lim}\frac{\val{z_j(s)}^2}{\alpha_s}=0,
\end{eqnarray*}
one gets $ \underset{s\to\infty}{\lim}\la^s_j=\sum \limits_{l=1}^{n-1}\mu_l\val{a_{l
j}}^2$. On the other hand,
$\underset{s\to\iy}{\lim}\la^s_{i_0}+\frac{\val{z_{i_0}(s)}^2}{\alpha_s}=0$,
which in turn implies that
$\underset{s\to\infty}{\lim}\big|\la^s_{i_0}\big|=\iy$.
This proves that $i_0$ can only take the value 1 if $\alpha_s<0$ and $n$ if $\alpha_s>0$. Otherwise, since $\la^s_{i_0-1}\geq\la^s_{i_0}
\geq\la^s_{i_0+1}$, one gets $\underset{s\to\infty}{\lim}\la^s_{i_0-1}=\iy$ if $\alpha_s<0$ and $\underset{s\to\infty}{\lim}\la^s_{i_0+1}=-\iy$ if $\alpha_s>0$ which contradicts the fact that $\underset{s\to\infty}{\lim}\la^s_j$ is finite for all $j\ne i_0$.\\
~\\
Case $i_0=n$: \\
In this case, one has $\underset{s\to\infty}{\lim}\alpha_s\la^s_{n}=-\frac{r^2}{2}$ and
$\underset{s\to\infty}{\lim}\la^s_j=\sum \limits_{l=1}^{n-1}\mu_l\val{a_{l
j}}^2$ for all $j \in \{1,...,n-1\}$. Furthermore, the matrices $A$ and $A^*J_\mu
A$ have the form
\begin{eqnarray*}
A=\left( \begin{array}{cccccccc}
 & & &0\\
 &\tilde{A}& &\vdots\\
 & & &0\\
0&\cdots&0&e^{-i\theta}\\
\end{array}\right) \:\:\: \textrm{and} \:\:\: A^*J_\mu A=\left(
 \begin{array}{ccccc}
*&\ldots&*& 0\\
\vdots&&\vdots&\vdots\\
*&\ldots&*&0\\
0&\ldots&0&0\\
\end{array}\right),
\end{eqnarray*}
where $\tilde{A}\in U(n-1)$. However, the limit matrix of the subsequence $\big(J_{\la^s} +\frac{i}{\alpha_s}z(s)z(s)^* \big)_{s\in I}$ has to be diagonal because
$\underset{s\to\infty}{\lim}\frac{z_m(s)\ol{z}_{j}(s)}{\alpha_s}=0$ for all $m\ne j$. This implies that
$$A^*J_\mu A= \left(
 \begin{array}{ccccc}
i \mu_1 &\ldots&0& 0\\
\vdots&&\vdots&\vdots\\
0&\ldots& i \mu_{n-1}&0\\
0&\ldots&0&0\\
\end{array}\right)$$
and consequently, $\la_j^s=\mu_j$ for large $s$ and $j \in \{1,...,n-1\}$.\\
~\\
Case $i_0=1$: \\
In this case, $\underset{s\to\infty}{\lim}\alpha_s\la^s_{1}=-\frac{r^2}{2}$ and $\underset{s\to\infty}{\lim}\la^s_j=\sum \limits_{l=1}^{n-1}\mu_l\val{a_{l
j}}^2$ for every $j \in \{2,...,n\}$. Moreover, there is an element $\tilde{A}\in U(n-1)$ such
that the matrix $A$ is given by
\begin{eqnarray*}
A=\left( \begin{array}{cccccccc}
0& & & \\
\vdots& &\tilde{A}& \\
0& & & & \\
e^{-i\theta}&0&\cdots&0\\
\end{array}\right) \:\:\: \textrm{and hence,} \:\:\: A^*J_\mu A=\left( \begin{array}{cccccccc}
0&0&\cdots&0\\
0&*&\cdots &* \\
\vdots&\vdots& &\vdots\\
0&* & \cdots&* \\
\end{array}\right).
\end{eqnarray*}
Using the same arguments as above, one has $\la_{j+1}^s=\mu_j$ for $s$ large enough and for every $j \in \{1,...,n-1\}$.\\
~\\
Conversely, suppose that
$\underset{k\to\infty}{\lim}\alpha_k=0$. If the regarded sequence of orbits satisfies
the first condition, one can take $z(k):=\left( \begin{array}{c}
0\\
\vdots\\
0\\
\sqrt{-\alpha_k\la_n^k}\\
\end{array}\right)$ and  $A_k:=\I$ for $k\geq N$ and $N \in \N$ large enough. In the other case, one lets
\begin{eqnarray*}
z(k):=\left( \begin{array}{c}
\sqrt{-\alpha_k\la_1^k}\\
0\\
\vdots\\
\vdots\\
0\\
\end{array}\right) \:\:\: \textrm{and} \:\:\: A_k:=\left( \begin{array}{cccccccc}
0&1&0&\cdots&0\\
0&0&1&\ddots&0\\
\vdots&\vdots&\ddots&\ddots&\vdots\\
0&0&0&\ddots&1\\
1&0&0&\cdots&0\\
\end{array}\right) \textrm{ for } k\geq N.
\end{eqnarray*}
Thus, $\underset{k\to\infty}{\lim}\Big(A_k \big(J_{\la^k} +\frac{i}{\alpha_k} z(k)z(k)^* \big)A_k^*,\sqrt{2}A_k z(k),\alpha_k \Big)=(J_\mu,v_r,0)$.\\
~\\
3) Suppose that
$\big(\mathcal{O}_{(\la^k,\alpha_k)}\big)_{k\in\mathbb{N}}$ converges to the
orbit $\mathcal{O}_\la$. Then, there exist a sequence
$(A_k)_{k\in\mathbb{N}}$ in $U(n)$ and a sequence
$\big(z(k)\big)_{k\in\mathbb{N}}$ in $\C^n$ such that
\begin{equation*}
\underset{k\to\infty}{\lim}\bigg(A_k\Big(J_{\la^k} +\frac{i}{\alpha_k}
z(k)z(k)^*\Big)A_k^*,\sqrt{2}A_k z(k),\alpha_k \bigg)=(J_\la,0,0).
\end{equation*}
It follows that $\underset{k\to\infty}{\lim}\alpha_k=0$ and that
$\big(z(k)\big)_{k \in \N}$ tends to $0$ in $\C^n$. Denote by $A=(a_{mj})_{1\leq
m,j\leq n}$ the limit matrix of a subsequence $(A_s)_{s\in I}$
for an index set $I\subset \N$. Then,
$$\underset{s\to\infty}{\lim}J_{\la^s}+\frac{i}{\alpha_s} z(s)z(s)^*=A^*J_\la A \:\:\: \text{with} \:\:\: (A^*J_\la A)_{mj}=i\sum \limits_{l=1}^{n}\la_l\ol{a}_{lm}a_{lj}.$$
Since $\underset{k\to\infty}{\lim}\alpha_k=0$, one can assume that $\al_s$ is either strictly positive for all $s \in I$ or strictly negative for all $s \in I$.\\
Let $\sqrt{|\alpha_s|}$ be the square root of $|\alpha_s|$. The fact
that $\underset{s\to\infty}{\lim}\frac{z_m(s)\ol{z}_j(s)}{\alpha_s}$
is finite for all $m,j \in \{1,...,n\}$ implies that there exists at most one integer $1\leq
i_0\leq n$ such that
$\underset{s\to\infty}{\lim}\frac{z_{i_0}(s)}{\sqrt{|\alpha_s|}}=\iy$.
Therefore,
\begin{equation*}
\underset{s\to\infty}{\lim}
\frac{z_{j}(s)}{\sqrt{|\alpha_s|}}
\end{equation*}
exists for all $j$ distinct from $i_0$. Hence, for the same
reasons as in the proof of 4), necessarily $i_0\in\{1, n\}$.\\

If there is no such $i_0 $, then there exists for all $j \in \{1,...,n\}$ an integer $\la'_j \in \Z$ such that
$\la'_j=\la^s_{j}$ for all $s \in I$ (by passing to a subsequence if necessary) and
$\tilde{z}_j:=\underset{s\to\infty}{\lim}
\frac{z_{j}(s)}{\sqrt{|\alpha_s|}}$ is
finite for all $j \in \{1,\ldots,n\}$. Thus, one gets
\begin{eqnarray*}
\begin{cases}
A^*J_\la A=J_{\la'}+i \tilde{z}(\tilde{z})^*, \:\:\: ~\text{if}~\al_s>0~\forall s \in I \\
A^*J_\la A=J_{\la'}-i \tilde{z}(\tilde{z})^*, \:\:\: ~\text{if}~\al_s<0~\forall s \in I.
\end{cases}
\end{eqnarray*}
It follows by Lemma
\ref{unspctrela} applied to $\tilde{z}$ and $\alpha=1$ or $\alpha=-1$, respectively, that
\begin{eqnarray*}
\begin{cases}
\la_1\geq\la'_1=\la^s_1\geq\la_2\geq\la'_2=\la^s_2\geq...\geq\la_n\geq\la'_n=\la^s_n,\:\:\:~\text{if}~\al_s>0~\forall s \in I \\
\la'_1=\la^s_1\geq\la_1 \geq \la'_2=\la^s_2\geq\la_2\geq...\geq\la'_n=\la^s_n\geq \la_n,\:\:\:~\text{if}~\al_s<0~\forall s \in I.
\end{cases}
\end{eqnarray*}
~\\
Case $i_0=n$: \\
Here, $\underset{s\to\infty}{\lim}\alpha_s\la^s_{n}=0$, as $\underset{s\to\infty}{\lim}\big|\la^s_{n}+
\frac{\val{z_{n}(s)}^2}{\alpha_s}\big|<\iy$. Furthermore,
$\underset{s\to\infty}{\lim}\la^s_{n}=-\iy $ and
$\underset{s\to\infty}{\lim}\la^s_j=\sum \limits_{l=1}^{n}\la_l\val{a_{l
j}}^2$ for all $j \in \{1,...,n-1\}$ and $\al_s$ has to be positive for large $s$. Since $\lim \limits_{s\to\infty} \frac{z_j(s)}{\al_s}$ exists and
$\lim \limits_{s\to\iy}\frac{\val{z_n(s)}}{\al_s}=\iy $, it follows that
$\lim \limits_{s\to\infty} \frac{z_j(s)}{\al_s}=0$ for every $j \in \{1,..., n-1\} $. \\
Now, choose
\begin{equation*}
x:=\underset{s\to\infty}{\lim}\la^s_{n}+
\frac{\val{z_{n}(s)}^2}{\alpha_s}, \:\:\: \la'_j:=\underset{s\to\infty}{\lim}\la^s_{j} \:\:\: 
\textrm{and} \:\:\: w_j:=-i\underset{s\to\infty}{\lim}
\frac{z_j(s)\ol{z}_n(s)}{\alpha_s} \:\:\: \forall j \in \{1,...,n-1\}.
\end{equation*}
Then, the limit matrix $A^*J_\la
A$ of the sequence $\big(J_{\la^s} +\frac{i}{\alpha_s} z(s)z(s)^* \big)_{s \in I}$ has
the form
\begin{equation*}
\left( \begin{array}{cccccccc}
i\la'_1&0& \ldots &0&-w_1\\
0&i\la'_2& \ldots &0&-w_2\\
\vdots &\vdots&\ddots&\vdots&\vdots\\
0&0&\ldots&i\la'_{n-1}&-w_{n-1}\\
\ol{w}_1&\ol{w}_2&\ldots &\ol{w}_{n-1}&ix\\
\end{array}\right).
\end{equation*}
By Lemma \ref{unlemmatrix}, one obtains
$\la_1\geq\la'_1\geq\la_2\geq\la'_2\geq...\geq\la_{n-1}\geq\la'_{n-1}\geq\la_n$,
and therefore, one has
$\la_1\geq\la^s_1\geq\la_2\geq\la^s_2\geq...\geq\la_{n-1}\geq\la^s_{n-1}\geq\la_n\geq\la^s_{n}$
for large $s$.\\
~\\
Case $i_0=1$: \\
In this case, $\underset{s\to\infty}{\lim}\alpha_s\la^s_{1}=0$, since
$\underset{s\to\infty}{\lim}\big|\la^s_{1}+
\frac{\val{z_{1}(s)}^2}{\alpha_s}\big|<\iy$. Moreover,
$$\underset{s\to\infty}{\lim}\la^s_{1}=\iy , \:\:\: \underset{s\to\infty}{\lim}\la^s_j=\sum \limits_{l=1}^{n}\la_l\val{a_{l
j}}^2~ \:\:\: \text{and} \:\:\: \lim \limits_{s\to\infty} \frac{z_j(s)}{\al_s}=0 \:\:\: \forall j \in \{2,...,n\}.$$
Hence, $\alpha_s<0$ for $s$ large enough.
If one sets
\begin{equation*}
x:=\underset{s\to\infty}{\lim}\la^s_{1}+\frac{\val{z_{1}(s)}^2}{\alpha_s}, \:\:\: \la'_j:=\underset{s\to\infty}{\lim}\la^s_{j+1} \:\:\: 
\textrm{and} \:\:\: w_j:=-i\underset{s\to\infty}{\lim}
\frac{\ol{z}_1(s)z_{j+1}(s)}{\alpha_s} \:\:\: \forall j \in \{1,...,n-1\},
\end{equation*}
the limit matrix $A^*J_\la A$ of
$\big(J_{\la^s}+\frac{i}{\alpha_s} z(s)z(s)\big)_{s \in I}$ can be written as follows:
\begin{equation}
\left( \begin{array}{cccccccc}\label{unmatremm}
ix&\ol{w}_1&\ol{w}_2&\ldots &\ol{w}_{n-1}\\
-w_1&i\la'_1&0& \ldots &0\\
-w_2&0&i\la'_2& \ldots &0\\
\vdots&\vdots &\vdots&\ddots&\vdots\\
-w_{n-1}&0&0&\ldots&i\la'_{n-1}\\
\end{array}\right)=\tilde{A}^*\left( \begin{array}{cccccccc}
i\la'_1&0& \ldots &0&-w_1\\
0&i\la'_2& \ldots &0&-w_2\\
\vdots &\vdots&\ddots&\vdots&\vdots\\
0&0&\ldots&i\la'_{n-1}&-w_{n-1}\\
\ol{w}_1&\ol{w}_2&\ldots &\ol{w}_{n-1}&ix\\
\end{array}\right)\tilde{A},
\end{equation}
where
\begin{equation*} 
\tilde{A}=\left( \begin{array}{cccccccc}
0&1&0&\cdots&0\\
0&0&1&\ddots&0\\
\vdots&\vdots&\ddots&\ddots&\vdots\\
0&0&0&\ddots&1\\
1&0&0&\cdots&0\\
\end{array}\right).
\end{equation*}
This proves that
$\la^s_1\geq\la_1\geq\la^s_2\geq\la_2\geq...\geq\la^s_{n-1}\geq\la_{n-1}\geq\la^s_n\geq\la_{n}$
for large $s$.\\

Conversely, suppose that the sequence
$\big(\mathcal{O}_{(\la^k,\alpha_k)}\big)_{k\in\mathbb{N}}$ satisfies the
first condition. \\
First, consider the case $\la^k_n \overset{k \to \infty}{\longrightarrow} - \infty$. \\
Then, there is a subsequence $(\la^s)_{s\in
I}$ for an index set $I \subset \N$ fulfilling $\la^s_j=\la'_j$ for every $j \in \{1,..., n-1\}$
and all $s\in I$. By Lemma \ref{unlemmatrix}, there exist
$w_1,w_2,\ldots,w_{n-1} \in \C$, $x\in \R$ and $A\in U(n)$ such that
\begin{equation*}
A^*J_\la A=\left( \begin{array}{cccccccc}
i\la'_1&0& \ldots &0&-w_1\\
0&i\la'_2& \ldots &0&-w_2\\
\vdots &\vdots&\ddots&\vdots&\vdots\\
0&0&\ldots&i\la'_{n-1}&-w_{n-1}\\
\ol{w}_1&\ol{w}_2&\ldots &\ol{w}_{n-1}&ix\\
\end{array}\right).
\end{equation*}
In this case, $\la^k\ne\la$ for large $k$, as $\la^k_n \overset{k \to \infty}{\longrightarrow} - \infty$. Choose
$x:=\sum \limits_{j=1}^{n}\la_j-\sum \limits_{j=1}^{n-1}\la'_j$. (Compare the proof of
Lemma \ref{unlemmatrix}.) It follows that
\begin{eqnarray*}
\alpha_s(x-\la^s_n)=\sum_{j=1}^{n}\alpha_s(\la_j-\la^s_j)>0.
\end{eqnarray*}
Furthermore, define the sequence $\big(z(s)\big)_{s\in I}$ in $\C^n$ by
$$z_n(s):=\sqrt{\alpha_s(x-\la^s_n)} \:\:\: \text{and} \:\:\: z_j(s):=i\frac{\alpha_s w_j}{\sqrt{\alpha_s(x-\la^s_n)}} \:\:\: \forall j \in \{1,2,...,n-1\}.$$
Then, one gets
\begin{eqnarray*}
\underset{s\to\infty}{\lim}z(s)&=&0,\\[2pt]
\la^s_n+\frac{\val{z_n(s)}^2}{\alpha_s}&=&x,\\[2pt]
\underset{s\to\infty}{\lim}\frac{\val{z_j(s)}^2}{\alpha_s}&=&\underset{s\to\infty}{\lim}\frac{\val{w_j}^2}{x-\la^s_n}=0 \:\:\: 
\forall j \in \{1,...,n-1\},\\[2pt]
\underset{s\to\infty}{\lim}\frac{z_m(s)\ol{z_j(s)}}{\alpha_s}&=&\underset{s\to\infty}{\lim}\frac{
w_m\ol{w}_j}{x-\la^s_n}=0 \:\:\:  \forall m \ne j \in \{1,...,n-1 \} \:\:\: \text{and}\\[2pt]
\underset{s\to\infty}{\lim}\frac{z_j(s)\ol{z_n(s)}}{\alpha_s}&=&iw_j  \:\:\: \forall j \in \{1,...,n-1\}.
\end{eqnarray*}
Hence, $\Big(A \big(J_{\la^s}+\frac{i}{\alpha_s}z(s)z(s)^*\big)A^*\Big)_{s\in I}$ converges to $J_\la$ and $\big(z(s)\big)_{s \in I}$ to $0$.\\

If $\lim \limits_{k \to \infty} \la^k_n \not= - \infty$, there is a subsequence $(\la^s)_{s \in I}$ for an index set $I \subset \N$ fulfilling $\la_j^s=\la'_j$ for all
$j \in \{1,...,n\}$ and all $s \in I$. Therefore,
$$\la_1 \geq \la'_1 \geq \la_2 \geq \la'_2 \geq...\geq \la_n \geq \la'_n$$
and thus, by Lemma \ref{unspctrela}(2), there exists $\tilde{z} \in \C^n$ such that $i \la_1,...,i \la_n$ are the eigenvalues of $J_{\la'}+i \tilde{z}(\tilde{z})^*$. \\
Let $z(s):=\tilde{z} \sqrt{\al_s}$, which is reasonable since $\al_s>0$ in this case. \\
As the matrices $J_{\la'}+i \tilde{z}(\tilde{z})^*$ and $J_{\la}$ are both skew-Hermitian and have the same eigenvalues, they are unitarily conjugated. Therefore, there exists an element $A \in U(n)$ such that $J_{\la'}+i \tilde{z}(\tilde{z})^*=A^*J_{\la}A$. Hence,
$$A^*J_{\la}A=J_{\la'}+i \tilde{z}(\tilde{z})^*=\lim \limits_{s \to \infty} J_{\la'}+i \tilde{z}(\tilde{z})^*
=\lim \limits_{s \to \infty} J_{\la'}+i \frac{z(s)z(s)^*}{\al_s},$$
i.e. $\Big(A \big(J_{\la'}+i \frac{z(s)z(s)^*}{\al_s} \big)A^*\Big)_{s \in I}$ converges to $J_{\la}$. Furthermore,
$$z(s)=\tilde{z} \sqrt{\al_s} \overset{s \to \infty}{\longrightarrow} 0,$$
as $\al_k \overset{k \to \infty}{\longrightarrow} 0$. \\
Thus, the claim is shown in this case. \\

Suppose now that for $k$ large $\alpha_k<0$,
$\la_1^k\geq\la_1\geq...\geq\la^k_{n-1}\geq\la_{n-1}\geq\la_n^k\geq\la_n$ and $\underset{k\to\infty}{\lim}\alpha_k\la^k_{1}=0$. First consider the case $\la_1^k \overset{k \to \infty}{\longrightarrow} \infty$. \\
In this case,
there is a subsequence $(\la^s)_{s\in I}$ for an index set $I \subset \N$ such that
$\la^s_j=\la'_{j-1}$ for all  $j \in \{2,...,n\}$ and all $s\in I$. By
Identity (\ref{unmatremm}) and Lemma \ref{unlemmatrix} ,
there exist $w_1,w_2,\ldots,w_{n-1} \in \C$, $x\in \R$ and $A\in
U(n)$ such that
\begin{equation*}
A^*J_\la A=\left( \begin{array}{cccccccc}\label{unmatrem}
ix&\ol{w}_1&\ol{w}_2&\ldots &\ol{w}_{n-1}\\
-w_1&i\la'_1&0& \ldots &0\\
-w_2&0&i\la'_2& \ldots &0\\
\vdots&\vdots &\vdots&\ddots&\vdots\\
-w_{n-1}&0&0&\ldots&i\la'_{n-1}\\
\end{array}\right).
\end{equation*}
Similarly to the last case, one takes
$x:=\sum \limits_{j=1}^{n}\la_j-\sum \limits_{j=1}^{n-1}\la'_j$ and thus gets
$$\alpha_s(x-\la^s_1)=\sum_{j=1}^{n}\alpha_s(\la_j-\la_j^s)>0.$$
Hence, one can define the sequence $\big(z(s)\big)_{s\in I}$ in $\C^n$ by
$$z_1(s):=\sqrt{\alpha_s(x-\la^s_1)} \:\:\: \text{and} \:\:\:  z_{j}(s):=-i\frac{\alpha_s
w_{j-1}}{\sqrt{\alpha_s(x-\la^s_1)}} \:\:\: \forall j \in \{2,...,n\}.$$
Here again, one gets
\begin{eqnarray*}
\underset{s\to\infty}{\lim}z(s)&=&0,\\[2pt]
\la^s_1+\frac{\val{z_1(s)}^2}{\alpha_s}&=&x,\\[2pt]
\underset{s\to\infty}{\lim}\frac{\val{z_j(s)}^2}{\alpha_s}&=&\underset{s\to\infty}{\lim}\frac{\val{w_{j-1}}^2}{x-\la^s_1}=0 \:\:\: \forall j \in \{2,...,n\},\\[2pt]
\underset{s\to\infty}{\lim}\frac{z_m(s)\ol{z_j(s)}}{\alpha_s}&=&\underset{s\to\infty}{\lim}\frac{ w_{m-1}\ol{w}_{j-1}}{x-\la^s_1}=0 \:\:\: 
\forall m \ne j \in \{2,...,n\} \:\:\: \text{and}\\[2pt]
\underset{s\to\infty}{\lim}\frac{z_j(s)\ol{z_1(s)}}{\alpha_s}&=&iw_{j-1}
 \:\:\:  \forall j \in \{2,...,n\}.
\end{eqnarray*}
Again, one can conclude that
$\bigg(\Big(A \big(J_{\la^s}+\frac{i}{\alpha_s}z(s)z(s)^*\big)A^*,\sqrt{2}Az(s),\al_s \Big)\bigg)_{s\in
I}$ converges to $(J_\la,0,0)$. \\

If $\lim \limits_{k \to \infty} \la^k_1 \not= \infty$, there is a subsequence $(\la^s)_{s \in I}$ for an index set $I \subset \N$ fulfilling $\la_j^s=\la'_j$ for all
$j \in \{1,...,n\}$ and all $s \in I$. Hence,
$$\la'_1 \geq \la_1 \geq \la'_2 \geq \la_2 \geq...\geq \la'_n \geq \la_n$$
and therefore, by Lemma \ref{unspctrela}(2), there exists $\tilde{z} \in \C^n$ in such a way that $i \la_1,...,i \la_n$ are the eigenvalues of $J_{\la'}-i \tilde{z}(\tilde{z})^*$. \\
Let now $z(s):=\tilde{z} \sqrt{-\al_s}$, which is reasonable since this time $\al_s<0$. \\
As above, there exists an element $A \in U(n)$ such that $J_{\la'}-i \tilde{z}(\tilde{z})^*=A^*J_{\la}A$ and thus,
$$A^*J_{\la}A=\lim \limits_{s \to \infty} J_{\la'}-i \frac{z(s)z(s)^*}{-\al_s}=\lim \limits_{s \to \infty} J_{\la'}+i \frac{z(s)z(s)^*}{\al_s},$$
i.e. $\Big(A \big(J_{\la'}+i \frac{z(s)z(s)^*}{\al_s} \big)A^* \Big)_{s \in I}$ converges to $J_{\la}$. Furthermore,
$$z(s)=\tilde{z} \sqrt{-\al_s} \overset{s \to \infty}{\longrightarrow} 0,$$
as $\al_k \overset{k \to \infty}{\longrightarrow} 0$. \\
Therefore, the assertion is also shown in this case.
~\\
\qed

\section{The continuity of the inverse of the Kirillov-Lipsman map $\K $}\label{top-spectrum}
In the  next two  sections, the topology of the spectrum  $\widehat{G_n} $ of the group $G_n=U(n)\ltimes\HH_n$ will be analyzed and the aim is to show that it is determined by the topology of its admissible quotient space. \\

\subsection{The representation $\pi_{(\mu,r)}$}
~\\

First, examine the representation
$\pi_{(\mu,r)}=\ind_{U(n-1)\ltimes\HH_n}^{G_n} \rho_\mu\otimes\chi_r$. Its Hilbert space $\H_{(\mu,r)} $ is given by the space
$$L^2 \Big(G_n/\big(U(n-1)\ltimes \HH_n \big),
\rho_\mu\otimes\chi_r \Big)\cong L^2 \big(U(n)/U(n-1),\rh_\mu \big).$$ Let $\xi$
be a unit vector in
 $\H_{(\mu,r)}$ and recall that $(z,w)_{\C^n}=\Re\langle z,w\rangle_{\C^n} $ for $z,w\in\C^n $. For
all $(z,t)\in\HH_n$ and all $A,B\in U(n)$,
\begin{eqnarray*} 
\pi_{(\mu,r)}(A,z,t)\xi (B)=e^{-i(Bv_r,z)_{\C^n}}\xi \big(A^{-1}B \big).
\end{eqnarray*}
 Therefore,
\begin{eqnarray*} 
C^{\pi_{(\mu,r)}}_{\xi}(A,z,t)&=&\Big\langle\pi_{(\mu,r)}(A,z,t)\xi,\xi \Big\rangle_{L^2 \big(U(n)/U(n-1),\rh_\mu \big)}\nn\\[2pt]
&=&\int \limits_{U(n)}e^{-i(Bv_r,z)_{\C^n}} \Big\langle\xi \big(A^{-1}B \big),\xi(B) \Big\rangle_{\mathcal{H}_{\rho_\mu}}dB.
\end{eqnarray*}
By (\ref{lamu}) in Section \ref{duun},
one has
\begin{eqnarray*} 
\pi_\mu:=\pi_{(\mu,r)}{\res{U(n)}}\cong
\ind_{U(n-1)}^{U(n)}\rho_\mu=\sum_{\substack{\tau_{\lambda}\in\widehat{U(n)} \\ \lambda_1\geq
 \mu_1\geq \lambda_2\geq \mu_2\geq...\geq\lambda_{n-1}\geq \mu_{n-1}\geq
 \lambda_n}}
\tau_{\lambda}.
\end{eqnarray*}
Every irreducible representation $\ta_\la $ of $U(n) $ can be realized as a subrepresentation of the left regular representation on $ L^2 \big(U(n) \big) $ via the
intertwining operator $$ U_\la: \H_\la\to L^2 \big(U(n) \big),~
U_\la(\xi)(A):=\big\langle{\xi},{\ta_\la(A)\xi_\la}\big\rangle_{\H_{\la}} \:\:\: \forall A\in U(n)~\forall \xi \in \H_{\la}$$
for a fixed unit vector $\xi_{\la} \in \H_{\la}$. \\

For $ \ta_\la \in \widehat {U(n)} $, consider the
orthonormal basis $ \B^\la =\big\{\ph^\la_j|~ j \in \{1, ... ,d_\la\} \big\}$ of $
\H_\la $ consisting of eigenvectors for $ \T_n $ of $\H_\la$.
%
~\\
Moreover, as a basis of the Lie algebra $\h_n$ of the Heisenberg group, one can take the left invariant vector fields $\big\{Z_1,Z_2,\ldots,Z_n$,
 $\ol{Z}_1,\ol{Z}_2,\ldots,\ol{Z}_n,T \big\}$, where
\begin{equation*}
Z_j:=2\frac{\partial}{\partial\ol{z}_j}+i\frac{z_j}{2}\frac{\partial}{\partial t}, \:\:\: \ol{Z}_j=2\frac{\partial}{\partial z_j}-i\frac{\ol{z}_j}{2}\frac{\partial}{\partial t} \:\:\: \text{and} \:\:\: T:=\frac{\partial}{\partial t}
\end{equation*}
and gets the Lie brackets $[Z_j,\ol{Z}_j]=-2iT$ for $j \in \{1,...,n\}$. \\
Now, regard the Heisenberg sub-Laplacian differential operator which is given by
\begin{equation*}
\LL=\frac{1}{2}\sum_{j=1}^n \big(Z_j\ol{Z}_j+\ol{Z}_j Z_j \big).
\end{equation*} 
This operator is $U(n)$-invariant.

\begin{lemma}\label{laplaceirr}
~\\
For every representation $ \pi_{(\mu,r)}$ for $r>0$ and
$\rh_\mu\in \widehat{U(n-1)} $,
\begin{equation}\label{}
\nonumber d\pi_{(\mu,r)}(\LL)=-r^ 2 \Id.
\end{equation}
\end{lemma}

Proof: \\
Since the representation $\pi_{(\mu,r)}$ is trivial on the center of $\h_n$, one has
\begin{eqnarray*}
d\pi_{(\mu,r)}(\LL)\xi(B)=2\sum_{j=1}^n \bigg(\frac{\partial^2}{\partial
z_j\partial\ol{z}_j}+\frac{\partial^2}{\partial\ol{z}_j\partial
z_j} \bigg)\Big(e^{-i( Bv_r,z)_{\C^n}}\Big)\xi(B).
\end{eqnarray*}
Let $\D=\{e_1,\ldots,e_n\}$ be an orthonormal basis
for $\C^n$. By writing $$(Bv_r,z)_{\C^n}=\frac{1}{2} \Big(\lan Bv_r,z\ran_{\C^n} +
\ol{\lan Bv_r,z\ran}_{\C^n}\Big),$$
one gets
\begin{eqnarray*}
d\pi_{(\mu,r)}(\LL)\xi(B)=-\sum_{j=1}^n\big|\lan
Bv_r,e_j\ran_{\C^n}\big|^2\xi(B)=-r^2\xi(B).
\end{eqnarray*}
\qed
~\\
~\\
In addition, the following theorem describes the convergence of sequences of representations $\big(\pi_{(\mu^k,r_k)} \big)_{k \in \N}$:

\begin{theorem}\label{rnenull}
~\\
Let $r>0$, $\rho_\mu\in\widehat{U(n-1)}$ and
$\ta_\la\in\widehat{U(n)}$.
\begin{enumerate}
\item A sequence $\big(\pi_{(\mu^k,r_k)} \big)_{k \in \N}$ of irreducible unitary representations
of $ G_n $ converges to $\pi_{(\mu,r)}$ in $ \wh{G_n} $ if and only
if $\underset{k\to\infty}{\lim}r_k=r$ and $\mu^k=\mu$ for  $k$ large enough.
\item A sequence $\big(\pi_{(\mu^k,r_k)} \big)_{k \in \N}$ of irreducible unitary representations
of $ G_n $ converges to $\ta_\la$ in $ \wh{G_n} $ if and only if
$\underset{k\to\infty}{\lim}r_k=0$ and $\ta_\la$ occurs in
$\pi_{\mu^k}$ for $k$ large enough.
\end{enumerate}
These are all possibilities for a sequence $\big(\pi_{(\mu^k,r_k)} \big)_{k \in \N}$ of irreducible unitary representations
of $ G_n $ to converge.
\end{theorem}
The proof of 1) and 2) of this theorem can be found in [\ref{Ba}], Theorem 6.2.A. \\
Furthermore, since the representations $\pi_{(\mu,r)}$ and $\ta_\la$ are trivial on $\big\{(\Id,0,t)|~t \in \R\big\}$, the center of $G_n$, while the representations $\pi_{(\la,\al)}$ are non-trivial there, the possibilities of convergence of a sequence $\big(\pi_{(\mu^k,r_k)} \big)_{k \in \N}$ listed above are the only ones that are possible.

\subsection{The representation $\tau_{\la}$}\label{tau-la}
~\\

As $\tau_{\la}$ only acts on $U(n)$ and $\widehat{U(n)}$ is discrete, every converging sequence $(\tau_{\la^k})_{k \in \N}$ has to be constant for large $k$. Hence,
$$\tau_{\la^k} \overset{k \to \iy}{\longrightarrow} \tau_{\la} ~\Longleftrightarrow~ \la^k=\la ~~ \text{for large}~k.$$
~\\

\subsection{The representation $\pi_{(\la,\alpha)}$}
~\\

Next, regard the representations $\pi_{(\la,\alpha)}$. \\
Consider the unit vector $\xi:=\sum \limits_{j=1}^{d_\la}
\phi_j^\la\otimes f_j$ in the Hilbert space
$\H_{(\la,\alpha)}=\H_\la\otimes \F_\alpha(n)$ of
$\pi_{(\la,\alpha)}$, where $f_1,\ldots, f_{d_\la}$ belong to the
Fock space $\F_\alpha(n)$. Then, for all $A\in U(n)$ and $(z,t)\in \HH_n$,
\begin{eqnarray*} 
\pi_{(\la,\alpha)}(A,z,t)\xi(w)&=&\sum \limits_{j=1}^{d_\la}
\ta_\la(A)\phi_j^\la\otimes e^{i\alpha
t-\frac{\alpha}{4}|z|^2-\frac{\alpha}{2}\lan w,z\ran_{\C^n}} f_j \big(A\inv
w+A\inv z \big) \:\:\: \textrm {if}~ \alpha>0 \:\:\: \text{and} \\
\pi_{(\la,\alpha)}(A,z,t)\xi(\ol{w})&=&\sum \limits_{j=1}^{d_\la} \ta_\la(A)\phi_j^\la\otimes e^{i\alpha
t+\frac{\alpha}{4}|z|^2+\frac{\alpha}{2}\lan\ol{w},\ol{z}\ran_{\C^n}} f_j \Big(\ol{A\inv
w}+\ol{A\inv z} \Big) \:\:\: \textrm{if}~ \alpha<0.
\end{eqnarray*}
It follows that
\begin{eqnarray*} 
&&C^{\pi_{(\la,\alpha)}}_{\xi}(A,z,t)~=~\big\lan\pi_{(\la,\alpha)}(A,z,t)\xi,\xi \big\ran_{\H_{(\la,\alpha)}} \\
=&&\begin{cases}
\sum \limits_{j,j'=1}^{d_\la} \big\lan\ta_\la(A)\phi_j^\la,\phi_{j'}^\la \big\ran_{\H_{\la}} \int \limits_{\C^n} e^{i\alpha t-\frac{\alpha}{4}|z|^2-
\frac{\alpha}{2}\lan w,z\ran_{\C^n}}f_j \big(A\inv w+A\inv z \big)\ol{f_{j'}(w)}e^{-\frac{\alpha}{2}\vert w\vert^2}dw \:\:\: \textrm {if}~ \alpha>0,\\[7pt]
\sum \limits_{j,j'=1}^{d_\la}
\big\lan\ta_\la(A)\phi_j^\la,\phi_{j'}^\la \big\ran_{\H_{\la}} \int \limits_{\C^n} e^{i\alpha
t+\frac{\alpha}{4}|z|^2+\frac{\alpha}{2}\lan\ol{w},\ol{z}\ran_{\C^n}}f_j \Big(\ol{A\inv
w}+\ol{A\inv z} \Big)\ol{f_{j'}(\ol{w})}e^{\frac{\alpha}{2}\vert
w\vert^2}dw \:\:\: \textrm {if}~ \alpha<0.\nn 
\end{cases} \\
\end{eqnarray*}

\begin{lemma}\label{uncenter}
~\\
For each representation $\pi_{(\la,\alpha)}$ for
$\alpha\in\R^*$ and $\ta_\la\in\widehat{U(n)}$, one has
\begin{eqnarray*}
d\pi_{(\la,\alpha)}(T)=i\alpha \I.
\end{eqnarray*}
\end{lemma}

Proof: \\
Let $\xi=\sum\limits_{j=1}^{d_\la} \phi_j^\la\otimes f_j$ be a unit vector in $\H_{(\la,\alpha)}$. Then,
\begin{eqnarray*}
\big\lan d\pi_{(\la,\alpha)}(T)\xi, \xi \big\ran_{\H_{(\la,\alpha)}}~=~\frac{d}{dt}\Big|_{t=0}\big\lan
\pi_{(\la,\alpha)}(\I,0,t)\xi,
\xi \big\ran_{\H_{(\la,\alpha)}}~=~\frac{d}{dt}\Big|_{t=0}e^{i\alpha t}~ \sum \limits_{j=1}^{d_\la} \Vert f_j\Vert_{\F_{\al}(n)}^2~=~i\alpha.
\end{eqnarray*}
\qed
~\\
~\\
If $\alpha$ is positive, the polynomials $\C[\C^n]$ are
dense in $\F_\alpha(n)$ and its multiplicity free decomposition is
\begin{equation*}
\C[\C^n]=\sum_{m=0}^{\infty}\P_m,
\end{equation*}
where $\P_m$ is the space of homogeneous polynomials of degree $m$.
Therefore, $p_m(z)=z_1^m$ is the highest weight vector in $\P_m$ with weight $(m,0,...,0)=:[m]$. Applying the classical Pieri's
rule (see [\ref{Fu-Ha}], Proposition 15.25), one obtains
\begin{eqnarray}\label{decoppositivealp}
\big(\ta_{\la}\otimes W_{\al}\big) \res{U(n)}=\sum_{m=0}^{\infty}\ta_\la\otimes\ta_{[m]}=\sum \limits_{\substack{\la'\in P_n \\ \la'_{1}\geq\la_{1}\geq....\geq\la'_{n}\geq\la_{n}}} \ta_{\la'},
\end{eqnarray}
where the definition of the operator $W_{\al}(A)$ can be found in Section 3 in the description of $\pi_{(\lambda, \alpha)}$. \\
If $\alpha$ is negative, one gets
\begin{eqnarray*}
\big(\ta_{\la}\otimes W_{\al}\big) \res{U(n)}=\sum_{m=0}^{\infty}\ta_\la\otimes\ta_{[m]}=\sum \limits_{\substack{\la'\in P_n \\ \la_{1}\geq\la'_{1}\geq....\geq\la_{n}\geq\la'_{n}}} \ta_{\la'}.
\end{eqnarray*}
Both of the sums again are multiplicity free. This follows from [\ref{knapp}], Chapter IV.11, since $W_{\al}$ is multiplicity free.\\

Furthermore, let
$\RR_\al:=\big\{h_{m,\alpha}|~\:m=(m_1,\ldots,m_n)\in\N^n \big\}$ be the
orthonormal basis of the Fock space $\F_\alpha(n)$ defined by the
Hermite functions
\begin{equation*}
h_{m,\alpha}(z)=\bigg(\frac{\val{\alpha}}{2\pi}\bigg)^{\frac{n}{2}}\sqrt{\frac{\val{\alpha}^{\val{m}}}{2^{\val{m}}m!}}~z^m
\end{equation*}
with $\val{m}=m_1+...+m_n$, $m!=m_1!\cdots m_n!$ and
$z^m=z_1^{m_1}\cdots z_n^{m_n}$ (see [\ref{Folland}], Chapter 1.7). \\
~\\
Now, one obtains the following theorem about the convergence of sequences of representations $\big(\pi_{(\la^k,\alpha_k)} \big)_{k \in \N}$:

\begin{theorem}\label{conv-la-al}
~\\
Let $\alpha\in\R^*$ and $\ta_\la\in\widehat{U(n)}$. Then, a sequence
$\big(\pi_{(\la^k,\alpha_k)} \big)_{k \in \N}$ of elements in $\widehat{G_n}$
converges to $\pi_{(\la,\alpha)}$ if and only if
$\underset{k\to\iy}{\lim}\alpha_k=\alpha$ and $\la^k=\la$ for large
$k$.
\end{theorem}

Proof: \\
First, consider the case where $\alpha$ is positive. Assume that $\alpha_k \overset{k \to \iy}{\longrightarrow}\alpha$ and that $\la^k=\la$ for $k$
large enough. Moreover, let $f\in C_0^\iy(G_n)$ and let $\xi$ be a unit vector in $\H_\la$. Then,
\begin{eqnarray*}
&& \Big\lan C^{\pi_{(\la^k,\alpha_k)}}_{\xi\otimes
h_{0,\alpha_k}},f \Big\ran_{\big(L^{\iy}(G_n),L^1(G_n)\big)} \\[2pt]
&=&\int \limits_{U(n)}\int \limits_{\HH_n}f(A,z,t) \big\lan\ta_{\la^k}(A)
\xi,\xi \big\ran_{\H_{\la^k}} e^{i\alpha_k
t-\frac{\alpha_k}{4}|z|^2}\int \limits_{\C^n}\Big(\frac{1}{2\pi}\Big)^n
e^{-\frac{1}{2}\sqrt{\alpha_k}\lan w,z\ran_{\C^n}-\frac{1}{2}|w|^2}dw d(z,t)dA
\end{eqnarray*}
tends to $\Big\lan C^{\pi_{(\la,\alpha)}}_{\xi\otimes h_{0,\alpha}},f \Big\ran_{\big(L^{\iy}(G_n),L^1(G_n)\big)}$. Hence, $\big(\pi_{(\la^k,\alpha_k)}\big)_{k \in \N}$ converges to $\pi_{(\la,\alpha)}$. \\
The same reasoning applies when $\al$ is negative.\\

Conversely, the fact that the sequence $\big(\pi_{(\la^k,\alpha_k)}\big)_{k \in \N}$
converges to the representation $\pi_{(\la,\alpha)}$ implies
that for $\xi\in \H^\iy_{(\la,\alpha)}$ of length 1, there is for every $k \in \N$ a unit vector
 $\xi_k\in\H^\iy_{(\la^k,\alpha_k)}$ such that $\Big(\big\lan
d\pi_{(\la^k,\alpha_k)}(T)\xi_k,\xi_k \big\ran_{\H_{(\la^k,\alpha_k)}} \Big)_{k \in \N}$ converges to $\big\lan
d\pi_{(\la,\alpha)}(T)\xi,\xi \big\ran_{\H_{(\la,\alpha)}}$. Thus, by Lemma \ref{uncenter}, we have
$\underset{k\to\iy}{\lim}\alpha_k=\alpha$. Hence, it remains to show that $\la^k=\la$ for $k$ large enough. \\
Let $\xi$ be a unit vector in $\H_{\la}$. Then for every $k \in \N$, there exists a vector $\xi_k=\underset{m\in\N^n}{\sum}\zeta^k_m\otimes
h_{m,\al_k}\in\mathcal{H}_{(\la^k,\alpha_k)}$ of length 1 such that
$\Big(C^{\pi_{(\la^k,\alpha_k)}}_{\xi_k} \Big)_{k \in \N}$ converges uniformly on
compacta to $C^{\pi_{(\la,\alpha)}}_{\xi\otimes h_{0,\alpha}}$.

Now, take $\delta \in \R_{>0}$ such that $0\not\in
I_{\alpha,\delta}=(\alpha-\delta,\alpha+\delta)$, as well as a
Schwartz function $\varphi$ on $\R$ fulfilling
$\varphi\res{I_{\alpha,\delta}}\equiv1$ and $\varphi\equiv0$ in a
neighbourhood of $0$. Then, there is a
Schwartz function $\psi$ on $\HH_n$ with the property
\begin{equation*}
\sigma_\beta(\psi)=\varphi(\beta)P_\beta~~~\: \: \forall \be \in\R^*,
\end{equation*}
where $\sigma_\beta$ is the $\mathbb{H}_n$-representation defined in Section 3 and $P_\beta :\F_\beta(n)\to \C$ is the orthogonal
projection onto the one-dimensional subspace $\C h_{0,\beta}$ of all
constant functions in $\F_\beta(n)$. On the other hand, there exists
$k_\delta\in\N$ such that $\alpha_k\in I_{\alpha,\delta}$ for all
$k\geq k_\delta$. One obtains
$\sigma_\alpha(\psi)h_{0,\alpha}=h_{0,\alpha}$ and $\sigma_{\alpha_k}(\psi)h_{0,\alpha_k}=h_{0,\alpha_k}$ for all $k\geq
k_\delta$ and thus, it follows that
\begin{eqnarray*}
\underset{k\rightarrow\iy}{\lim}\big\|\zeta^k_0\big\|_{\H_{\la^k}}^2
&=&\underset{k\rightarrow\iy}{\lim}\underset{m,m'\in\N^n}{\sum}\big\lan\zeta^k_m,
\zeta^k_{m'}\big\ran_{\H_{\la^k}} \big\lan\sigma_{\alpha_k}(\psi)h_{m,\al_k},h_{m',\al_k}\big\ran_{\F_{\al_k}(n)}\\[4pt]
&=&\underset{k\rightarrow\iy}{\lim}\bigg\lan C^{\pi_{(\la^k,\alpha_k)}}_{\underset{m\in\N^n}{\sum}\zeta^k_m\otimes h_{m,\al_k}}(\I,.,.),\overline{\psi}\bigg\ran_{\big(L^{\iy}(\HH_n),L^1(\HH_n)\big)}\\[4pt]
&=&\big\lan\sigma_\alpha(\psi)h_{0,\alpha},h_{0,\alpha}\big\ran_{\F_{\al}(n)}=1.
\end{eqnarray*}
Hence, one gets
\begin{equation*}
\underset{k\rightarrow\iy}{\lim}\big\|\xi_k-\zeta^k_0\otimes
h_{0,\alpha_k}\big\|_{\mathcal{H}_{(\la^k,\alpha_k)}}=0
\end{equation*}
and one can deduce that
\begin{eqnarray*} 
\underset{k\to\iy}{\lim}\big\lan\ta_{\la^k}(A)\zeta^k_0,\zeta^k_0\big\ran_{\H_{\la^k}}=\big\lan
\ta_\la(A)\xi,\xi\big\ran_{\H_{\la}}
\end{eqnarray*}
uniformly in $A\in U(n)$. Therefore, for all $k\in\N$, one can take the unit vector $\ph_k=\frac{\zeta^k_0}{\Vert\zeta^k_0\Vert_{\H_{\la^k}}}$ in
$\H_{\la^k}$ to finally obtain the uniform convergence on
compacta of $\Big(C^{\ta_{\la^k}}_{\ph_k}\Big)_{k \in \N}$ to $C^{\ta_\la}_\xi$. Thus, $\la^k=\la$ for $k$ large enough.
~\\
\qed

\begin{lemma}\label{unlap2}
~\\
For each representation $\pi_{(\la,\alpha)}$ for
$\alpha\in\R^*$ and $\ta_\la\in\widehat{U(n)}$,
\begin{eqnarray*}
\big\lan
d\pi_{(\la,\alpha)}(\LL)h_{m,\al},h_{m,\al}\big\ran_{\F_{\al}(n)}=-\vert\al\vert \big(n
+2\vert m\vert \big) \:\:\: \forall m\in\N^n.
\end{eqnarray*}
\end{lemma}
The proof follows from [\ref{BJR}], Proposition 3.20 together with [\ref{BJRW}], Lemma 3.4.

\begin{theorem}\label{unconrep}
~\\
Let $r>0$, $\rho_\mu\in\widehat{U(n-1)}$ and $\ta_\la\in\widehat{U(n)}$.
\begin{enumerate}
\item If a sequence $\big(\pi_{(\la^k,\alpha_k)}\big)_{k\in\mathbb{N}}$ of
elements of $\wh{G_n}$ converges to the representation
$\pi_{(\mu,r)}$  in $\wh{G_n}$, then
$\underset{k\to\infty}{\lim}\alpha_k=0$ and the sequence
$\big(\pi_{(\la^k,\alpha_k)}\big)_{k\in\mathbb{N}}$ satisfies one of the
following conditions:
\begin{enumerate}[(i)]
\item For $k$ large enough, $\alpha_k>0$, $\la_j^k=\mu_j$
for all $j \in \{1,..., n-1\}$ and
$\underset{k\to\infty}{\lim}\alpha_k\la^k_{n}=-\frac{r^2}{2}$.
\item For $k$ large enough, $\alpha_k<0$,
$\la_j^k=\mu_{j-1}$ for all $j \in \{2,..., n\}$ and
$\underset{k\to\infty}{\lim}\alpha_k\la^k_{1}=-\frac{r^2}{2}$.
\end{enumerate}
\item If a sequence $\big(\pi_{(\la^k,\alpha_k)}\big)_{k\in\mathbb{N}}$ of
elements of $\wh{G_n}$ converges to the representation $\ta_\la$ in
$\wh{G_n}$, then $\underset{k\to\infty}{\lim}\alpha_k=0$ and the
sequence $\big(\pi_{(\la^k,\alpha_k)}\big)_{k\in\mathbb{N}}$ satisfies one
of the following conditions:
\begin{enumerate}[(i)]
\item For $k$ large enough, $\alpha_k>0$,
$\la_1\geq\la_1^k\geq...\geq\la_{n-1}\geq\la_{n-1}^k\geq\la_n\geq\la_n^k$ and $\underset{k\to\infty}{\lim}\alpha_k\la^k_{n}=0$.
\item For $k$ large enough, $\alpha_k<0$,
$\la_1^k\geq\la_1\geq\la_2^k\geq\la_2\geq...\geq\la_{n-1}\geq\la_{n}^k\geq\la_n$ and $\underset{k\to\infty}{\lim}\alpha_k\la^k_{1}=0$. \\
\end{enumerate}
\end{enumerate}
\end{theorem}

Proof: \\
1) Let $\tilde{\mu}^s=(\mu_1,\ldots,\mu_s,\mu_s,\mu_{s+1},\ldots,\mu_{n-1})$ for $s \in \{1,...,n-1\}$. By
hypothesis, the sequence $\big(\pi_{(\la^k,\al_k)}\big)_{k \in \N}$ converges to the
representation $\pi_{(\mu,r)}$ in $\widehat{G_n}$. Thus, for the unit vector $\xi^s=\sqrt{d_{\tilde{\mu}^s}}\ol{C^{\tilde{\mu}^s}_{\phi_1^{\tilde{\mu}^s},\phi_1^{\tilde{\mu}^s}}}\in
\H_{(\mu,r)}^\iy$, there is a sequence of unit vectors
$(\xi_k^s )_{k \in \N} \subset \H_{(\la^k,\al_k)}^\iy$ such that
\begin{eqnarray*}
\big\lan d\pi_{(\la^k,\al_k)}(T)\xi_k^s,\xi_k^s \big\ran_{\H_{(\la^k,\al_k)}} \overset{k \to \iy}{\longrightarrow}\big\lan
d\pi_{(\mu,r)}(T)(\xi^s),\xi^s \big\ran_{\H_{(\mu,r)}}=0 \:\:\: \forall T\in\t_n \:\:\: \text{and}
\end{eqnarray*}
\begin{eqnarray*}
\big\lan d\pi_{(\la^k,\al_k)}(\LL)\xi_k^s,\xi_k^s \big\ran_{\H_{(\la^k,\al_k)}} \overset{k \to \iy}{\longrightarrow} \big\lan
d\pi_{(\mu,r)}(\LL)(\xi^s),\xi^s \big\ran_{\H_{(\mu,r)}}=-r^2.
\end{eqnarray*}
As by Lemma \ref{uncenter} one gets $\big\lan d\pi_{(\la^k,\al_k)}(T)\xi_k^s,\xi_k^s \big\ran_{\H_{(\la^k,\al_k)}}=\big\lan i \al_k \xi_k^s,\xi_k^s \big\ran_{\H_{(\la^k,\al_k)}}$, it follows that $\underset{k\to\iy}{\lim}\al_k=0$. Therefore, one can assume without restriction that $\al_k>0$ for large $k$ (by passing to a subsequence if necessary). The case $\al_k<0$ is very similar. \\
~\\
On the other hand, the sequence $\Big(\big\lan \ta_{\la^k}\otimes
W_{\al_k}(A)\xi_k^s,\xi_k^s \big\ran_{\H_{(\la^k,\al_k)}}\Big)_{k \in\N}$ \mbox{converges to the matrix} coefficient
$C^{\pi_{(\mu,r)}}_{\xi^s}(A,0,0)=C^{\tilde{\mu}^s}_{\phi_1^{\tilde{\mu}^s},\phi_1^{\tilde{\mu}^s}}(A)$
uniformly in each $A\in U(n)$. Hence, from this convergence, Orthogonality Relation (\ref{e1}) and the fact that $\Big\|C^{\tilde{\mu}^s}_{\phi_1^{\tilde{\mu}^s},\phi_1^{\tilde{\mu}^s}} \Big\|_{L^2\big(U(n)\big)}=\frac{1}{\sqrt{d_{\tilde{\mu}^s}}}$,
follows
\begin{eqnarray}\label{orthrelfolge}
\lim_{k\to\iy} \int \limits_{U(n)} \big\lan \ta_{\la^k}\otimes
W_{\al_k}(A)\xi_k^s,\xi_k^s \big\ran_{\H_{(\la^k,\al_k)}} \ol{ \big\lan
\ta_{\tilde{\mu}^s}(A)\phi_1^{\tilde{\mu}^s},\phi_1^{\tilde{\mu}^s} \big\ran_{\H_{(\mu,r)}}}~dA=\frac{1}{d_{\tilde{\mu}^s}}\ne0.
\end{eqnarray}
By (\ref{decoppositivealp}), one can write the expression $\big(\ta_{\la^k}\otimes W_{\al_k}\big) \res{U(n)}$ as
$$\big(\ta_{\la^k}\otimes W_{\al_k}\big) \res{U(n)}= \sum \limits_{\substack{\tilde{\la}^k\in P_n \\ \tilde{\la}_1^k \geq \la_1^k \geq...\geq \tilde{\la}_n^k \geq \la_n^k}} \tau_{ \tilde{\la}^k}$$
and, since for $k$ large enough the above integral is not $0$, again by the orthogonality relation, there has to be one $\tilde{\la}^k\in P_n$ with $\tilde{\la}_1^k \geq \la_1^k \geq...\geq \tilde{\la}_n^k \geq \la_n^k$ such that $\tilde{\la}^k=\tilde{\mu}^s$. But as $\tilde{\la}^k_s=\tilde{\mu}^s_s=\tilde{\mu}^s_{s+1}=\tilde{\la}^k_{s+1}$, one obtains that $\la^k_s=\tilde{\la}^k_s=\tilde{\mu}^s_s=\mu_s$ for $k$ large enough. As this is true for all $s \in \{1,..., n-1\}$, one gets $\la^k_j=\mu_j$ for all $j \in \{1,..., n-1\}$.\\
~\\
So, it remains to show that $\underset{k\to\infty}{\lim}\alpha_k\la^k_{n}=-\frac{r^2}{2}$. \\
Again, by the decomposition of $\big(\ta_{\la^k}\otimes W_{\al_k}\big) \res{U(n)}$ in (\ref{decoppositivealp}), one can decompose $\H_{(\la^k,\al_k)}$ as follows
$$\H_{(\la^k,\al_k)}= \sum \limits_{\substack{\tilde{\la}^k\in P_n \\ \tilde{\la}_1^k \geq \la_1^k \geq...\geq \tilde{\la}_n^k \geq \la_n^k}} \H_{\tilde{\la}^k}$$
and thus, for every $k \in \N$, the vector $\xi_k^s$ can be written as
$$\xi_k^s=\sum \limits_{\substack{\tilde{\la}^k\in P_n \\ \tilde{\la}_1^k \geq \la_1^k \geq...\geq \tilde{\la}_n^k \geq \la_n^k}} \xi_{\tilde{\la}^k}^s \:\:\:  \text{for}~\xi_{\tilde{\la}^k}^s \in \H_{\tilde{\la}^k} \:\:\: \forall k \in \N.$$
Let
\begin{eqnarray*}
 C_k:=\int \limits_{U(n)} \big\lan \ta_{\la^k}\otimes W_{\al_k}(A)\xi_k^s,\xi_k^s \big\ran_{\H_{(\la^k,\al_k)}}\ol{ \big\lan \ta_{\tilde{\mu}^s}(A)\phi_1^{\tilde{\mu}^s},\phi_1^{\tilde{\mu}^s} \big\ran_{\H_{(\mu,r)}}}~dA \:\:\: \forall k\in\N. \
\end{eqnarray*}
Then, with the Orthogonality Relation (\ref{e1}),
\begin{eqnarray*}\displaystyle
C_k&=& \sum \limits_{\substack{\tilde{\la}^k\in P_n \\ \tilde{\la}_1^k \geq \la_1^k \geq...\geq \tilde{\la}_n^k \geq \la_n^k}}~
\sum \limits_{\substack{\tilde{\ga}^k\in P_n \\ \tilde{\ga}_1^k \geq \la_1^k \geq...\geq \tilde{\ga}_n^k \geq \la_n^k}}
~\int \limits_{U(n)} \big\lan \ta_{\la^k}\otimes W_{\al_k}(A)\xi_{\tilde{\la}^k}^s,\xi_{\tilde{\ga}^k}^s \big\ran_{\H_{\tilde{\la}^k}} \ol{ \big\lan \ta_{\tilde{\mu}^s}(A)\phi_1^{\tilde{\mu}^s},\phi_1^{\tilde{\mu}^s} \big\ran_{\H_{(\mu,r)}}}~dA \\[5pt]
&=&{\small \sum \limits_{\substack{\tilde{\la}^k\in P_n \\ \tilde{\la}_1^k \geq \la_1^k \geq...\geq \tilde{\la}_n^k \geq \la_n^k}}~
\sum \limits_{\substack{\tilde{\ga}^k\in P_n \\ \tilde{\ga}_1^k \geq \la_1^k \geq...\geq \tilde{\ga}_n^k \geq \la_n^k}}
~\int \limits_{U(n)} \sum \limits_{\substack{\tilde{\nu}^k\in P_n \\ \tilde{\nu}_1^k \geq \la_1^k \geq...\geq \tilde{\nu}_n^k \geq \la_n^k}}\big\lan\tau_{\tilde{\nu}^k}(A)
\xi_{\tilde{\la}^k}^s,\xi_{\tilde{\ga}^k}^s \big\ran_{\H_{\tilde{\la}^k}}\ol{ \big\lan \ta_{\tilde{\mu}^s}(A)\phi_1^{\tilde{\mu}^s},\phi_1^{\tilde{\mu}^s} \big\ran_{\H_{(\mu,r)}}}~dA} \\[5pt]
&=& \sum \limits_{\substack{\tilde{\la}^k\in P_n \\ \tilde{\la}_1^k \geq \la_1^k \geq...\geq \tilde{\la}_n^k \geq \la_n^k}}~
\sum \limits_{\substack{\tilde{\ga}^k\in P_n \\ \tilde{\ga}_1^k \geq \la_1^k \geq...\geq \tilde{\ga}_n^k \geq \la_n^k}}
~\int \limits_{U(n)} \big\lan \ta_{\tilde{\la}^k}(A)\xi_{\tilde{\la}^k}^s,\xi_{\tilde{\ga}^k}^s \big\ran_{\H_{\tilde{\la}^k}} \ol{ \big\lan \ta_{\tilde{\mu}^s}(A)\phi_1^{\tilde{\mu}^s},\phi_1^{\tilde{\mu}^s} \big\ran_{\H_{(\mu,r)}}}~dA \\[2pt]
&&~\\
&=& \sum \limits_{\substack{\tilde{\ga}^k\in P_n \\ \tilde{\ga}_1^k \geq \la_1^k \geq...\geq \tilde{\ga}_n^k \geq \la_n^k}}
\frac{\big\lan \xi_{\tilde{\mu}^s}^s,\phi_1^{\tilde{\mu}^s} \big\ran_{\H_{\tilde{\mu}^s}} \big \lan \phi_1^{\tilde{\mu}^s}, \xi_{\tilde{\ga}^k}^s\big\ran_{\H_{\tilde{\mu}^s}}}{d_{\tilde{\mu}^s}} \\[1pt]
&&~\\
&=& \frac{\big\lan \xi_{\tilde{\mu}^s}^s,\phi_1^{\tilde{\mu}^s} \big\ran_{\H_{\tilde{\mu}^s}} \big \lan \phi_1^{\tilde{\mu}^s}, \xi_{\tilde{\mu}^s}^s\big\ran_{\H_{\tilde{\mu}^s}}}{d_{\tilde{\mu}^s}}
~=~ \frac{ \big|\big\lan \xi_{\tilde{\mu}^s}^s,\phi_1^{\tilde{\mu}^s} \big\ran_{\H_{\tilde{\mu}^s}} \big|^2}{d_{\tilde{\mu}^s}}.
\end{eqnarray*}
From (\ref{orthrelfolge}) it follows that
$$\big|\big\lan \xi_{\tilde{\mu}^s}^s,\phi_1^{\tilde{\mu}^s} \big\ran_{\H_{\tilde{\mu}^s}} \big|^2 \overset{k \to \iy}{\longrightarrow}1.$$
As
$$1=\| \xi_k^s \|_{\H_{(\la^k,\al_k)}}^2=\sum \limits_{\substack{\tilde{\la}^k\in P_n \\ \tilde{\la}_1^k \geq \la_1^k \geq...\geq \tilde{\la}_n^k \geq \la_n^k}} \big\|\xi_{\tilde{\la}^k}^s \big\|_{\H_{\tilde{\la}^k}}^2,$$
one can assume that $\xi_k^s= \phi_1^{\tilde{\mu}^s}$ for large $k \in \N$. Since $\la_j^k=\mu_j$ for all $k \in \N$ and for all $j \in \{1,...,n-1\}$, one gets for $s=n-1$
$$\tilde{\mu}^{n-1}=\la^k+m_k \:\:\: \text{for}~ m_k= \big(0,...,0,\mu_{n-1}-\la_n^k \big).$$
From now on, consider only $k$ large enough such that $\xi_k^{n-1}= \phi_1^{\tilde{\mu}^{n-1}}$. Then, $\xi_k^{n-1}$ is the \mbox{highest} weight vector with weight $\tilde{\mu}^{n-1}$ of length 1. Moreover,
$$\sum \limits_{\substack{\tilde{\la}^k\in P_n \\ \tilde{\la}_1^k \geq \la_1^k \geq...\geq \tilde{\la}_n^k \geq \la_n^k}} \H_{\tilde{\la}^k}=\H_{(\la^k,\al_k)}=\H_{\la^k} \otimes \F_{\al_k}(n)=\sum \limits_{m=0}^{\iy} \H_{\la^k} \otimes \P_m.$$
Every weight in the decomposition on the left hand side has multiplicity one, as mentioned in (\ref{decoppositivealp}), and therefore, this is the case for every weight appearing in the sum on the right hand side as well. From this, one can deduce that there exists one unique $M_k$ such that $\tilde{\mu}^{n-1}$, the weight of $\xi_k^{n-1} \in \H_{(\la^k,\al_k)}$, appears in $\H_{\la^k} \otimes \P_{M_k}$. \\
By [\ref{knapp}], Chapter IV.11,
every highest weight appearing in $\H_{\la^k} \otimes \P_{M_k}$ is the sum of the highest weight of $\H_{\la^k}$ and a weight of $\P_{M_k}$. Hence, $\tilde{\mu}^{n-1}$ is the sum of $\la^k$
and a weight of $\P_{M_k}$. From this follows that the mentioned weight of $\P_{M_k}$ has the same length as $\tilde{\mu}^{n-1}-\la^k=m_k$. Therefore, $M_k=|m_k|$, i.e. $\P_{M_k}=\P_{|m_k|}$. \\
Let $\big(\phi_j^{\la^k} \big)_{j \in \N} $ be an orthogonal weight vector basis of $\H_{\ta_{\la^k}}$, the corresponding Hilbert space of $\ta_{\la^k}$, and let  $\ga_j^k $ denote the weight of $\phi_j^{\la^k} $. Then, one gets
\begin{eqnarray*}
 \xi_k^{n-1}=\sum \limits_{(\gamma^k, \tilde{m}_k)\in \Omega_{\la^k}^{\tilde{\mu}^{n-1}}}~\sum_{j:\ga_j^k=\ga^k}  ~c_j^{k}\phi^{\ga^k}_j\otimes h_{\tilde{m}_k, \al_k}= \sum \limits_{(\gamma^k, \tilde{m}_k) \in \Omega_{\la^k}^{\tilde{\mu}^{n-1}}} ~\phi^{\ga^k} \otimes h_{\tilde{m}_k, \al_k},
\end{eqnarray*}
where $\phi^{\ga^k}=\sum \limits_{j:\ga_j^k=\ga^k}  ~c_j^{k}\phi^{\ga^k}_j$ is a uniquely determined eigenvector for $\T_n$ of the space $\H_{\la^k}$ \mbox{with weight $\ga^k$,} $\Omega_{\la^k}^{\tilde{\mu}^{n-1}}$ is the set of all pairs
$(\gamma^k, \tilde{m}_k)$ such that $\tilde{m}_k \in \N^n$ with $|\tilde{m}_k|=|m_k|$ and $\ga^k$ is a weight that appears in the representation $\tau_{\la^k}$ fulfilling $\gamma^k+ \tilde{m}_k=\tilde{\mu}^{n-1}$.
Furthermore,
$$\sum \limits_{(\gamma^k, \tilde{m}_k) \in \Omega_{\la^k}^{\tilde{\mu}^{n-1}}} \big|\phi^{\ga^k} \big|^2=1.$$
Let
\begin{eqnarray*}
 c_k:= \big\lan d\pi_{(\la^k,\al_k)}(\LL)\xi_k^{n-1},\xi_k^{n-1} \big\ran_{\H_{(\la^k,\al_k)}} \:\:\:  \forall k\in \N.
\end{eqnarray*}

Then, as seen above, $\lim \limits_{k \to \iy} c_k=-r^2 $ and  from Lemma \ref{unlap2} and the $U(n)$-invariance of $\LL$, it follows that
\begin{eqnarray*}
c_k&=& \big\lan d\pi_{(\la^k,\al_k)}(\LL)\xi_k^{n-1},\xi_k^{n-1} \big\ran_{\H_{(\la^k,\al_k)}} \\[5pt]
&=& \sum \limits_{(\gamma^k, \tilde{m}_k) \in \Omega_{\la^k}^{\tilde{\mu}^{n-1}}}~ \sum \limits_{(\tilde{\gamma}^k, \tilde{\tilde{m}}_k) \in \Omega_{\la^k}^{\tilde{\mu}^{n-1}}} \Big\lan d\pi_{(\la^k,\al_k)}(\LL)\phi^{\ga^k} \otimes h_{\tilde{m}_k, \al_k},\phi^{\tilde{\ga}^k} \otimes h_{\tilde{\tilde{m}}_k, \al_k} \Big\ran_{\H_{(\la^k,\al_k)}} \\[5pt]
&=& \sum \limits_{(\gamma^k, \tilde{m}_k) \in \Omega_{\la^k}^{\tilde{\mu}^{n-1}}} \big|\phi^{\ga^k} \big|^2  ~\big\lan d\pi_{(\la^k,\al_k)}(\LL)h_{\tilde{m}_k, \al_k},h_{\tilde{m}_k, \al_k} \big\ran_{\F_{\al_k}(n)} \\[2pt]
&=& \sum \limits_{(\gamma^k, \tilde{m}_k) \in \Omega_{\la^k}^{\tilde{\mu}^{n-1}}} \big|\phi^{\ga^k} \big|^2  \Big(- \al_k \big(n+2 | \tilde{m_k}| \big)\Big) \\[2pt]
&=& - \al_k \big(n+2 | m_k| \big) \sum \limits_{(\gamma^k, \tilde{m}_k) \in \Omega_{\la^k}^{\tilde{\mu}^{n-1}}} \big|\phi^{\ga^k} \big|^2  \\[2pt]
&=& - \al_k \big(n+2 | m_k| \big)~=~- \al_k \big(n+2 \mu_{n-1}-2 \la_n^k \big).
\end{eqnarray*}
As $\al_k \overset{k \to \iy}{\longrightarrow}0$, also $\al_k (n+2 \mu_{n-1})\overset{k \to \iy}{\longrightarrow}0$ and thus, $\underset{k\to\infty}{\lim}\alpha_k\la^k_{n}=-\frac{r^2}{2}$. \\

2) The fact that the sequence
$\big(\pi_{(\la^k,\alpha_k)}\big)_{k \in \N}$ converges to $\ta_{\la}$ in $\widehat{G_n}$ implies that for the unit vector $\phi_1^\la \in \H_{\la}^\iy$, there is a sequence of unit vectors
$(\xi_k)_{k \in \N} \subset \H_{(\la^k,\al_k)}^\iy$ such that
\begin{eqnarray*}
\big\lan d\pi_{(\la^k,\al_k)}(T)\xi_k,\xi_k \big\ran_{\H_{(\la^k,\al_k)}} \overset{k \to \iy}{\longrightarrow} \big\lan
d\ta_{\la}(T)\phi_1^\la,\phi_1^\la \big\ran_{\H_{\la}} \:\:\: \forall T\in\t_n \:\:\: \text{and}
\end{eqnarray*}
\begin{eqnarray}\label{dpi L=0}
\big\lan d\pi_{(\la^k,\al_k)}(\LL)\xi_k,\xi_k \big\ran_{\H_{(\la^k,\al_k)}} \overset{k \to \iy}{\longrightarrow} \big\lan
d\ta_{\la}(\LL)\phi_1^\la,\phi_1^\la \big\ran_{\H_{\la}}=0.
\end{eqnarray}
As above in the first part, by Lemma \ref{uncenter}, from the first convergence it follows that $\underset{k\to\iy}{\lim}\al_k=0$ and one can assume without restriction that $\al_k>0$ for large $k$. \\
On the other hand, $\Big(\big\lan \ta_{\la^k}\otimes
W_{\al_k}(A)\xi_k,\xi_k \big\ran_{\H_{(\la^k,\al_k)}}\Big)_{k \in \N}$ converges to
$C^{\la}_{\phi_1^{\la},\phi_1^{\la}}(A)$ uniformly in $A\in
U(n)$. Hence, as above one gets
\begin{eqnarray*}
\lim_{k\to\iy}\int \limits_{U(n)} \big\lan \ta_{\la^k}\otimes
W_{\al_k}(A)\xi_k,\xi_k \big\ran_{\H_{(\la^k,\al_k)}} \ol{ \big\lan
\ta_{\la}(A)\phi_1^{\la},\phi_1^{\la} \big\ran_{\H_{\la}}}~dA=\frac{1}{d_{\la}}\ne0.
\end{eqnarray*}
Again, like in the first part above, by (\ref{decoppositivealp}) and the orthogonality relation, one can deduce that
$\la_{1}\geq\la^k_{1}\geq....\geq\la_{n}\geq\la^k_{n}$ for large $k$. \\
~\\
So again, it remains to show that $\underset{k\to\infty}{\lim}\alpha_k\la^k_{n}=0$. \\
In the same manner as above, by replacing $\tilde{\mu}^{n-1}$ by $\la$, one can now show that for large $k \in \N$, it is possible to assume $\xi_k=\phi_1^{\la}$. So consider $k$ large enough in order for this equality to be true. Then $\xi_k$ is the highest weight vector of length 1 with weight $\la$. \\
Now,
$$\la=\la^k+m_k \:\:\: \text{for}~ m_k= \big(\la_1-\la_1^k,...,\la_n-\la_n^k \big),$$
where the sequences $\big(\la_1-\la_1^k \big)_{k \in \N}$,...,$\big(\la_{n-1}-\la_{n-1}^k \big)_{k \in \N}$ are bounded, because \mbox{$\la_{1}\geq\la^k_{1}\geq....\geq\la_{n}\geq\la^k_{n}$} for large $k$. \\
Again, by replacing $\tilde{\mu}^{n-1}$ by $\la$ in the proof of the first part above, one can also write $\xi_k$ as
$$\xi_k= \sum \limits_{(\gamma^k, \tilde{m}_k) \in \Omega_{\la^k}^{\la}} ~\phi^{\ga^k} \otimes h_{\tilde{m}_k, \al_k},$$
where $\phi^{\ga^k}$ is a uniquely determined eigenvector for $\T_n$ of $\H_{\la^k}$ with weight $\ga^k$ and $\Omega_{\la^k}^{\la}$ is the set of all pairs
$(\gamma^k, \tilde{m}_k)$ such that $\tilde{m}_k \in \N^n$ with $|\tilde{m}_k|=|m_k|$ and $\ga^k$ is a weight that appears in the representation $\tau_{\la^k}$ fulfilling $\gamma^k+ \tilde{m}_k=\la$.
Furthermore, again
$$\sum \limits_{(\gamma^k, \tilde{m}_k) \in \Omega_{\la^k}^{\la}} \big|\phi^{\ga^k} \big|^2=1.$$
Now, like in the first part above, by Lemma \ref{unlap2} and the $U(n)$-invariance of $\LL$,
\begin{eqnarray*}
 \lan d\pi_{(\la^k,\al_k)}(\LL)\xi_k,\xi_k \big\ran_{\H_{(\la^k,\al_k)}}
&=& - \al_k \big(n+2 | m_k| \big)\\[1pt]
&=&- \al_k \Big(n+2 \big(\la_1- \la_1^k \big)+...+2 \big(\la_{n-1}- \la_{n-1}^k \big)+2\la_n- 2\la_n^k \Big).
\end{eqnarray*}
By (\ref{dpi L=0}), as $\al_k \overset{k \to \iy}{\longrightarrow}0$, also $\al_k \Big(n+2 \big(\la_1- \la_1^k \big)+...+2\big(\la_{n-1}- \la_{n-1}^k \big)+2 \la_n \Big)\overset{k \to \iy}{\longrightarrow}0$ because of the boundedness of the sequences $\big(\la_1-\la_1^k \big)_{k \in \N}$,...,$\big(\la_{n-1}-\la_{n-1}^k \big)_{k \in \N}$. Therefore,  $\underset{k\to\infty}{\lim}\alpha_k\la^k_{n}=0$.
~\\
\qed
\section{The continuity of $\K $}
In this section, it will be shown that the inverse of the Kirillov-Lipsman mapping $\K $ is also continuous.
By Theorem \ref{rnenull}, it suffices to  consider  converging sequences  of orbits $(\O_{(\la^k,\al_k)})_{k \in \N} $ and to show that the corresponding  representations $(\pi_{(\la^k,\al_k)})_{k \in \N} $  converge in the same way.
\begin{theorem}\label{rueckweg}
~\\
Let $r>0$, $\rho_\mu\in\widehat{U(n-1)}$ and $\ta_\la\in\widehat{U(n)}$. \\
If $\underset{k\to\infty}{\lim}\alpha_k=0$ and the sequence $\big(\O_{(\la^k,\al_k)} \big)_{k\in\mathbb{N}}$ of
elements of the admissible orbit space $\g_n^\ddagger/G_n  $ satisfies one of the
following conditions:
\begin{enumerate}[(i)]
\item for $k$ large enough, $\alpha_k>0$, $\la_j^k=\mu_j$
for all $j \in \{1,..., n-1\}$ and
$\underset{k\to\infty}{\lim}\alpha_k\la^k_{n}=-\frac{r^2}{2}$,
\item for $k$ large enough, $\alpha_k<0$,
$\la_j^k=\mu_{j-1}$ for all $j \in \{2,..., n\}$ and
$\underset{k\to\infty}{\lim}\alpha_k\la^k_{1}=-\frac{r^2}{2}$,
\end{enumerate}
then the sequence $\big(\pi_{(\la^k,\alpha_k)}\big)_{k\in\mathbb{N}}$ converges to the representation $\pi_{(\mu,r)}$ in $\wh{G_n}$.
\end{theorem}

In order to prove this theorem, one needs the following proposition:
\begin{proposition}\label{prep-conj}
~\\
Let $r>0$, $\rho_\mu\in\widehat{U(n-1)}$ and $\ta_\la\in\widehat{U(n)}$. \\
Furthermore, let $\underset{k\to\infty}{\lim}\alpha_k=0$, $\alpha_k>0$ for large $k$ and consider the sequence $\big(\la^k)_{k\in\mathbb{N}}$ in $P_n$ fulfilling $\la_j^k=\mu_j$ for all $j \in \{1,..., n-1\}$ and $\underset{k\to\infty}{\lim}\alpha_k\la^k_{n}=-\frac{r^2}{2}$. \\
Denote $\tilde{\mu}:=\tilde{\mu}^{n-1}=( \mu_1,..., \mu_{n-1}, \mu_{n-1})$, $N_k:= \mu_{n-1}- \la_n^k$ and let $\ol{\P_{N_k}}$ be the space of conjugated homogeneous polynomials of degree $N_k$. \\
Moreover, define the representation $\pi^{(\tilde{\mu}, \al_k)}$ of $G_n$ on the subspace $\H_{\tilde{\mu}} \otimes \ol{\P_{N_k}} \otimes \P_{N_k}$ of the Hilbert space $\H_{\tilde{\mu}} \otimes \ol{\F_{\al_k}(n)} \otimes \F_{\al_k}(n)$ by
$$\pi^{(\tilde{\mu}, \al_k)}(A,z,t):= \tau_{\tilde{\mu}}(A) \otimes \ol{W_{\al_k}}(A)\otimes \big( \sigma_{\al_k}(z,t) \circ W_{\al_k}(A) \big) \:\:\: \forall (A,z,t) \in G_n.$$
Then, for any $\th\in \H_{\ti\mu} $ and for each $k \in \N$, there exist  vectors $\xi_k \in \H_{\tilde{\mu}} \otimes \ol{\P_{N_k}} \otimes \P_{N_k}$ such that for all $(A,z,t) \in G_n$ and for $\xi^\th:= \th\otimes 1 \in \H_{\tilde{\mu}} \otimes \H_{(0,r)}$,
\begin{eqnarray*}
\Big\langle \pi^{(\tilde{\mu}, \al_k)} (A,z,t) \xi_k, \xi_k \Big\rangle_{\H_{\tilde{\mu}} \otimes \ol{\F_{\al_k}(n)} \otimes \F_{\al_k}(n)} \overset{k \to \iy}{\longrightarrow}
\Big\langle \big(\tau_{\tilde{\mu}} \otimes \pi_{(0,r)}\big)(A,z,t) \xi^\th, \xi^\th \Big\rangle_{\H_{\tilde{\mu}} \otimes \H_{(0,r)}}
\end{eqnarray*}
uniformly on compacta.
\end{proposition}

Proof: \\
Let $m_k:=(0,...,0,N_k)$. Then, $\la^k= \tilde{\mu}+m_k$. Moreover, since $\underset{k\to\infty}{\lim}\alpha_k\la^k_{n}=-\frac{r^2}{2}$, $\underset{k\to\infty}{\lim}\alpha_k=0$ and $\al_k >0$, one gets $N_k \overset{k \to \iy}{\longrightarrow} \iy$. \\
Let $\phi\in \H_{\tilde{\mu}}$ and
let
\begin{eqnarray*}
 R_k:=\Big\vert\{q\in \N^n; \val {q}=N_k\}\Big\vert \:\:\: \forall k\in \N.
 \end{eqnarray*}
Then, $R_k $ is the dimension of the space $\P_{N_k}(n) =\P_{N_k}$ of complex polynomials of degree $N_k $ in $n $ variables.
Now, define
$$\xi_k:= \ph \otimes \bigg(\frac{1}{{R_k}^{\frac{1}{2}}}~\sum \limits_{\substack{q\in \N^n:\\ |q|=N_k}}\ol{h_{q,\al_k}} \otimes h_{q,\al_k} \bigg)=\frac{1}{{R_k}^{\frac{1}{2}}}~\sum \limits_{\substack{q\in \N^n:\\ |q|=N_k}} \ph \otimes \ol{h_{q,\al_k}} \otimes h_{q,\al_k} \in \H_{\tilde{\mu}} \otimes \ol{\P_{N_k}} \otimes \P_{N_k}.$$
Since ${R_k}^{\frac{1}{2}}$ is the norm of $\sum \limits_{\substack{q\in \N^n:\\ |q|=N_k}}\ol{h_{q,\al_k}} \otimes h_{q,\al_k}$, the vector $\frac{1}{{R_k}^{\frac{1}{2}}}~\sum \limits_{\substack{q\in \N^n:\\ |q|=N_k}}\ol{h_{q,\al_k}} \otimes h_{q,\al_k}$ has norm $1$. \\
Let $(A,z,t) \in G_n$. One has
\begin{eqnarray*}
 c_k(A,z,t)&:= &\Big\langle \pi^{(\tilde{\mu}, \al_k)}(A,z,t) \xi_k, \xi_k \Big\rangle_{\H_{\tilde{\mu}} \otimes \ol{\F_{\al_k}(n)} \otimes \F_{\al_k}(n)} \\[5pt]
&=& \Bigg\langle \pi^{(\tilde{\mu}, \al_k)}(A,z,t) \bigg(\frac{1}{{R_k}^{\frac{1}{2}}}~\sum \limits_{\substack{q\in \N^n:\\ |q|=N_k}} \ph \otimes \ol{h_{q,\al_k}} \otimes h_{q,\al_k} \bigg), \\
&&\:\:\:\:\:\: \frac{1}{{R_k}^{\frac{1}{2}}}~\sum \limits_{\substack{\tilde{q}\in \N^n:\\ |\tilde{q}|=N_k}} \ph \otimes \ol{h_{\tilde q,\al_k}} \otimes h_{\tilde q,\al_k} \Bigg\rangle_{\H_{\tilde{\mu}} \otimes \ol{\F_{\al_k}(n)} \otimes \F_{\al_k}(n)} \\[2pt]
&=& \Bigg\langle \pi^{(\tilde{\mu}, \al_k)}\big((\Id,z,t)(A,0,0) \big) \bigg(\frac{1}{{R_k}^{\frac{1}{2}}}~\sum \limits_{\substack{q\in \N^n:\\ |q|=N_k}} \ph \otimes \ol{h_{q,\al_k}} \otimes h_{q,\al_k} \bigg), \\[1pt]
&&\:\:\:\:\:\: \frac{1}{{R_k}^{\frac{1}{2}}}~\sum \limits_{\substack{\tilde{q}\in \N^n:\\ |\tilde{q}|=N_k}} \ph \otimes \ol{h_{\tilde q,\al_k}} \otimes h_{\tilde q,\al_k} \Bigg\rangle_{\H_{\tilde{\mu}} \otimes \ol{\F_{\al_k}(n)} \otimes \F_{\al_k}(n)} \\[2pt]
&=& \Bigg\langle \frac{1}{{R_k}^{\frac{1}{2}}}~\sum \limits_{\substack{q\in \N^n:\\ |q|=N_k}} \ph \otimes \ol{W_{\al_k}(A) h_{q,\al_k}} \otimes \Big(\sigma_{\al_k}(z,t)\circ W_{\al_k}(A) h_{q,\al_k} \Big), \\
&&\:\:\:\:\:\: \frac{1}{{R_k}^{\frac{1}{2}}}~\tau_{\tilde{\mu}}(A\inv)\ph \otimes \bigg(\sum \limits_{\substack{\tilde{q}\in \N^n:\\ |\tilde{q}|=N_k}}  \ol{h_{\tilde q,\al_k}} \otimes h_{\tilde q,\al_k}\bigg) \Bigg\rangle_{\H_{\tilde{\mu}} \otimes \ol{\F_{\al_k}(n)} \otimes \F_{\al_k}(n)}.
\end{eqnarray*}

Now, one can write
\begin{eqnarray*}
W_{\al_k}(A) h_{q,\al_k}&=&\sum \limits_{\substack{m \in \N^n:\\ |m|=N_k}} w_{m,q}^k(A) h_{m, \al_k}\:\:\:\text{and} \\[2pt]
\ol{W_{\al_k}(A) h_{q,\al_k}}&=&\sum \limits_{\substack{m \in \N^n:\\ |m|=N_k}} \ol{w_{m,q}^k(A) h_{m, \al_k}}
\end{eqnarray*}
with $w_{m,q}^k(A) \in \C$. Because of the unitarity of the matrix $W_{\al_k}(A)$, one gets for $m,m' \in \N^n$ with $|m|=|m'|=N_k$,
\begin{eqnarray*}
\sum \limits_{\substack{q \in \N^n:\\ |q|=N_k}} w_{m,q}^k(A) \ol{w_{m',q}^k(A)}=
\begin{cases}
0 & \textrm{if}~ m \not=m', \\
1 & \textrm{if}~ m=m'.
\end{cases}
\end{eqnarray*}
Hence,
\begin{eqnarray*}
c_k(A,z,t)&=&\Bigg\langle \frac{1}{{R_k}^{\frac{1}{2}}}~\sum \limits_{\substack{m\in \N^n:\\ |m|=N_k}} \ph \otimes \ol{h_{m,\al_k}} \otimes \Big(\sigma_{\al_k}(z,t)h_{m,\al_k} \Big), \\[1pt]
&&\:\:\:\:\:\: \frac{1}{{R_k}^{\frac{1}{2}}}~\tau_{\tilde{\mu}}(A\inv)\ph \otimes \bigg(\sum \limits_{\substack{\tilde{q}\in \N^n:\\ |\tilde{q}|=N_k}}  \ol{h_{\tilde q,\al_k}} \otimes h_{\tilde q,\al_k}\bigg) \Bigg\rangle_{\H_{\tilde{\mu}} \otimes \ol{\F_{\al_k}(n)} \otimes \F_{\al_k}(n)} \\
&=& \Big\langle \tau_{\tilde{\mu}}(A)\phi,\ph \Big\rangle_{\H_{\tilde{\mu}}}~ \\[3pt]
&& \frac{1}{R_k}~
\Bigg\langle \sum \limits_{\substack{m\in \N^n:\\ |m|=N_k}} \ol{h_{m,\al_k}} \otimes \Big(\sigma_{\al_k}(z,t)h_{m,\al_k} \Big),  \sum \limits_{\substack{\tilde{q}\in \N^n:\\ |\tilde{q}|=N_k}}  \ol{h_{\tilde q,\al_k}} \otimes h_{\tilde q,\al_k} \Bigg\rangle_{\ol{\F_{\al_k}(n)} \otimes \F_{\al_k}(n)} \\[3pt]
&=& \Big\langle \tau_{\tilde{\mu}}(A)\phi,\ph \Big\rangle_{\H_{\tilde{\mu}}}~ \frac{1}{R_k}~\sum \limits_{\substack{m, \tilde{q} \in \N^n:\\ |m|=|\tilde{q}|=N_k}}  \Big\langle \ol{h_{m,\al_k}},\ol{h_{\tilde q,\al_k}} \Big\rangle_{\ol{\F_{\al_k}(n)}} \Big\langle \sigma_{\al_k}(z,t) h_{m,\al_k},h_{\tilde q,\al_k} \Big\rangle_{\F_{\al_k}(n)} \\[3pt]
&=& \Big\langle \tau_{\tilde{\mu}}(A)\phi,\ph \Big\rangle_{\H_{\tilde{\mu}}}~ \frac{1}{R_k}~\sum \limits_{\substack{q\in \N^n:\\ |q|=N_k}} \Big\langle \sigma_{\al_k}(z,t) h_{q,\al_k},h_{ q,\al_k} \Big\rangle_{\F_{\al_k}(n)} \\[3pt]
&=& \Big\langle \tau_{\tilde{\mu}}(A)\phi,\ph \Big\rangle_{\H_{\tilde{\mu}}}~ \frac{1}{R_k}~\sum \limits_{\substack{q\in \N^n:\\ |q|=N_k}}
~\int \limits_{\C^n} e^{i \al_k t - \frac{\al_k}{4} |z|^2} e^{- \frac{\al_k}{2} \langle w,z \rangle_{\C^n}} h_{q,\al_k} (z+w) \overline{h_{q,\al_k}(w)} e^{- \frac{\al_k}{2}|w|^2}dw \\[3pt]
&=& \Big\langle \tau_{\tilde{\mu}}(A)\phi,\ph \Big\rangle_{\H_{\tilde{\mu}}}~ \frac{1}{R_k}~\sum \limits_{\substack{q\in \N^n:\\ |q|=N_k}}
\Big(\frac{\al_k}{2 \pi} \Big)^n ~\frac{\al_k^{N_k}}{2^{N_k}}\frac{1}{q!}~e^{i \al_k t - \frac{\al_k}{4} |z|^2} \\[1pt]
&&~~~~~~~~~~~~~~~~~~~~~~~~~~~~~~~~~~~~~~
\int \limits_{\C^n}  e^{- \frac{\al_k}{2} \langle w,z \rangle_{\C^n}}(z+w)^{q} \overline{w}^{q} e^{- \frac{\al_k}{2}|w|^2}dw.
\end{eqnarray*}
Now, by the binomial theorem, letting $\binom{q}{l}:=\binom{q_1}{l_1} \cdots \binom{q_n}{l_n}$ for $q=(q_1,...,q_n) \in \N^n$ and \mbox{$l=(l_1,...,l_n) \in \N^n$},
$$(z+w)^q= \sum \limits_{l_1=0}^{q_1} \binom{q_1}{l_1} z_1^{q_1-l_1}w_1^{l_1} \cdots \sum \limits_{l_n=0}^{q_n} \binom{q_n}{l_n} z_n^{q_n-l_n}w_n^{l_n}=\sum \limits_{\substack{l:=(l_1,...,l_n) \in \N^n:\\ l_1 \leq q_1,...,l_n \leq q_n}} \binom{q}{l} z^{q-l}w^{l}.$$
Thus, one gets for $q \in \N^n$ with $|q|=N_k$,
\begin{eqnarray*}
&&\Big(\frac{\al_k}{2 \pi} \Big)^n ~\frac{\al_k^{N_k}}{2^{N_k}}\frac{1}{q!}~e^{i \al_k t - \frac{\al_k}{4} |z|^2}
\int \limits_{\C^n}  e^{- \frac{\al_k}{2} \langle w,z \rangle_{\C^n}}(z+w)^{q} \overline{w}^{q} e^{- \frac{\al_k}{2}|w|^2}dw 
\\[1pt]
&=&\Big(\frac{\al_k}{2 \pi} \Big)^n ~\frac{\al_k^{N_k}}{2^{N_k}}\frac{1}{q!}~e^{i \al_k t - \frac{\al_k}{4} |z|^2}
\sum \limits_{\substack{l:=(l_1,...,l_n) \in \N^n:\\ l_1 \leq q_1,...,l_n \leq q_n}} \binom{q}{l} z^{q-l} \int \limits_{\C^n}  e^{- \frac{\al_k}{2} \langle w,z \rangle_{\C^n}}
{w}^{l} \overline{w}^{q} e^{- \frac{\al_k}{2}|w|^2} dw.
\end{eqnarray*}

The integrals in $w_m$ for $m \in \{1,...,n\}$ can be written as follows:
$$\sum \limits_{j_m=0}^{\iy}~ \int \limits_{\C} \frac{{w_m}^{j_m} (-\overline{z_m} )^{j_m}}{j_m!} ~\Big(\frac{\al_k}{2} \Big)^{j_m} e^{- \frac{\al_k}{2}|w_m|^2} {w_m} ^{l_m} \overline{w_m}^{q_m} dw_m.$$

Therefore,
\begin{eqnarray*}
\et_{k,q}(z)&:=& \Big(\frac{\al_k}{2 \pi} \Big)^n ~\frac{\al_k^{N_k}}{2^{N_k}}\frac{1}{q!}~e^{i \al_k t - \frac{\al_k}{4} |z|^2}
\sum \limits_{\substack{l:=(l_1,...,l_n) \in \N^n:\\ l_1 \leq q_1,...,l_n \leq q_n}} \binom{q}{l} z^{q-l} \int \limits_{\C^n}  e^{- \frac{\al_k}{2} \langle w,z \rangle_{\C^n}} {w}^{l} \overline{w}^{q} e^{- \frac{\al_k}{2}|w|^2} dw \\[2pt]
&=& \Big(\frac{\al_k}{2 \pi} \Big)^n ~\frac{\al_k^{N_k}}{2^{N_k}}\frac{1}{q!}~e^{i \al_k t - \frac{\al_k}{4} |z|^2}
\sum \limits_{\substack{l:=(l_1,...,l_n) \in \N^n:\\ l_1 \leq q_1,...,l_n \leq q_n}} ~\sum \limits_{j \in \N^n} \Big(\frac{\al_k}{2} \Big)^{|j|} \binom{q}{l} z^{q-l}~ \frac{(-\overline{z} )^{j}}{j!} \\
&& \int \limits_{\C^n} {w}^{j+l} \overline{w}^{q} e^{- \frac{\al_k}{2}|w|^2} dw.
\end{eqnarray*}
Because of the orthogonality of the functions $ \C^n \to \C^n,x \mapsto x^a$ and $\C^n \to \C^n,x \mapsto x^b$ for $a,b \in \N^n$ with respect to the scalar product of the Fock space, $j+l=q$, i.e. $l=q-j$. As $\no{h_{q,\al_k}}{\F_{\al_k}(n)}=1 $,
 the norm of the function $z \mapsto  z^{q} $ is
$\sqrt{\frac{1}{\big(\frac{\al_k}{2 \pi} \big)^n \frac{\al_k^{N_k}}{2^{N_k}q!}}}$ and hence,
\begin{eqnarray}\label{computations for q}
\nn\et_{k,q}(z)&=&\Big(\frac{\al_k}{2 \pi} \Big)^n ~\frac{\al_k^{N_k}}{2^{N_k}}\frac{1}{q!}~e^{i \al_k t - \frac{\al_k}{4} |z|^2}
\sum \limits_{\substack{j:=(j_1,...,j_n) \in \N^n:\\ j_1 \leq q_1,...,j_n \leq q_n}}  \Big(\frac{\al_k}{2} \Big)^{|j|} \binom{q}{q-j} z^{j}~ \frac{(-\overline{z})^{j}}{j!} \big\| \cdot ^{~q} \big\|_{\F_{\al_k}(n)}^2 \\[4pt]
&=&e^{i \al_k t - \frac{\al_k}{4} |z|^2}
\sum \limits_{\substack{j:=(j_1,...,j_n) \in \N^n:\\ j_1 \leq q_1,...,j_n \leq q_n}} \Big(\frac{\al_k}{2} \Big)^{|j|}
\frac{q!}{(q-j)!} \frac{z^j(-\ol{z})^{j}}{\big(j!\big)^2}.
\end{eqnarray}

Thus,
\begin{eqnarray*}
c_k(A,z,t)&=&\Big\langle \pi^{(\tilde{\mu}, \al_k)}(A,z,t) \xi_k, \xi_k \Big\rangle_{\H_{\tilde{\mu}} \otimes \ol{\F_{\al_k}(n)} \otimes \F_{\al_k}(n)} \\[5pt]
&=&\Big\langle \tau_{\tilde{\mu}}(A)\big( \ph \big), \ph \Big\rangle_{\H_{\tilde{\mu}}}~ e^{i \al_k t - \frac{\al_k}{4} |z|^2} \frac{1}{R_k} \sum \limits_{\substack{q\in \N^n:\\ |q|=N_k}}~ \sum \limits_{\substack{j:=(j_1,...,j_n) \in \N^n:\\ j_1 \leq q_1,...,j_n \leq q_n}} \Big(\frac{\al_k}{2} \Big)^{|j|}
\frac{q!}{(q-j)!} \frac{z^j(-\ol{z})^{j}}{\big(j!\big)^2}.
\end{eqnarray*}

Now, regard
\begin{eqnarray*}
\zeta_k(z)&:=&\frac{1}{R_k} \sum \limits_{\substack{q\in \N^n:\\ |q|=N_k}}~ \sum \limits_{\substack{j:=(j_1,...,j_n) \in \N^n:\\ j_1 \leq q_1,...,j_n \leq q_n}} \Big(\frac{\al_k}{2} \Big)^{|j|} \frac{q!}{(q-j)!} \frac{z^j(-\ol{z})^{j}}{\big(j!\big)^2} \\[3pt]
&=&\frac{1}{R_k} \sum \limits_{\substack{q_1,...,q_n \in \N:\\ q_1+...+q_n=N_k}}~ \sum \limits_{\substack{j:=(j_1,...,j_n) \in \N^n:\\ j_1 \leq q_1,...,j_n \leq q_n}} \Big(\frac{\al_k}{2} \Big)^{j_1+...+j_n} \Big( q_1 (q_1-1) \cdots (q_1-j_1+1) \Big) \\
&&\:\:\:\:\:\:\:\:\:\:\:\:\:\:\:\:\:\:\:\:\:\:\:\:\:\:\:\:\:\:\:\:\:\:\:\:\:\:\:\:\:\:\:\:\:\:\:\:\:\:\:\:\:\:\:\:\:\:\:\:\:\:\:\: \cdots \Big( q_n (q_n-1) \cdots (q_n-j_n+1) \Big) \frac{z^j(-\ol{z})^{j}}{\big(j!\big)^2}.
\end{eqnarray*}

Then, fixing large $k \in \N$, since $\underset{k\to\infty}{\lim}\alpha_k N_k=\frac{r^2}{2}$, one gets for $j=(j_1,...,j_n) \in \N^n$,
\begin{eqnarray}\label{sum}
\nn\vert \zeta_k(z)\vert&=& \Bigg|\frac{1}{R_k}~ \sum \limits_{\substack{q_1 \in \N_{\geq j_1},...,q_n \in \N_{\geq j_n}:\\ q_1+...+q_n=N_k}} \Big(\frac{\al_k N_k}{2} \Big)^{j_1+...+j_n}~ \frac{q_1 (q_1-1) \cdots (q_1-j_1+1)}{{N_k}^{j_1}} \\[2pt]
&&\:\:\:\: \cdots \frac{q_n (q_n-1) \cdots (q_n-j_n+1)}{{N_k}^{j_n}}  \frac{z^j(-\ol{z})^{j}}{\big(j!\big)^2} \Bigg| \\[4pt]
\nn&\leq& \bigg| \frac{1}{R_k} ~\sum \limits_{\substack{q_1 \in \N_{\geq j_1},...,q_n \in \N_{\geq j_n}:\\ q_1+...+q_n=N_k}}  \Big(\frac{r^2}{4}+1 \Big)^{j_1+...+j_n}  \frac{z^j(-\ol{z})^{j}}{\big(j!\big)^2} \bigg|~=~ \bigg( \Big(\frac{r^2}{4} +1 \Big)z \bigg)^j \ol{z}^j \frac{1}{(j!)}.
\end{eqnarray}
 The above expression does not depend on $k$ and
$$\sum \limits_{j:=(j_1,...,j_n) \in \N^n} \bigg( \Big(\frac{r^2}{4} +1 \Big)z \bigg)^j \ol{z}^j \frac{1}{(j!)} ~=~ \exp{ \bigg( \Big(\frac{r^2}{4} +1 \Big)z \ol{z} \bigg)}  < \iy.$$
So, by the theorem of Lebesgue, the sum in (\ref{sum}) converges and it suffices to regard the limit of each summand by itself. Hence, for $j=(j_1,...,j_n) \in \N^n$, one has
\begin{eqnarray*}
\zeta_k(z)&\cong&\frac{1}{R_k}~ \sum \limits_{\substack{q_1 \in \N_{\geq j_1},...,q_n \in \N_{\geq j_n}:\\ q_1+...+q_n=N_k}} \Big(\frac{r^2}{4} \Big)^{j_1+...+j_n}~ \frac{q_1 (q_1-1) \cdots (q_1-j_1+1)}{{N_k}^{j_1}}  \\
&&\:\:\:\:\:\:\:\:\:\:\:\:\:\:\:\:\:\:\:\:\:\:\:\:\:\:\:\:\:\:\:\:\:\:\:\:\:\:\:\:\:\:\:\:\:\:\: \cdots \frac{q_n (q_n-1) \cdots (q_n-j_n+1)}{{N_k}^{j_n}}~ \frac{z^j(-\ol{z})^{j}}{\big(j!\big)^2} \\[4pt]
&=&\frac{1}{R_k}~ \sum \limits_{\substack{q_1 \in \N_{\geq j_1},...,q_n \in \N_{\geq j_n}:\\ q_1+...+q_n=N_k}} \Big(\frac{r^2}{4} \Big)^{j_1+...+j_n}~
\frac{q_1}{N_k} \Big(\frac{q_1}{N_k}-\frac{1}{R_k} \Big) \cdots \Big(\frac{q_1}{N_k}-\frac{j_1-1}{N_k} \Big) \\
&&\:\:\:\:\:\:\:\:\:\:\:\:\:\:\:\:\:\:\:\:\:\:\:\:\:\:\:\:\:\:\:\:\:\:\:\:\:\:\:\:\:\:\:\:\:\:\: \cdots \frac{q_n}{N_k} \Big(\frac{q_n}{N_k}-\frac{1}{R_k} \Big) \cdots \Big(\frac{q_n}{N_k}-\frac{j_n-1}{N_k} \Big)~ \frac{z^j(-\ol{z})^{j}}{\big(j!\big)^2} \\[4pt]
&=&\frac{1}{R_k}~ \sum \limits_{\substack{q_1 \in \N_{\geq j_1},...,q_{n-1} \in \N_{\geq j_{n-1}}:\\ q_1+...+q_{n-1}\leq N_k-j_n}} \Big(\frac{r^2}{4} \Big)^{j_1+...+j_n}~
\frac{q_1}{N_k} \Big(\frac{q_1}{N_k}-\frac{1}{R_k} \Big) \cdots \Big(\frac{q_1}{N_k}-\frac{j_1-1}{N_k} \Big) \\
&&\:\:\:\:\:\:\:\:\:\:\:\:\:\:\:\:\:\:\:\:\:\:\:\:\:\:\:\:\:\:\:\:\:\:\:\:\:\:\:\:\:\:\:\:\:\:\:~~~~ \cdots \frac{q_{n-1}}{N_k} \Big(\frac{q_{n-1}}{N_k}-\frac{1}{R_k} \Big) \cdots \Big(\frac{q_{n-1}}{N_k}-\frac{j_{n-1}-1}{N_k} \Big) \\[2pt]
&&\:\:\:\:\:\:\:\:\:\:\:\:\:\:\:\:\:\:\:\:\:\:\:\:\:\:\:\:\:\:\:\:\:\:\:\:\:\:\:\:\:\:\:\:\:\:\:~~~~ \cdot \Big(1- \frac{q_1+...+q_{n-1}}{N_k}\Big) \bigg(\Big(1- \frac{q_1+...+q_{n-1}}{N_k}\Big)-\frac{1}{R_k} \bigg) \\[2pt]
&&\:\:\:\:\:\:\:\:\:\:\:\:\:\:\:\:\:\:\:\:\:\:\:\:\:\:\:\:\:\:\:\:\:\:\:\:\:\:\:\:\:\:\:\:\:\:\:~~~~\cdots \bigg(\Big(1- \frac{q_1+...+q_{n-1}}{N_k}\Big)-\frac{j_{n}-1}{N_k} \bigg)~\frac{z^j(-\ol{z})^{j}}{\big(j!\big)^2}.
\end{eqnarray*}
Now, define for $k \in \N$ the function $F_k:[0,1]^{n-1} \to \R$ by
\begin{eqnarray*}
&&F_k(s_1,...,s_{n-1}):=\Big(\frac{r^2}{4} \Big)^{j_1+...+j_n} s_1 \Big(s_1- \frac{1}{R_k} \Big) \cdots \Big( s_1- \frac{j_1-1}{N_k} \Big) \cdots s_{n-1} \Big(s_{n-1}- \frac{1}{R_k} \Big) \cdots \Big( s_{n-1}- \frac{j_{n-1}-1}{N_k} \Big) \\[1pt]
&&\:\:\:\:\:\:\:\: \:\:\:\:\:\:\:\:\:\:\:\:\:\:\:\:\:\:\:\:\:\:\:\:\cdot \big(1-(s_1+...+s_{n-1})\big) \Big(\big(1-(s_1+...+s_{n-1})\big)- \frac{1}{R_k} \Big)
\\
&& \:\:\:\:\:\:\:\: \:\:\:\:\:\:\:\:\:\:\:\:\:\:\:\:\:\:\:\:\:\:\:\:\cdots \Big(\big(1-(s_1+...+s_{n-1})\big)- \frac{j_n-1}{N_k} \Big)~\frac{z^j(-\ol{z})^{j}}{\big(j!\big)^2}.
\end{eqnarray*}

Then, for $\varepsilon >0$ and large $k \in \N$,
$$\Bigg| \frac{1}{R_k}~ \sum \limits_{q_1,...,q_{n-1} \in \N_{\leq N_k}} F_k \Big(\frac{q_1}{N_k},...,\frac{q_{n-1}}{N_k} \Big)~-~\int \limits_0^1 \cdots \int \limits_0^1 F_k(s_1,...,s_{n-1}) ds_1...ds_{n-1} \Bigg|  < \varepsilon.$$
Since $F_k\Big(\frac{q_1}{N_k},...,\frac{q_{n-1}}{N_k} \Big) =0$, if $q_1 <j_1$, $q_2 < j_2$ or ... or $q_{n-1} < j_{n-1}$ or $q_1+...+q_{n-1}> N_k-j_n$, it follows that
$$\Bigg| \frac{1}{R_k}~ \sum \limits_{\substack{q_1 \in \N_{\geq j_1},...,q_{n-1} \in \N_{\geq j_{n-1}}:\\ q_1+...+q_{n-1}\leq N_k-j_n}} F_k \Big(\frac{q_1}{N_k},...,\frac{q_{n-1}}{N_k} \Big)~-~\int \limits_0^1 \cdots \int \limits_0^1 F_k(s_1,...,s_{n-1}) ds_1...ds_{n-1} \Bigg|  < \varepsilon.$$
Furthermore, $F_k$ converges pointwise to the function $F: [0,1]^{n-1} \to \R$ defined by
$$F(s_1,...,s_{n-1}):=\Big(\frac{r^2}{4} \Big)^{j_1+...+j_n} {s_1}^{j_1} \cdots {s_{n-1}}^{j_{n-1}} \cdot  \big(1-(s_1+...+s_{n-1}) \big)^{j_n}~ \frac{z^j(-\ol{z})^{j}}{\big(j!\big)^2}$$
and thus, by the theorem of Lebesgue for integrals,
$$\lim \limits_{k \to \iy} \int \limits_0^1 \cdots \int \limits_0^1 F_k(s_1,...,s_{n-1}) ds_1...ds_{n-1}~ =~ \int \limits_0^1 \cdots \int \limits_0^1 F(s_1,...,s_{n-1}) ds_1...ds_{n-1}.$$
From these observations, one can now deduce that
\begin{eqnarray*}
&&\frac{1}{R_k}~ \sum \limits_{\substack{q_1 \in \N_{\geq j_1},...,q_{n-1} \in \N_{\geq j_{n-1}}:\\ q_1+...+q_{n-1}\leq N_k-j_n}} \Big(\frac{r^2}{4} \Big)^{j_1+...+j_n}~
\frac{q_1}{N_k} \Big(\frac{q_1}{N_k}-\frac{1}{R_k} \Big) \cdots \Big(\frac{q_1}{N_k}-\frac{j_1-1}{N_k} \Big) \\
&&\:\:\:\:\:\:\:\:\:\:\:\:\:\:\:\:\:\:\:\:\:\:\:\:\:\:\:\:\:\:\:\:\:\:\:\:\:\:\:\:\:\:\:\:\:\:\:~~~~ \cdots \frac{q_{n-1}}{N_k} \Big(\frac{q_{n-1}}{N_k}-\frac{1}{R_k} \Big) \cdots \Big(\frac{q_{n-1}}{N_k}-\frac{j_{n-1}-1}{N_k} \Big) \\[4pt]
&&\:\:\:\:\:\:\:\:\:\:\:\:\:\:\:\:\:\:\:\:\:\:\:\:\:\:\:\:\:\:\:\:\:\:\:\:\:\:\:\:\:\:\:\:\:\:\:~~~~ \cdot \Big(1- \frac{q_1+...+q_{n-1}}{N_k}\Big) \bigg(\Big(1- \frac{q_1+...+q_{n-1}}{N_k}\Big)-\frac{1}{R_k} \bigg) \\[4pt]
&&\:\:\:\:\:\:\:\:\:\:\:\:\:\:\:\:\:\:\:\:\:\:\:\:\:\:\:\:\:\:\:\:\:\:\:\:\:\:\:\:\:\:\:\:\:\:\:~~~~ \cdots \bigg(\Big(1- \frac{q_1+...+q_{n-1}}{N_k}\Big)-\frac{j_{n}-1}{N_k} \bigg)~\frac{z^j(-\ol{z})^{j}}{\big(j!\big)^2} \\
&=& \frac{1}{R_k}~ \sum \limits_{\substack{q_1 \in \N_{\geq j_1},...,q_{n-1} \in \N_{\geq j_{n-1}}:\\ q_1+...+q_{n-1}\leq N_k-j_n}} F_k \Big(\frac{q_1}{N_k},...,\frac{q_{n-1}}{N_k} \Big)
~\overset{k \to \iy}{\longrightarrow}~ \int \limits_0^1 \cdots \int \limits_0^1 F(s_1,...,s_{n-1}) ds_1...ds_{n-1} \\
\end{eqnarray*}
\begin{eqnarray*}
&=& \int \limits_{0}^1 \cdots \int \limits_{0}^1 \Big(\frac{r^2}{4} \Big)^{j_1+...+j_n} {s_1}^{j_1} \cdots {s_{n-1}}^{j_{n-1}} \big(1- (s_1+...+s_{n-1})\big)^{j_n}~ \frac{z^j(-\ol{z})^{j}}{\big(j!\big)^2} ~d s_1... ds_{n-1}.~~~~
\end{eqnarray*}
Therefore,
\begin{eqnarray}\label{mass}
\nn c_k(A,z,t) &\overset{k \to \iy}{\longrightarrow}& \Big\langle \tau_{\tilde{\mu}}(A)\big( \ph \big), \ph \Big\rangle_{\H_{\tilde{\mu}}}~ \sum \limits_{j:=(j_1,...,j_n) \in \N^n}\: \int \limits_{0}^1 \cdots \int \limits_0^{1} \Big(\frac{r^2}{4} \Big)^{j_1+...+j_n} {s_1}^{j_1} \cdots {s_{n-1}}^{j_{n-1}} \\
\nn && \:\:\:\:\:\:\:\:\:\:\:\:\:\:\:\:\:\:\:\:\:\:\:\:\:\:\:\:\:\:\:\:\:\:\:\:\:\:\:\:\:\:\:\:\:\:\:\:\:\:\:\:\:\:\:\:\:\:\:\:\:\:\:\:\:\:\:\: \big((1- (s_1+...+s_{n-1})\big)^{j_n}~ \frac{z^j(-\ol{z})^{j}}{\big(j!\big)^2} ~d s_1 ... d s_{n-1} \\
\nn &=&\Big\langle \tau_{\tilde{\mu}}(A)\big( \ph \big), \ph \Big\rangle_{\H_{\tilde{\mu}}} ~\int \limits_{0}^{1} \cdots \int \limits_0^{1}~ \underbrace{\Bigg(\sum \limits_{j_1 \in \N} \frac{\Big(\frac{- |z_1|^2 s_1 r^2}{4} \Big)^{j_1}}{\big(j_1!\big)^2}\Bigg)}_{\text{Bessel function}} \\[5pt]
   \nn &&\cdots \underbrace{\Bigg(\sum \limits_{j_{n-1} \in \N} \frac{\Big(\frac{- |z_{n-1}|^2 s_{n-1} r^2}{4} \Big)^{j_{n-1}}}{\big(j_{n-1}!\big)^2}\Bigg)}_{\text{Bessel function}}
 \underbrace{\Bigg(\sum \limits_{j_n \in \N} \frac{\Big(\frac{- |z_n|^2 \big(1-(s_1+...+s_{n-1} )\big) r^2}{4} \Big)^{j_n}}{\big(j_n!\big)^2}\Bigg)}_{\text{Bessel function}} ~d s_1... ds_{n-1} \\[5pt]
\nn &=&\Big\langle \tau_{\tilde{\mu}}(A)\big( \ph \big), \ph \Big\rangle_{\H_{\tilde{\mu}}} ~\Big(\frac{1}{2 \pi} \Big)^n \int \limits_{0}^{1} \cdots \int \limits_0^{1} \int \limits_0^{2 \pi} \cdots \int \limits_0^{2 \pi}
e^{-ir \Re{\big(e^{ia_1} \sqrt{s_1}~ \ol{z_1}\big)}} \\[3pt]
\nn && \:\:\:\:\:\:\:\:\:\:\:\:\:
 \cdots e^{-ir \Re{\big(e^{ia_{n-1}} \sqrt{s_{n-1}}~ \ol{z_{n-1}} \big)}} e^{-ir \Re{\big(e^{ia_n} \sqrt{1- (s_1+...+s_{n-1})} ~\ol{z_n}\big)}} da_1... da_n d s_1...d s_{n-1} \\[3pt]
\nn &=&\Big\langle \tau_{\tilde{\mu}}(A)\big( \ph \big), \ph \Big\rangle_{\H_{\tilde{\mu}}} ~\Big(\frac{1}{2 \pi} \Big)^n \int \limits_{0}^{1} \cdots \int \limits_0^{1} \int \limits_0^{2 \pi} \cdots \int \limits_0^{2 \pi}\\[3pt]
\nn && ~~~~~~~~~
e^{-ir \big( \big(\sqrt{s_1} e^{ia_1},...,\sqrt{s_{n-1}} e^{ia_{n-1}}, \sqrt{1- (s_1+...+s_{n-1})} e^{ia_n}\big),(z_1,...,z_n) \big)_{\C^n}} da_1... da_n d s_1...d s_{n-1} \\[3pt]
\nn &=&\Big\langle \tau_{\tilde{\mu}}(A)\big( \ph \big), \ph \Big\rangle_{\H_{\tilde{\mu}}}~\int \limits_{S^{n}} e^{-i (rv,z)_{\C^n}} d\sigma(v) \\[2pt]
\nn &=& \Big\langle \tau_{\tilde{\mu}}(A)\big( \ph \big), \ph \Big\rangle_{\H_{\tilde{\mu}}} ~ \int \limits_{U(n)} e^{-i(B v_r,z)_{\C^n}} dB \\[2pt]
\nn &=& \Big\langle \big(\tau_{\tilde{\mu}} \otimes \pi_{(0,r)}\big)(A,z,t) \big( \ph \otimes 1  \big),\ph \otimes 1  \Big\rangle_{\H_{\tilde{\mu}} \otimes \H_{(0,r)}},
\end{eqnarray}
where the measure $d\sigma(v) $ is the invariant measure on the complex
sphere $S^n $ in $\C^{n} $ defined in \mbox{Corollary 8.2}.

Choosing $\xi:= \ph \otimes 1 \in \H_{\tilde{\mu}} \otimes \H_{(0,r)}$, the claim is shown.
~\\
\qed

\begin{definition}\label{nulargerequalmu}
~\\
\rm   Let $\nu=(\nu_1,\cdots, \nu_n),\mu=(\mu_1,\cdots, \mu_n)\in P_n $ and write
\begin{eqnarray*}
 \mu\precsim\nu \text{ or }\nu\succsim\mu
 \end{eqnarray*}
if
\begin{eqnarray*}
 \nu_1\geq\mu_1\geq\nu_2\geq\mu_2\geq\cdots \geq\mu_{n-1}\geq\nu_n\geq\mu_n.
 \end{eqnarray*}
Let also for $\nu\in P_n $
\begin{eqnarray*}
 s(\nu):=\nu_1+\nu_2+\cdots +\nu_n.
 \end{eqnarray*}

 \end{definition}
 
\begin{lemma}\label{mulessnu}
~\\
Let $\ti\mu\in P_n $ and let $\ch_{\ti\nu} $ be a weight of $\T^n $ appearing in $\H_{\ti\mu} $ with $\ti\mu\ne\ti\nu $. Then
\begin{eqnarray*}
 \ti\nu\not\precsim\ti\mu
 \end{eqnarray*}
 and $\mu:=(\mu_1,\cdots, \mu_{n-1})\ne\nu:=(\nu_1,\cdots, \nu_{n-1}) $.
 \end{lemma}

Proof: \\
Since $\ti\nu $ is  a weight of $\H_{\ti\mu} $, one has
\begin{eqnarray*}
 \ti\nu=\ti\mu-\sum_{i=1}^{n-1}m_il_i,
 \end{eqnarray*}
where $m_i\in\N$ for $i \in \{1,\cdots, n-1\}$ and $l_i $ is the fundamental weight of $\t^n $ defined by
\begin{eqnarray*}
\langle l_i,X\rangle:=x_i-x_{i+1} \:\:\: \text{for}~ X=\sum \limits_{j=1}^{n-1} x_j E_{j,j}\in\t^n.
 \end{eqnarray*}
Since $\ti\mu\ne\ti\nu $, there exists $j\in\{1,\cdots, n\} $ such that
$\nu_1=\mu_1,\cdots, \nu_{j-1}=\nu_{j-1}$ and $\nu_j<\mu_j $. In particular,
$m_1=\cdots m_{j-1}=0$ and $m_j\ne0 $. But since $\sum \limits_{i=1}^n\mu_i=\sum \limits_{i=1}^n\nu_i $, one cannot have $\nu_i\leq \mu_i $ for all $i $. Hence, for some smallest $ k>j $, one has $\nu_k>\mu_k $. Necessarily $k<n $ since otherwise $m_i=0 $ for all $i\leq n-1 $ and then $\ti\mu=\ti\nu $. Therefore, $\nu\ne\mu $ and $\ti\nu\not\precsim\ti\mu $. \\
\qed
~\\
~\\

Proof of Theorem \ref{rueckweg}: \\
Without restriction, one can assume that the sequence $\big(\pi_{(\la^k,\alpha_k)}\big)_{k\in\mathbb{N}}$ fulfills Condition (i). The case of a sequence $\big(\pi_{(\la^k,\alpha_k)}\big)_{k\in\mathbb{N}}$ fulfilling the second condition is very similar. \\
~\\
For $ \tilde{n} \in \N$ and $\nu \in P_{\tilde{n}}$, let $\phi^{\nu}$ be the highest weight vector of $\tau_{\nu}$ in the Hilbert space $\H_{\nu}$.
Let $\tilde{\mu}:=(\mu_1,...,\mu_{n-1}, \mu_{n-1})$ and define the representation $\si_{(\tilde{\mu},r)} $ of $G_n $ by
\begin{eqnarray*}
\sigma_{(\tilde{\mu},r)}:=\tau_{\tilde{\mu}}\ot \pi_{(0,r)}.
 \end{eqnarray*}
The Hilbert space
$\H_{\sigma_{(\tilde{\mu},r)}} $ of the representation $\sigma_{(\tilde{\mu},r)}$ is the space
\begin{eqnarray*}
 \H_{\sigma_{(\tilde{\mu},r)}}=L^2 \big( U(n)/ U(n-1),\H_{\tilde{\mu}} \big)
 \end{eqnarray*}
 and $G_n $ acts on $\H_{\sigma_{(\tilde{\mu},r)}}$ by
\begin{eqnarray*}
\sigma_{(\tilde{\mu},r)}(A,z,t)(\xi)(B)=e^{-i( B v_r ,z)_{\C^n}} \tau_{\tilde{\mu}}(A)\big(\xi(A\inv B) \big) \:\:\: \forall A,B \in  U(n) ~\forall (z,t)\in \mathbb{H}_n~\forall \xi \in \H_{\sigma_{(\tilde{\mu},r)}}.
 \end{eqnarray*}
One decomposes the representation ${\tau_{\tilde{\mu}}}\res{U(n-1)}$ into the direct sum of irreducible representations of the group $ U(n-1) $ as follows:
\begin{eqnarray*}
 {\tau_{\tilde{\mu}}}\res{U(n-1)} =\sum_{\nu\in S(\tilde{\mu})}\rho_\nu,
 \end{eqnarray*}
where  $S(\tilde{\mu})$ denotes the support of  ${\tau_{\tilde{\mu}}}\res{U(n-1)} $ in $\wh{U(n-1)}$. Furthermore, let $p_\nu $ be the orthogonal projection of $\H_{\tilde{\mu}} $ onto its $ U(n-1) $-invariant component $\H_\nu $.\\
The representation $\rho_{\mu}$ is one of the representations appearing in this sum, since the highest weight vector $\phi^{\tilde{\mu}}$ of $\tau_{\tilde{\mu}} $ is also the highest weight vector of the representation $\rho_\mu $. \\


Defining for $\tilde{\nu} \in  P_n$ and the highest weight vector $\phi^{\tilde\nu}$ in $\H_{\tilde\nu} $ the function $c^{\tilde{\nu}}_{\et,\ph^{\tilde{\nu}}} $
by
\begin{eqnarray*}
 c^{\tilde{\nu}}_{\et,\ph^{\tilde{\nu}}}(A):=\big\langle \tau_{\tilde{\nu}}(A\inv)\et,\ph^{\tilde{\nu}} \big\rangle_{\H_{\tilde{\nu}}} \:\:\: \forall A \in  U(n)~\forall \et\in \H_{\tilde{\nu}},
 \end{eqnarray*}
one can identify for any $\tau_{\tilde{\nu}}\in \wh{U(n)} $ the Hilbert space $\H_{\tilde{\nu}} $ with the subspace
$L^2_{\tilde{\nu}} $ of $L^2 \big(U(n) \big) $ given by
\begin{eqnarray*}
 L^2_{\tilde{\nu}}=\big\{c^{\tilde{\nu}}_{\et,\ph^{\tilde{\nu}}}|~\et\in\H_{\tilde{\nu}} \big\}.
 \end{eqnarray*}
Now, it will be shown that
\begin{eqnarray}\label{decomtimur}
 \sigma_{(\tilde{\mu},r)}\cong \sum_{\nu\in S(\tilde{\mu})}\pi_{(\nu,r)}.
 \end{eqnarray}
In particular, one then gets
\begin{eqnarray*}
 \H_{\sigma_{(\tilde{\mu},r)}} \cong \sum_{\nu\in S(\tilde{\mu})}L^2 \big( U(n)/ U(n-1),\rho_{\nu}\big).
 \end{eqnarray*}
For $\nu \in S(\tilde\mu)$, $\xi \in  L^{2} (U (n)/U (n- 1), \H_{\ti\mu})$, $A \in U (n) $ and $A' \in  U (n - 1)$, let
\begin{eqnarray*}
U_{\tilde{\mu}}^\nu(\xi)(A)(A'):=\sqrt{d_\nu} \big\langle\tau_{\tilde{\mu}}(A\inv)\xi(A),\rho_{\nu}(A')\ph^{\nu} \big\rangle_{\H_{\nu}}.
 \end{eqnarray*}
Moreover, 
\begin{eqnarray*}
\big\langle\tau_{\tilde{\mu}}(A\inv)\xi(A),\rho_{\nu}(A')\ph^{\nu} \big\rangle_{\H_{\nu}}&=&\Big\langle\tau_{\tilde{\mu}}(A\inv)\xi(A),p_{\nu} \big(\rho_{\nu}(A')\ph^{\nu} \big) \Big\rangle_{\H_{\nu}}\\
&=&\Big\langle p_{\nu} \big(\tau_{\tilde{\mu}}(A\inv)\xi(A) \big),\rho_{\nu}(A')\ph^{\nu} \Big\rangle_{\H_{\nu}},
\end{eqnarray*}
i.e. one has a scalar product on the space $\H_{\nu}$. Furthermore, for all
 $\xi \in  L^{2} (U (n)/U (n - 1), \H_{\ti\mu})$, $A \in U (n) $ and $A' \in  U (n - 1)$,
\begin{eqnarray*}
 U_{\tilde{\mu}}^{\nu}\xi(AB')(A')&=&\sqrt{d_\nu} \big\langle \tau_{\tilde{\mu}} \big({B'}\inv A\inv \big)\xi(A),\rho_{\nu}(A')\ph^{\nu} \big\rangle_{\H_{\nu}}\\
\nn  &= &
\sqrt{d_\nu} \big\langle\tau_{\tilde{\mu}}( A\inv)\xi(A),\rho_{\nu}(B'A')\ph^{\nu} \big\rangle_{\H_{\nu}}\\[1pt]
\nn  &= &
\rho_{\nu}(B')\inv U_{\tilde{\mu}}^{\nu}\xi(A)(A').
 \end{eqnarray*}
Hence, each vector  $U_{\tilde{\mu}}^{\nu}(\xi) $ fulfills the covariance condition of the space $L^2 \big( U(n)/ U(n-1),\rho_\nu \big)$. Therefore, $U_{\tilde\mu}(\xi):=\sum \limits_{\nu\in S(\tilde\mu)}U_{\tilde\mu}^{\nu}(\xi) $ is an element of the space $\sum \limits_{\nu\in S(\tilde\mu)}L^2 \big( U(n)/ U(n-1),\rho_\nu \big) $:
\begin{eqnarray*}
 U_{\tilde\mu}:L^2 \big( U(n)/ U(n-1),\H_{\tilde\mu}\big)\to \sum_{\nu\in S(\tilde\mu)}L^2 \big( U(n)/ U(n-1),\rho_\nu \big).
 \end{eqnarray*}
 Furthermore,

\begin{eqnarray*}
\|U_{\tilde{\mu}}\xi \|_2^2\nn
\nn&= & \sum_{\nu \in S(\tilde{\mu})}~ \int \limits_{ U(n)}\|U_{\tilde{\mu}}^\nu(\xi)(A)\|_{\nu}^2 dA\\[1pt]
\nn  &= &
\sum_{\nu \in S(\tilde{\mu})}~ \int \limits_ {U(n)} \int \limits_{U(n-1)} \sqrt{d_\nu}\big\langle\tau_{\tilde{\mu}}(A\inv)\xi(A),\rho_{\nu}(A')\ph^{\nu}\big\rangle_{\H_{\nu}} \\[-2pt]
&&~~~~~~~~~~~~~~~~~~~~~~~~\ol{ \sqrt{d_{\nu}} \big\langle\tau_{\tilde{\mu}}( A\inv)\xi(A),\rho_{\nu}(A')\ph^{\nu} \big\rangle_{\H_{\nu}}}dA' dA\\[2pt]
\nn  &= &
\sum_{\nu \in S(\tilde{\mu})}~ \int \limits_{U(n)} \int \limits_{U(n-1)} d_\nu \big\langle p_{\nu}\big(\tau_{\tilde{\mu}}(A\inv)\xi(A) \big),\rho_{\nu}(A')\ph^{\nu} \big\rangle_{\H_{\nu}} \\[-2pt]
&&~~~~~~~~~~~~~~~~~~~~~~~~ \ol{ \big\langle p_{\nu} \big(\tau_{\tilde{\mu}}( A\inv)\xi(A) \big),\rho_{\nu}(A')\ph^{\nu} \big\rangle_{\H_{\nu}}}dA' dA\\[3pt]
\nn  &= &
\sum_{\nu \in S(\tilde{\mu})}~ \int \limits_{U(n)} \int \limits_{U(n-1)} d_\nu \Big|\big\langle p_{\nu}\big(\tau_{\tilde{\mu}}(A\inv)\xi(A) \big),\rho_{\nu}(A')\ph^{\nu} \big\rangle_{\H_{\nu}} \Big|^2dA' dA\\[1pt]
\end{eqnarray*}
\begin{eqnarray*}
\nn  & =&
\sum_{\nu \in S(\tilde{\mu})}~\int \limits_{ U(n)} \big\|p_\nu\big(\tau_{\tilde{\mu}}(A\inv)\xi(A)\big) \big\|_{\nu}^2~ \big| \phi^{\nu} \big|^2 dA \:\:\:\:\:\:\:\:\:\:\:\:\:\:\:\:\:\:\:\:\:\:\:\:\:\: \\[1pt]
\nn  & =&
\sum_{\nu \in S(\tilde{\mu})}~\int \limits_{ U(n)} \big\|p_\nu\big(\tau_{\tilde{\mu}}(A\inv)\xi(A)\big) \big\|_{\nu}^2 dA\\[1pt]
\nn  & =&
\int \limits_{ U(n)} \bigg\| \sum_{\nu \in S(\tilde{\mu})} p_\nu\big(\tau_{\tilde{\mu}}(A\inv)\xi(A)\big) \bigg\|_{\tilde{\mu}}^2 dA\\[2pt]
\nn  & =&
 \big\| \tau_{\tilde{\mu}}(\cdot \inv)\xi( \cdot) \big\|_{2}^2
~= ~\|\xi \|_2^2.
 \end{eqnarray*}
Moreover, for all $(A,z,t) \in G_n$, all $\xi \in L^2 \big( U(n)/ U(n-1),\H_{\tilde{\mu}} \big)$, all $B\in U(n)$ and all $A' \in U(n-1)$, one gets
\begin{eqnarray*}
\sum_{\nu\in S(\tilde{\mu})}\pi_{(\nu,r)}(A,z,t)(U_{\tilde{\mu}}\xi)(B)(A')
 \nn  & = &
e^{-i( B v_r ,z)_{\C^n} }\sum_{\nu\in S(\tilde{\mu})}(U_{\tilde{\mu}}^\nu\xi)(A\inv B)(A')\\[1pt]
\nn  & = &
e^{-i( B v_r ,z)_{\C^n} }\sum_{\nu\in S(\tilde{\mu})}\sqrt{d_\nu}\big\langle\tau_{\tilde{\mu}}(B\inv A)\xi(A\inv B),\rho_{\nu}(A')\ph^{\nu} \big\rangle_{\H_{\nu}}\\[1pt]
\nn  & = &
e^{-i( B v_r ,z)_{\C^n} }\sum_{\nu\in S(\tilde{\mu})}\sqrt{d_\nu} \Big\langle\tau_{\tilde{\mu}}(B\inv) \big(\tau_{\tilde{\mu}}(A)\xi(A\inv B) \big),\rho_{\nu}(A')\ph^{\nu} \Big\rangle_{\H_{\nu}}\\ [1pt]
\nn  & = &
e^{-i( B v_r ,z)_{\C^n} }\sum_{\nu\in S(\tilde{\mu})} U_{\tilde{\mu}}^{\nu} \Big( \tau_{\tilde{\mu}}(A) \xi \big( A^{-1} \cdot \big) \Big)(B)(A')\\ [1pt]
\nn  & = &
U_{\tilde{\mu}} \big(\tau_{\tilde{\mu}}\otimes\pi_{(0,r)}(A,z,t)\xi\big)(B)(A')
~=~ U_{\tilde{\mu}} \big(\sigma_{(\tilde{\mu},r)}(A,z,t)\xi\big)(B)(A').
 \end{eqnarray*}
Therefore, (\ref{decomtimur}) holds. \\

For  $\xi=\phi^{\tilde{\mu}} \otimes 1\in L^2 \big(U(n)/ U(n-1),\H_{\tilde{\mu}} \big) $, such that $\xi(A)=\phi^{\tilde \mu} $ for all $A$, one has for all $B \in U(n)$ and all $A' \in U(n-1)$,
\begin{eqnarray*}
 U_{\tilde{\mu}}\xi(B)(A')=\sum_{\nu\in S(\tilde{\mu})}\sqrt{d_\nu} \big\langle \tau_{\tilde{\mu}}(B\inv)\ph^{\tilde{\mu}} ,\rho_\nu(A')\ph^\nu \big\rangle_{\H_{\nu}}=:\sum_{\nu\in S(\tilde{\mu})} \Phi_{\tilde{\mu}}^\nu(B)(A').
 \end{eqnarray*}
In particular, since $\phi^{\mu}=\phi^{\tilde{\mu}}$,
$$ \Phi_{\tilde{\mu}}^\mu(B)(I) =\sqrt{d_\mu} \big\langle \tau_{\tilde{\mu}}(B\inv)\ph^{\tilde{\mu}} ,\ph^\mu \big\rangle_{\H_{\mu}}
= \sqrt{d_\mu} \big\langle \ph^{\tilde{\mu}} ,\tau_{\tilde{\mu}}(B)\ph^{\tilde{\mu}} \big\rangle_{\H_{\tilde{\mu}}} \ne  0.$$
From Theorem \ref{rnenull} follows that the subset $\big\{\pi_{(\nu,r)}|~ \mu\ne\nu \in P_{n-1} \big\} $ is closed in $\wh{G_n} $. Hence, there exists $F_{\mu} =(F_\mu)^*$ of norm 1 in  $C^*(G_n) $ whose Fourier transform at
$\pi_{(\nu,r)} $ is 0 if $\mu\ne\nu \in P_{n-1} $ and for which
\begin{eqnarray*}
 \pi_{(\mu,r)}(F_{\mu})=:P_{\Phi_{\tilde{\mu}}^\mu}
 \end{eqnarray*}
is the orthogonal projection onto the space $\C \Phi_{\tilde{\mu}}^\mu\subset \H_{(\mu,r)} $.  In particular,
\begin{eqnarray}\label{phi}
\nn U_{\tilde{\mu}} \big(\sigma_{(\tilde{\mu},r)}(F_\mu)(\phi^{\tilde{\mu}} \otimes 1) \big)
&=& \sum_{\nu\in S(\tilde{\mu})} \pi_{(\nu,r)}(F_\mu)\big( U_{\tilde{\mu}}(\phi^{\tilde{\mu}} \otimes 1) \big)
~=~ \pi_{(\mu,r)}(F_\mu)\big( U_{\tilde{\mu}}(\phi^{\tilde{\mu}} \otimes 1) \big) \\
&=& P_{\Phi_{\tilde{\mu}}^\mu} \big( U_{\tilde{\mu}}(\phi^{\tilde{\mu}} \otimes 1) \big)~=~  \Phi_{\tilde{\mu}}^\mu.
 \end{eqnarray}
 \\
Define the coefficient $c_{\mu} $ of $\l1{G_n} $ by
\begin{eqnarray}\label{cucom}
 \nn  c_\mu(F)& := &
\big\langle \sigma_{(\tilde{\mu},r)}(F_\mu\ast F\ast F_\mu)(\phi^{\tilde{\mu}} \otimes 1) ,\phi^{\tilde{\mu}} \otimes 1 \big\rangle_{\H_{\sigma_{(\tilde{\mu},r)}}} \\[1pt]
\nn&=&
\big\langle \sigma_{(\tilde{\mu},r)}(F \ast F_\mu)(\phi^{\tilde{\mu}} \otimes 1) , \sigma_{(\tilde{\mu},r)}(F_\mu)(\phi^{\tilde{\mu}} \otimes 1) \big\rangle_{\H_{\sigma_{(\tilde{\mu},r)}}} \\[1pt]
\nn&=&
\Big\langle  U_{\tilde{\mu}} \big( \sigma_{(\tilde{\mu},r)}(F \ast F_\mu)(\phi^{\tilde{\mu}} \otimes 1)\big) , U_{\tilde{\mu}} \big( \sigma_{(\tilde{\mu},r)}(F_\mu)(\phi^{\tilde{\mu}} \otimes 1) \big) \Big\rangle_{\H_{(\mu,r)}} \\[1pt]
\nn&\overset{(\ref{phi})}{=}&
\bigg\langle \sum_{\nu\in S(\tilde{\mu})} \pi_{(\nu,r)} (F) \circ \pi_{(\nu,r)} ( F_\mu) \big(U_{\tilde{\mu}}(\phi^{\tilde{\mu}} \otimes 1) \big), \Phi_{\tilde{\mu}}^{\mu} \bigg\rangle_{\H_{(\mu,r)}} \\[1pt]
\nn&=&
\big\langle \pi_{(\mu,r)} (F) \circ \pi_{(\mu,r)} ( F_\mu) \big(U_{\tilde{\mu}}(\phi^{\tilde{\mu}} \otimes 1) \big),\Phi_{\tilde{\mu}}^\mu \big\rangle_{\H_{(\mu,r)}} \\[1pt]
&\overset{(\ref{phi})}{=} &
\big\langle \pi_{(\mu,r)}(F)\Phi_{\tilde{\mu}}^\mu ,\Phi_{\tilde{\mu}}^\mu \big\rangle_{\H_{(\mu,r)}}
 \end{eqnarray}
for all $F \in L^1(G_n)$. \\
Let $X(\tilde{\mu}, \ol{\al_k}) $ be the collection of all $\tilde{\nu} = (\nu_1,..., \nu_n)\in  P_n$ such that $\ch_{\tilde{\nu}} $ is a character of $\T_n $ appearing in $\H_{\tilde{\mu}} $ and such that $\tau_{\tilde{\nu}_k}$ is contained in the representation $\tau_{\tilde{\mu}}\ot \ol{W_{\al_k}}$ for $\tilde{\nu}_k:=(\nu_1,\nu_2,..., \nu_{n-1},\nu_{n} + \la_n^k)$. \\
Then, for $\pi^{(\tilde{\mu},\al_k)}$ defined as in Proposition \ref{prep-conj}, by [\ref{knapp}], Chapter IV.11,
\begin{eqnarray*}\label{mukalkdecom}
 \pi^{(\tilde{\mu},\al_k)}=\sum_{\tilde{\nu}\in X(\tilde{\mu}, \ol{\al_k})}\pi_{(\tilde{\nu}_k,\al_k)}.
\end{eqnarray*}

Furthermore, decompose the vector
$$\xi_k= \phi^{\tilde{\mu}} \otimes \bigg(\frac{1}{{R_k}^{\frac{1}{2}}}~\sum \limits_{\substack{q\in \N^n:\\ |q|=N_k}}\ol{h_{q,\al_k}} \otimes h_{q,\al_k} \bigg)$$
for every $k \in \N$ into the orthogonal sum
\begin{eqnarray*}
 \xi_k=\sum_{\tilde{\nu}\in X(\tilde{\mu},\ol{\al_k})}\xi_k^{\tilde{\nu}}
 \end{eqnarray*}
for $\xi_k^{\tilde{\nu}}\in\H_{(\tilde{\nu}_k,\al_k)} $.
This gives a decomposition
\begin{eqnarray*}\label{decomlaal}
 {\big\langle \pi^{(\tilde{\mu},\al_k)}(\cdot) \xi_k, \xi_k \big \rangle_{\H_{\tilde{\mu}} \otimes \ol{\F_{\al_k}(n)} \otimes \F_{\al_k}(n)}} = c^{\pi^{(\tilde{\mu},\al_k)}}_{\xi_k}=\sum_{\tilde{\nu}\in X(\tilde{\mu}, \ol{\al_k})}c^{\tilde{\nu}}_{\xi_k^{\tilde{\nu}}}.
 \end{eqnarray*}
 \\
Let $c_{\xi^{\tilde{\nu}}}^{\tilde\nu} $ be the weak$^* $-limit of a subsequence of $\big(c_{\xi_k^{\tilde{\nu}}}^{\tilde\nu} \big)_{k \in \N} $ and let for   $c_{\xi^{\tilde{\nu}}}^{\tilde\nu} \ne 0 $ the representation $\pi\in \wh{G_n} $ be an element of the support of $c_{\xi^{\tilde{\nu}}}^{\tilde\nu} $. From Theorem \ref{unconrep} follows that $\pi=\lim \limits_{k \to \iy} \pi_{(\tilde{\nu}_k,\al_k)}=\pi_{(\nu,r)}$ for $\nu=(\nu_1,..., \nu_{n-1})$.
Furthermore, one observes that for $\tilde{\mu} \ne \tilde{\nu}:=(\nu_1,..., \nu_{n-1}, \nu_{n}) \in X(\tilde{\mu},\ol{\al_k})$, one has $\nu\ne\mu $. Hence, $\pi_{(\nu,r)}(F_\mu)=0 $.
Thus,
\begin{eqnarray*}
\lim \limits_{k \to \iy} \big\langle\pi_{(\tilde{\nu}_k,\al_k)} (F_\mu)\xi_k^{\tilde{\nu}}, \pi_{(\tilde{\nu}_k,\al_k)} (F_\mu)\xi_k^{\tilde{\nu}} \big\rangle_{\H_{(\tilde{\nu}_k,\al_k)}}
=\big\langle \pi_{(\nu,r)}(F_\mu)\xi^{\tilde{\nu}},\pi_{(\nu,r)}(F_\mu)\xi^{\tilde{\nu}} \big\rangle_{\H_{(\nu,r)}}=0
 \end{eqnarray*}
 and therefore,
 \begin{eqnarray}\label{fornuzero}
\lim \limits_{k \to \iy} \big\langle \pi_{(\tilde{\nu}_k,\al_k)}(F)\circ \pi_{(\tilde{\nu}_k,\al_k)} (F_\mu)\xi_k^{\tilde{\nu}}, \pi_{(\tilde{\nu}_k,\al_k)} (F_\mu)\xi_k^{\tilde{\nu}} \big\rangle_{\H_{(\tilde{\nu}_k,\al_k)}}
=0 \:\:\: \forall F\in C^*(G_n).
 \end{eqnarray}
Now, $\tilde{\mu}_k= (\mu_1,...,\mu_{n-1},\mu_{n-1}+ \la_n^k)$ and $\la^k=(\mu_1,...,\mu_{n-1}, \la_n^k)$. Since $\la_n^k \overset{k \to \iy}{\longrightarrow} - \iy$,  their behavior for $k \to \iy$ is the same. Hence, by Proposition \ref{prep-conj} and its proof, for all $F \in C^*(G_n)$,
\begin{eqnarray*}
\big\langle \pi_{(\mu,r)}(F)\Phi_{\tilde{\mu}}^{\mu},\Phi_{\tilde{\mu}}^{\mu} \big\rangle_{\H_{(\mu,r)}}
&=&c_\mu(F_\mu)\\
\nn  &\overset{(\ref{cucom})}=&
\big\langle \tau_{\tilde{\mu}}\ot\pi_{(0,r)}(F_\mu\ast F\ast F_\mu)(\phi^{\tilde{\mu}} \otimes 1) ,\phi^{\tilde{\mu}} \otimes 1 \big\rangle_{\H_{\sigma_{(\tilde{\mu},r)}}} \\[1pt]
\nn  &\overset{Proposition}{\underset{\ref{prep-conj}}{=}} &
\lim \limits_{k \to \iy} \big\langle \pi^{(\tilde{\mu},\al_k)} (F_\mu\ast F\ast F_\mu)\xi_k,\xi_k \big\rangle_{\H_{\tilde{\mu}} \otimes \ol{\F_{\al_k}(n)} \otimes \F_{\al_k}(n)}\\[1pt]
\nn  &= &
\lim \limits_{k \to \iy} ~
 \sum_{\tilde{\nu}\in X(\tilde{\mu}, \ol{\al_k})}\big\langle\pi_{(\tilde{\nu}_k,\al_k)}(F_\mu\ast F\ast F_\mu)\xi_k^{\tilde{\nu}},\xi_k^{\tilde{\nu}}\big\rangle_{\H_{(\tilde{\nu}_k,\al_k)}}\\
\nn  &\overset{(\ref{fornuzero})}= &
\lim \limits_{k \to \iy} \big\langle \pi_{(\tilde{\mu}_k,\al_k)}(F_\mu\ast F\ast F_\mu)\xi_k^{\tilde{\mu}},\xi_k^{\tilde{\mu}}\big\rangle_{\H_{(\tilde{\mu}_k,\al_k)}}+0\nn\\[1pt]
\nn  &= &
\lim \limits_{k \to \iy} \Big\langle \pi_{(\tilde{\mu}_k,\al_k)}( F) \big(\pi_{(\tilde{\mu}_k,\al_k)}(F_\mu)\xi_k^{\tilde{\mu}} \big),\pi_{(\tilde{\mu}_k,\al_k)}(F_\mu)\xi_k^{\tilde{\mu}} \Big\rangle_{\H_{(\tilde{\mu}_k,\al_k)}} \\[1pt]
\nn  &= &
{\lim \limits_{k \to \iy} \Big\langle \pi_{(\la^k,\al_k)}( F) \big(\pi_{(\tilde{\mu}_k,\al_k)}(F_\mu)\xi_k^{\tilde{\mu}} \big),\pi_{(\tilde{\mu}_k,\al_k)}(F_\mu)\xi_k^{\tilde{\mu}} \Big\rangle_{\H_{(\la^k,\al_k)}}}.
 \end{eqnarray*}
Choosing $\tilde{\xi}:=\Phi_{\tilde{\mu}}^{\mu}$ and {$\tilde{\xi}_k:=\pi_{(\tilde{\mu}_k,\al_k)}(F_\mu)\xi_k^{\tilde{\mu}}$}, one has for any $F\in C^*(G_n) $
\begin{eqnarray}\label{lim xi mu k}
  \lim \limits_{k \to \iy} \big\langle \pi_{(\la^k,\al_k)}( F) \tilde{\xi}_k,\tilde{\xi}_k \big\rangle_{\H_{(\la^k,\al_k)}}=  \big\langle \pi_{(\mu,r)}(F)\tilde{\xi},\tilde{\xi} \big\rangle_{\H_{(\mu,r)}}
\end{eqnarray}

and hence,
\begin{eqnarray*}
 \pi_{(\mu,r)}=\lim \limits_{k \to \iy} \pi_{(\la^k,\al_k)}.
 \end{eqnarray*}
 \qed
~\\

\begin{remark}\label{xi mu k}
~\\
\rm   It follows from (\ref{lim xi mu k})
that
\begin{eqnarray*}
 \limk \big\|\xi^{\ti\mu}_k \big\|_{\H_{(\tilde{\mu}_k,\al_k)}}=1.
 \end{eqnarray*}

 \end{remark}

 \begin{theorem}\label{ris0}
 ~\\
Let $\ta_{\ti\rh}\in\widehat{U(n)}$.

If $\underset{k\to\infty}{\lim}\alpha_k=0$ and the sequence $\big(\pi_{(\ti{\la}^k,\alpha_k)}\big)_{k\in\mathbb{N}}$ of
elements of $\wh{G_n}$ satisfies one
of the following conditions:
\begin{enumerate}[(i)]
\item for $k$ large enough, $\alpha_k>0$, $\rh_1\geq\la_1^k\geq...\geq\rh_{n-1}\geq\la_{n-1}^k\geq\rh_n\geq\la_n^k$
and $\underset{k\to\infty}{\lim}\alpha_k\la^k_{n}=0$,
\item for $k$ large enough, $\alpha_k<0$,
$\la_1^k\geq\rh_1\geq\la_2^k\geq\rh_2\geq...\geq\la_{n}^k\geq\rh_n$ and $\underset{k\to\infty}{\lim}\alpha_k\la^k_{1}=0$,
\end{enumerate}
then the sequence $\big(\pi_{(\ti\la^k,\alpha_k)}\big)_{k\in\mathbb{N}}$ converges to the representation $\ta_{\ti\rh}$ in
$\wh{G_n}$. \\
\end{theorem}

Proof: \\
Again, only consider the case $\al_k>0$ for all $k\in\N $. \\

Let $\ti\rh=(\rh_1,\cdots,\rh_n )\in P_n$ satisfy the conditions of the theorem, i.e.
\mbox{$\rh_1\geq\la_1^k\geq...\geq\rh_{n-1}\geq\la_{n-1}^k\geq\rh_n\geq\la_n^k$}. Passing to a subsequence, one can assume that $\mu_1:=\la_1^k,\cdots, \mu_{n-1}:=\la^k_{n-1} $ for all $k\in\N $.

Let
\begin{eqnarray*}
\mu:=(\mu_1,\cdots, \mu_{n-1},\mu_{n-1}),\:\:\:  \ti\mu_k=(\mu_1,\cdots, \mu_{n-1},\la_n^k) \:\:\: \text{and} \:\:\: N_k:=\mu_{{n-1}}-\la_n^{k} \:\:\: \forall k\in\N.
 \end{eqnarray*}
Then
\begin{eqnarray*}
 \ti\rh=\ti\mu+r=\ti\mu_k+r_{k}
 \end{eqnarray*}
for some $r=(r_1,\cdots, r_n)\in \N^n $ and $r_k=r+(0,\cdots, 0, N_k) $. Let
\begin{eqnarray*}
 m:=\sum_{i=1}^n r_i.
 \end{eqnarray*}
Hence, by Pieri's rule, one obtains
\begin{eqnarray*}
 \ta_{\rh}\in \ta_{\ti\mu_k}\ot \ta_{N_k+m} \in \pi_{(\mu_k,\al_k)}{\res{U(n)}}
 \end{eqnarray*}
      for $k $ large enough.
We take the highest weight vector $\ph ^{\ti\rh}_k $ of the representation $\ta_{\ti\rh} $ considered as a subrepresentation of $U(n) $ on the Hilbert space $\H_{\ti\mu_k}\ot \P_{N_k+m,\al_k}(n) $, where $\P_{N_k+m,\al_k}(n) $ is the space of all polynomials of degree $N_k+m$ in the Fock space $\F_{\al_k}(n)$.

Recall also that the polynomials
\begin{eqnarray*}
h_{q,\al}(w)&:=&
{\bigg(\frac{ \al }{2 \pi} \bigg)^{\frac{n}{2}} \sqrt{\frac{\al^{\val q}}{2^{\val{q}}q!}}
}~w^q \:\:\: \forall w\in\C^n \:\:\: \text{for}~ \val q=N_k+m
\end{eqnarray*}
form a Hilbert space  basis of $\P_{N_k+m,\al_k}(n)  $.
Hence, one can write
\begin{eqnarray*}
 \ph ^{\ti\rh}_k=\sum \limits_{\substack{q\in \N^n:\\ |q|=N_k+m}}a^k_q\ot h_{q,\al_k},
 \end{eqnarray*}
where for any $q $ and $k $, the vector $a^k_q $ is contained in the $\T^n $-eigenspace of $\H_{\ti\mu_k} $ for the weight $\ch_{\ti\rh-q} $, since $h_{q,\al_k} $ is in the eigenspace for the weight $\ch_q $. In particular, one has
\begin{eqnarray*}
 \big\langle a^k_q,a^k_{q'} \big\rangle_{\H_{\tilde\mu_k}}=0 \:\:\: \text{for}~q\ne q'.
 \end{eqnarray*}

Let $(\xi_k)_{k \in \N}= \Bigg(\sum \limits_{\substack{q\in \N^n:\\ |q|=N_k}} a_q^k\ot h_{q, \al_k} \Bigg)_{k \in \N}$  be a sequence of vectors of length 1 in $ \big(\H_{\ti\mu_k}\ot \P_{N_k}\big)_{k \in \N}$ such that 
\begin{eqnarray*}
 \big\langle a^k_q,a^k_{q'} \big\rangle_{\H_{\tilde\mu_k}}=0 \:\:\: \text{for}~q\ne q'.
 \end{eqnarray*}
 It will now be shown that 
\begin{eqnarray}\label{pi ti mu z is almost one }
 \limk \big\|\pi_{(\ti\mu_k,\al_k)}(z,t)(\xi_k)-\xi_k \big\|_{\H_{(\tilde\mu_k,\al_k)}}=0
 \end{eqnarray}
 uniformly on compacta in $(z,t)\in \HH_n $: \\
~\\

Since $\sum \limits_{\substack{q\in \N^n:\\ |q|=N_k}} \big\|a_q^k \big\|^2_{\H_{\tilde\mu_k}}=\|\xi_k\|_{\H_{(\tilde\mu_k,\al_k)}}^2=1 $ for all $k\in\N $, it suffices to prove that
\begin{eqnarray*}
 \limk \big\langle \pi_{(\ti\mu_k,\al_k)}(z,t)(\xi_k)-\xi_k,\xi_k \big\rangle_{\H_{(\tilde\mu_k,\al_k)}}=0
 \end{eqnarray*}
uniformly on compacta.

One gets
\begin{eqnarray*}
\big\langle \pi_{(\ti\mu_k,\al_k)}(z,t)(\xi_k)-\xi_k,\xi_k \big\rangle_{\H_{(\tilde\mu_k,\al_k)}}
\nn  & = &
\sum \limits_{\substack{q\in \N^n:\\ |q|=N_k}} \sum \limits_{\substack{q'\in \N^n:\\ |q'|=N_k}}
 \big\langle a^k_q,a^k_{q'} \big\rangle_{\H_{\tilde\mu_k}} \big\langle \tau_{\al_k}(z,t)(h_{q,\al_k}),
 h_{q',\al_k} \big\rangle_{\H_{\al_k}}-1\\
 \nn  & = &
\sum \limits_{\substack{q\in \N^n:\\ |q|=N_k}}
 \big\|a^k_q\big\|^2_{\H_{\tilde\mu_k}} \big\langle \tau_{\al_k}(z,t)(h_{q,\al_k}),
 h_{q',\al_k} \big\rangle_{\H_{\al_k}}-1\\
 \nn  & = &
\sum \limits_{\substack{q\in \N^n:\\ |q|=N_k}}
 \big\|a^k_q \big\|^2_{\H_{\tilde\mu_k}} \Big( \big\langle \tau_{\al_k}(z,t)(h_{q,\al_k}),
 h_{q',\al_k} \big\rangle_{\H_{\al_k}} -1 \Big),
 \end{eqnarray*}
 since for $q\ne q' $
\begin{eqnarray*}
 \big\langle a^k_q,a^k_{q'} \big\rangle_{\H_{\tilde\mu_k}}=0.
 \end{eqnarray*}

Now, for any $q\in \N^n$ and any $k\in\N $, by (\ref{computations for q}), 
\begin{eqnarray*}
\big\langle \tau_{\al_k}(z,t)h_{q,\al_k},h_{q,\al_k} \big\rangle_{\H_{\al_k}}
 = 
e^{i \al_k t - \frac{\al_k}{4} |z|^2}
\sum \limits_{\substack{j:=(j_1,...,j_n) \in \N^n:\\ j_1 \leq q_1,...,j_n \leq q_n}} \Big(\frac{\al_k}{2} \Big)^{|j|}
\frac{q!}{(q-j)!} \frac{z^j(-\ol{z})^{j}}{\big(j!\big)^2}
 \end{eqnarray*}
and thus, for $k $ large enough and $(z,t) $ in some compact set (i.e. $\val{\frac{\al_k N_k}{2}}<e^{-\val z^2}$), one has
\begin{eqnarray*}
\Big|\big\langle \pi^{(\ti\mu,\al_k)}(z,t)h_{q,\al_k},h_{q,\al_k} \big\rangle_{\H_{\al_k}}-1 \Big|
&=&
\Bigg|{e^{i \al_k t - \frac{\al_k}{4} |z|^2}}{}~ \sum \limits_{\substack{q_1 \in \N_{\geq j_1},...,q_n \in \N_{\geq j_n}:\\ q_1+...+q_n=N_k}} \Big(\frac{\al_k N_k}{2} \Big)^{j_1+...+j_n}~ \frac{q_1 (q_1-1) \cdots (q_1-j_1+1)}{{N_k}^{j_1}} \\[2pt]
&&\:\:\:\: \cdots \frac{q_n (q_n-1) \cdots (q_n-j_n+1)}{{N_k}^{j_n}}  \frac{z^j(-\ol{z})^{j}}{\big(j!\big)^2}-1 \Bigg| \\[4pt]
&\leq&
\big|
e^{i\al_k t - \frac{\al_k}{4} |z|^2}-1 \big|
+ \val{\al_k N_k }e^{\val z^2}.
\end{eqnarray*}

 Therefore,
\begin{eqnarray*}
 \nn  \Big| \big\langle \pi_{(\ti\mu_k,\al_k)}(z,t)(\xi_k)-\xi_k,\xi_k \big\rangle_{\H_{(\tilde\mu_k,\al_k)}}\Big|
\nn  & \leq &
\sum \limits_{\substack{q\in \N^n:\\ |q|=N_k}}
 \big\|a_q^k\big\|^2_{\H_{\tilde\mu_k}} \Big( \big|
e^{i\al_k t - \frac{\al_k}{4} |z|^2}-1 \big|
+ \big|\al_k N_k e^{\val z^2}\big| \Big)\\
\nn  & = &
\big|
e^{i\al_k t - \frac{\al_k}{4} |z|^2}-1 \big|
+ \big|\al_k N_k e^{\val z^2} \big|.
\end{eqnarray*}

This proves Claim (\ref{pi ti mu z is almost one }). \\

To finish the proof of Theorem \ref{ris0},
by (\ref{pi ti mu z is almost one }) for any $(A,z,t)\in G_n $, one has uniformly on compacta
\begin{eqnarray*}
 \nn
\limk \big\langle \pi_{(\la^k,\al_k)}(A,z,t)\xi^{\ti\rh}_k,\xi^{\ti\rh}_k \big\rangle_{\H_{(\la^k,\al_k)}}\nn  &=  &
\limk \big\langle \pi_{(\la^k,\al_k)}(A,0,0)\xi^{\ti\rh}_k,\xi^{\ti\rh}_k \big\rangle_{\H_{\tilde\rh}}\\
\nn  & = &
  \limk \big\langle \ta_{\ti\rh}(A)\ph ^{\ti\rh},\ph^{\ti\rh} \big\rangle_{\H_{\tilde\rh}}\\
\nn  & = &
\big\langle \ta_{\ti\rh}(A)\ph ^{\ti\rh},\ph^{\ti\rh} \big\rangle_{\H_{\tilde\rh}}
\end{eqnarray*}
for the highest weight vector $\ph ^{\ti\rh } $ of the representation $\ta_{\ti\rh} $.
This shows that
\begin{eqnarray*}
 \limk \pi_{(\la^k,\al_k)}=\ta_{\ti\rh}.
 \end{eqnarray*}
~\\
\qed

\section{The final result}
Together with the Theorems \ref{rnenull}, \ref{conv-la-al}, \ref{unconrep}, \ref{rueckweg} and \ref{ris0} and the result in Subsection \ref{tau-la}, we obtain the final result below:

\begin{theorem}\label{homeomorphism}
~\\
For general $n \in \N^*$, the  spectrum of the group $G_n=U(n)\ltimes\HH_n$ is homeomorphic to its space of admissible coadjoint orbits
$\g_n^\ddagger/G_n$.
\end{theorem}

\section{Appendix}

\begin{lemma}\label{det}
~\\
Let $n \in \N^*$, let $B_{\R}^{2n}$ be the $2n$-dimensional real unit ball and define the mapping
\begin{eqnarray*}
&&\psi:[0,1]^{n-1} \times [0,2 \pi)^{n} \times (0,1] \to B_{\R}^{2n}, \\[2pt]
&&\psi(s_1,...,s_{n-1},t_1,...,t_n,\rho):=\\[3pt]
&&\rho \Big( \sqrt{s_1} \cos (t_1), \sqrt{s_1} \sin (t_1),...,\sqrt{s_{n-1}} \cos (t_{n-1}),\sqrt{s_{n-1}} \sin (t_{n-1}), \\[1pt]
&&~~~~ \sqrt{1-s} \cos (t_n), \sqrt{1-s} \cos (t_n) \Big),
\end{eqnarray*}
where $s= \sum \limits_{i=1}^{n-1} s_i$.\\
Then, the absolute value of the determinant of the Jacobian of $\psi$ equals $\frac{1}{2^{n-1}} \cdot \rho^{2n-1}$.
\end{lemma}

Proof: \\
Denote for $i \in \{1,...,2n\}$ by $C_i$ the $i$-th column and by $R_i$ the $i$-th row \mbox{of the Jacobian of $\psi$.} \\
For $i \in \{1,...,n-1 \}$, one has
\begin{eqnarray*}
C_{2i-1}&=& \left(
\begin{matrix}
 0\\
 \vdots \\
 0 \\
 \frac{\rho \cos(t_i)}{2 \sqrt{s_i}} \\
 0 \\
 \vdots \\
 0 \\
 - \rho \sqrt{s_i} \sin(t_i) \\
 0 \\
 \vdots \\
 0 \\
 \sqrt{s_i} \cos(t_i) \\
\end{matrix}
\right),
~~
\begin{matrix}
 ~\\
 ~ \\
 ~ \vspace{0.2cm}\\
 \leftarrow~i-\text{th row}~ \rightarrow \\
 ~ \\
 ~ \\
 ~\\
 ~\vspace{-0.2cm}\\
 \leftarrow~(n-1+i)-\text{th row}~\rightarrow \\
 ~ \\
 ~ \\
 ~\\
 ~ \\
\end{matrix}
~~
C_{2i}~=~ \left(
\begin{matrix}
 0\\
 \vdots \\
 0 \\
 \frac{\rho \sin(t_i)}{2 \sqrt{s_i}} \\
 0 \\
 \vdots \\
 0 \\
 \rho \sqrt{s_i} \cos(t_i) \\
 0 \\
 \vdots \\
 0 \\
 \sqrt{s_i} \sin(t_i) \\
\end{matrix}
\right)
\end{eqnarray*}
and
\begin{eqnarray*}
C_{2n-1}&=& \left(
\begin{matrix}
 -\frac{\rho \cos(t_n)}{2 \sqrt{1-s}}\\
 \vdots \\
 -\frac{\rho \cos(t_n)}{2 \sqrt{1-s}}\\
 0 \\
 \vdots \\
 0 \\
 - \rho \sqrt{1-s} \sin(t_n) \vspace{0.1cm}\\
 \sqrt{1-s} \cos(t_n) \\
\end{matrix}
\right),
~~
\begin{matrix}
 ~\\
 ~ \\
 ~\\
 ~\vspace{-0.3cm}\\
 \leftarrow~(n-1)-\text{th row}~ \rightarrow \\
  ~ \\
 ~ \\
 ~ \\
 ~\\
 ~\\
 ~ \\
\end{matrix}
~~
C_{2n}~=~ \left(
\begin{matrix}
 -\frac{\rho \sin(t_n)}{2 \sqrt{1-s}}\\
 \vdots \\
 -\frac{\rho \sin(t_n)}{2 \sqrt{1-s}}\\
 0 \\
 \vdots \\
 0 \\
 \rho \sqrt{1-s} \cos(t_n) \vspace{0.1cm}\\
 \sqrt{1-s} \sin(t_n) \\
\end{matrix}
\right).
\end{eqnarray*}
Now, in several steps, this matrix will be transformed into a new matrix whose determinant can easily be calculated. For simplicity, the columns and rows of the matrices appearing in each step will also be denoted by $C_i$ and $R_i$ for $i \in \{1,...,2n\}$. \\
First, one takes out the factor $\rho$ in the rows $R_i$ for $i \in \{1,...,2n-1\}$, the factor $\frac{1}{2}$ in the rows $R_i$ for $i \in \{1,...,n-1\}$, the factor $\sqrt{s_i}$ in the columns $C_{2i-1}$ and $C_{2i}$ for every $i \in \{1,...,n-1\}$ and the factor $\sqrt{1- s}$ in the columns $C_{2n-1}$ and $C_{2n}$. Hence, one has the prefactor $\rho^{2n-1} \frac{1}{2^{n-1}} s_1 \cdots s_{n-1} (1-s)$ and the columns of the remaining matrix have the shape
\begin{eqnarray*}
C_{2i-1}&=& \left(
\begin{matrix}
 0\\
 \vdots \\
 0 \vspace{0.1cm} \\
 \frac{\cos(t_i)}{s_i} \vspace{0.1cm}\\
 0 \\
 \vdots \\
 0 \\
 - \sin(t_i) \\
 0 \\
 \vdots \\
 0 \\
 \cos(t_i) \\
\end{matrix}
\right),
~~
\begin{matrix}
 ~\\
 ~ \\
 ~ \vspace{0.2cm}\\
 \leftarrow~i-\text{th row}~ \rightarrow \\
 ~ \\
 ~ \\
 ~\\
 ~\vspace{-0.2cm}\\
 \leftarrow~(n-1+i)-\text{th row}~\rightarrow \\
 ~ \\
 ~ \\
 ~\\
 ~ \\
\end{matrix}
~~
C_{2i}~=~ \left(
\begin{matrix}
 0\\
 \vdots \\
 0 \vspace{0.1cm}\\
 \frac{\sin(t_i)}{s_i} \vspace{0.1cm}\\
 0 \\
 \vdots \\
 0 \\
 \cos(t_i) \\
 0 \\
 \vdots \\
 0 \\
 \sin(t_i) \\
\end{matrix}
\right)
\end{eqnarray*}
for all $i \in \{1,...,n-1\}$ and
\begin{eqnarray*}
C_{2n-1}&=& \left(
\begin{matrix}
 -\frac{\cos(t_n)}{1-s}\\
 \vdots \\
 -\frac{\cos(t_n)}{1-s}\\
 0 \\
 \vdots \\
 0 \\
 -\sin(t_n) \vspace{0.1cm}\\
 \cos(t_n) \\
\end{matrix}
\right),
~~
\begin{matrix}
 ~\\
 ~ \\
 ~\\
  ~\\
 \leftarrow~(n-1)-\text{th row}~\rightarrow\\
  ~ \\
  ~\\
 ~ \\
 ~ \\
 ~\\
 ~\\
 ~ \\
\end{matrix}
~~
C_{2n}~=~ \left(
\begin{matrix}
 -\frac{\sin(t_n)}{1-s}\\
 \vdots \\
 -\frac{\sin(t_n)}{1-s}\\
 0 \\
 \vdots \\
 0 \\
 \cos(t_n) \vspace{0.1cm}\\
 \sin(t_n) \\
\end{matrix}
\right).
\end{eqnarray*}
Next, for every $i \in \{1,...,n\}$, the column $C_{2i-1}$ shall be replaced by $\sin(t_i) C_{2i-1}-\cos(t_i) C_{2i}$. Then, the prefactor changes to $\rho^{2n-1} \frac{1}{2^{n-1}} \frac{s_1 \cdots s_{n-1} (1-s)}{\sin(t_1) \cdots \sin(t_n)}$ and for every $i \in \{1,...,n\}$,
\begin{eqnarray*}
C_{2i-1}&=& \left(
\begin{matrix}
 0\\
 \vdots \\
 0 \\
 -1\\
 0 \\
 \vdots \\
 0 \\
\end{matrix}
\right).
~~
\begin{matrix}
 ~\\
 ~ \\
 ~ \\
  ~ \\
 \leftarrow~(n-1+i)-\text{th row} \\
 ~ \\
 ~ \\
 ~\\
 ~ \\
\end{matrix}
\end{eqnarray*}
The columns $C_{2i}$ for $i \in \{1,...,n\}$ stay the same. \\
Now, for all $i \in \{1,...,n-1\}$, the rows $R_i$ and $R_{n-1+i}$ will be interchanged. Therefore, the prefactor is multiplied by $(-1)^{n-1}$ and for every $i \in \{1,...,n-1\}$,
\begin{eqnarray*}
C_{2i-1}&=& \left(
\begin{matrix}
 0\\
 \vdots \\
 0 \\
 -1\\
 0 \\
 \vdots \\
 0 \\
\end{matrix}
\right),
~~
\begin{matrix}
 ~\\
 ~ \\
 ~ \\
 \leftarrow~i-\text{th row}~~~ \\
 ~ \\
 ~\\
 ~ \\
\end{matrix}
~~
C_{2n-1}~=~ \left(
\begin{matrix}
 0\\
 \vdots \\
 0 \\
 -1\\
 0 \\
\end{matrix}
\right)
\end{eqnarray*}
and
\begin{eqnarray*}
C_{2i}&=& \left(
\begin{matrix}
 \vdots \\
 0 \\
 \cos(t_i) \\
 0 \\
 \vdots \\
 0 \vspace{0.1cm}\\
 \frac{\sin(t_i)}{s_i} \vspace{0.1cm}\\
 0 \\
 \vdots \\
 0 \\
 \sin(t_i) \\
\end{matrix}
\right),
~~
\begin{matrix}
 ~\\
 ~ \\
 ~ \vspace{0.05cm}\\
 \leftarrow~i-\text{th row}~ \\
 ~ \\
 ~n-\text{th row}~\rightarrow~
 ~\\
 ~\vspace{0.3cm}\\
 \leftarrow~(n-1+i)-\text{th row}~ \\
 ~ \\
 ~ \\
 ~\\
 ~ \\
\end{matrix}
~~
C_{2n}~=~ \left(
\begin{matrix}
 \vdots \\
 0 \\
 -\frac{\sin(t_n)}{1-s}\\
 \vdots \\
 -\frac{\sin(t_n)}{1-s} \vspace{0.1cm}\\
  \cos(t_n) \vspace{0.1cm}\\
 \sin(t_n) \\
\end{matrix}
\right).
\end{eqnarray*}
In the next step, for every $i \in \{1,...,n-1\}$, the matrix will be developed with respect to the $i$-th row, which has only one non-zero entry, namely the entry $-1$ in the \mbox{$(2i-1)$-th} column. One develops with respect to the $(2n-1)$-th row which also only consists of one non-zero entry, $-1$, in the $(2n-1)$-th column. The prefactor is then multiplied by
$(-1)^n (-1)^{2n-1+2n-1}\prod \limits_{i=1}^{n-1} (-1)^{i+2i-1}= \prod \limits_{i=1}^{n-1}(-1)^{n+i-1}$, i.e. the prefactor now equals
$$(-1)^{n-1} \prod \limits_{i=1}^{n-1}(-1)^{n+i-1} \rho^{2n-1} \frac{1}{2^{n-1}} \frac{s_1 \cdots s_{n-1} (1-s)}{\sin(t_1) \cdots \sin(t_n)}= \prod \limits_{i=1}^{n-1}(-1)^{i} \rho^{2n-1} \frac{1}{2^{n-1}} \frac{s_1 \cdots s_{n-1} (1-s)}{\sin(t_1) \cdots \sin(t_n)}.$$
 One has a $n \times n$-matrix left, whose columns have the shape
\begin{eqnarray*}
C_{i}&=& \left(
\begin{matrix}
 \vdots \\
 0 \\
 \frac{\sin(t_i)}{s_i} \\
 0 \\
 \vdots \\
 0 \\
 \sin(t_i) \\
\end{matrix}
\right),
~~
\begin{matrix}
 ~\\
 ~ \\
 ~ \\
  ~ \\
 \leftarrow~i-\text{th row}~~~ \\
 ~ \\
 ~ \\
 ~\\
 ~\\
 ~ \\
\end{matrix}
~~
C_{n}~=~ \left(
\begin{matrix}
 -\frac{\sin(t_n)}{1-s}\\
 \vdots \\
 -\frac{\sin(t_n)}{1-s} \vspace{0.1cm}\\
 \sin(t_n) \\
\end{matrix}
\right)
\end{eqnarray*}
for all $i \in \{1,...,n-1\}$. Now, in every column $C_i$ for $i \in \{1,...,n\}$, one can take out the factor $\sin (t_i)$. Then, the prefactor changes to $\prod \limits_{i=1}^{n-1}(-1)^{i} \rho^{2n-1} \frac{1}{2^{n-1}} s_1 \cdots s_{n-1} (1-s)$ and one has the following columns:
\begin{eqnarray*}
C_{i}&=& \left(
\begin{matrix}
 \vdots \\
 0 \vspace{0.1cm}\\
 \frac{1}{s_i} \vspace{0.1cm}\\
 0 \\
 \vdots \\
 0 \\
 1\\
\end{matrix}
\right),
~~
\begin{matrix}
 ~\\
 ~ \\
 ~ \\
  ~ \vspace{0.05cm}\\
 \leftarrow~i-\text{th row}~~~ \\
 ~ \\
 ~ \\
 ~\\
 ~\\
 ~ \\
\end{matrix}
~~
C_{n}~=~ \left(
\begin{matrix}
 -\frac{1}{1-s}\\
 \vdots \\
 -\frac{1}{1-s}\vspace{0.1cm}\\
 1\\
\end{matrix}
\right)
\end{eqnarray*}
for all $i \in \{1,...,n-1\}$. In the last step, the column $C_n$ will be replaced by $C_n + \frac{1}{1-s} \sum \limits_{i=1}^{n-1}s_i C_i$. Since
$$1+\frac{1}{1-s}~ \sum \limits_{i=1}^{n-1}s_i=\frac{1-s}{1-s}+\frac{1}{1-s}~s=\frac{1}{1-s},$$
one obtains the columns
\begin{eqnarray*}
C_{i}&=& \left(
\begin{matrix}
 \vdots \\
 0 \vspace{0.1cm}\\
 \frac{1}{s_i} \vspace{0.1cm}\\
 0 \\
 \vdots \\
 0 \\
 1\\
\end{matrix}
\right),
~~
\begin{matrix}
 ~\\
 ~ \\
 ~ \\
  ~ \vspace{0.05cm}\\
 \leftarrow~i-\text{th row}~~~ \\
 ~ \\
 ~ \\
 ~\\
 ~\\
 ~ \\
\end{matrix}
~~
C_{n}~=~ \left(
\begin{matrix}
 0\\
 \vdots \\
 0 \vspace{0.1cm}\\
 \frac{1}{1-s}\\
\end{matrix}
\right)
\end{eqnarray*}
for $i \in \{1,...,n-1\}$ and the prefactor stays the same, i.e. $\prod \limits_{i=1}^{n-1}(-1)^{i} \rho^{2n-1} \frac{1}{2^{n-1}} s_1 \cdots s_{n-1} (1-s)$. \\
Since the remaining matrix is a triangular matrix, one can easily calculate its determinant and gets
$$\prod \limits_{i=1}^{n-1}(-1)^{i} \rho^{2n-1} \frac{1}{2^{n-1}}~ s_1 \cdots s_{n-1} (1-s)~ \frac{1}{s_1 \cdots s_{n-1} (1-s)}=\prod \limits_{i=1}^{n-1}(-1)^{i} \frac{1}{2^{n-1}} \cdot \rho^{2n-1}.$$
\qed

\begin{corollary}\label{measure on Sn}
~\\
For the map $\ps $ defined in Lemma \ref{det}, the measure defined on the complex sphere $S^{n} $ in $\C^{n} $ by
\begin{eqnarray*}
 \int \limits_{S^n}f(v)d\sigma(v)=\int \limits_{[ 0,1]\times [ 0,2\pi)}f \big(\ps(s_1,\cdots, s_{n-1},t_1,\cdots, t_n,1) \big)ds_1\cdots ds_{n-1}dt_1\cdots t_n\end{eqnarray*}
is the $U(n) $-invariant measure such that for each function $f$ which is continuous on the unit ball in $\C^{n} $,
\begin{eqnarray}\label{second formula}
 \int \limits_{B^{n}}f(z)dz=\int \limits_{0}^{1}\rh^{2n-1}d\rh\int_{S^n}f(\rh v)\frac{d\si(v)}{2^{n-1}}.
 \end{eqnarray}
\end{corollary}

Proof: \\
(\ref{second formula}) gives the decomposition of the Lebesgue measure on $B^{n}\simeq  [0, 1] \times S^{n}$
through $z = \rh v $. Thus, it defines the $U (n) $-invariant measure $d\sigma(v) $ on the unit sphere
$S^{n} $. Moreover, by Lemma \ref{det}
\begin{eqnarray*}
 \int \limits_{[ 0,1]^{n-1}\times [ 0,2\pi)^{n}\times[ 0,1]}f(\ps(s_1,\cdots, s_{n-1},t_1,\cdots, t_n,\rh))ds_1\cdots ds_{n-1}dt_1\cdots t_n\frac{\rh ^{2n-1}}{2^{n-1}}d\rh
 =\int \limits_{B^{n}
 }f(z)dz.
 \end{eqnarray*}

Since $\ps(s_1 , . . . , s_{n-1}, t_1 , . . . , t_n , \rh) = \rh\ps(s_1, . . . , s_{n-1}, t_1, . . . , t_n , 1) $ and $ \ps(s_1, . . . , s_{n-1}, t_1, . . . , t_n,1)$ is normed, this proves the corollary.
~\\
\qed 
~\\
~\\
\textbf{Acknowledgements:} \\
The authors would like to thank the referee for his / her careful reading of our paper and many valuable suggestions. \\
Janne-Kathrin G\"unther was supported for this work by the Fonds National de la Recherche, Luxembourg (Project Code 3964572).

\section{References}
\begin{enumerate}
\item \label{Ba} L.W.Baggett, A description of the topology on the dual spaces of certain locally compact groups, Transactions of the American Mathematical Society 132, pp.175-215, 1968.
\item \label{bel-bel-lud} D.Beltita, I.Beltita and J.Ludwig, Fourier transforms of $C^*$-algebras of nilpotent Lie groups, to appear in: International Mathematics Research Notices, doi:10.1093/imrn/rnw040.
\item \label{BGLR} C.Benson, J.Jenkins, R.Lipsman and G.Ratcliff, A geometric criterion for Gelfand pairs associated with the Heisenberg group,
Pacific Journal of Mathematics 178, \mbox{no.1}, \mbox{pp.1-36}, 1997.
\item\label{BJR} C.Benson, J.Jenkins and G.Ratcliff, Bounded $K$-spherical functions on Heisenberg groups, Journal of Functional Analysis 105, pp.409-443, 1992.
\item \label{BJRW} C.Benson, J.Jenkins, G.Ratcliff and T.Worku, Spectra for Gelfand pairs associated with the Heisenberg group, Colloquium Mathematicae 71, pp.305-328, 1996.
\item \label{brown} I.Brown, Dual topology of a nilpotent Lie group, Annales scientifiques de l'\'{E}.N.S. $4^e$ s\'{e}rie, tome 6, no.3, pp.407-411, 1973.
\item \label{cohn} L.Cohn, Analytic Theory of the Harish-Chandra C-Function, Springer-Verlag, Berlin-Heidelberg-New York, 1974.
\item \label{cor-green} L.Corwin and F.P.Greenleaf, Representations of nilpotent Lie groups and their applications. Part I. Basic theory and examples, Cambridge Studies in Advanced Mathematics, vol.18, Cambridge University Press, Cambridge, 1990.
\item \label{dix} J.Dixmier, $C^*$-algebras. Translated from French by Francis Jellett, North-Holland Mathe\-ma\-ti\-cal Library, vol.15, North-Holland Publishing Company, Amsterdam-New York-Oxford, 1977.
\item \label{Dix-Mal} J.Dixmier and P.Malliavin, Factorisations de fonctions et de vecteurs ind\'efiniment diff\'erentiables, Bulletin des Sciences Math\'ematiques (2) 102, no.4, pp.307-330, 1978.
\item \label{elloumi-dr} M.Elloumi, Espaces duaux de certains produits semi-directs et noyaux associ\'es aux orbites plates, Ph.D. thesis at the Universit\'e de Lorraine, 2009.
\item \label{El-Lu} M.Elloumi and J.Ludwig, Dual topology of the motion groups $SO(n)\ltimes\R^n$, Forum Mathematicum 22, pp.397-410, 2010.
\item \label{fell} J.M.G.Fell, The structure of algebras of operator fields, Acta Mathematica 106, \mbox{pp.233-280}, 1961.
\item \label{Folland} G.B.Folland, Harmonic Analysis in Phase Space, Princeton University Press, 1989.
\item \label{Fu-Ha} W.Fulton and J.Harris, Representation theory, Readings in Mathematics, Springer-Verlag, 1991.
\item \label{janne2} J.-K.G\"unther, The $C^*$-algebra of $SL(2, \R)$, arXiv:1605.09256, 2016.
\item \label{janne1} J.-K.G\"unther and J.Ludwig, The $C^*$-algebras of connected real two-step nilpotent Lie groups, Revista Matem\'atica Complutense 29(1), pp.13-57, 10.1007/s13163-015-0177-7, 2016.
\item \label{hol-bie} W.J.Holman III and L.C.Biedenharn, The Representations and Tensor Operators of the Unitary Groups $U(n)$, in: Group Theory and its Applications, Volume 2, edited by E.M.Loebl, Academic Press, Incorporation, London, 1971.
\item \label{Howe} R.Howe, Quantum mechanics and partial differential equations, Journal of Functional Analysis 38, pp.188-255, 1980.
\item \label{knapp} A.Knapp, Representation Theory of Semisimple Groups. An Overview Based on Examples, Princeton University Press, Princeton, New Jersey, 1986.
\item \label{lahi} R.Lahiani, Analyse Harmonique sur certains groupes de Lie \`{a} croissance polynomiale, Ph.D. thesis at the University of Luxembourg and the Universit\'{e} Paul Verlaine-Metz, 2010.
\item \label{lang} S.Lang, $SL_2(\R)$, Graduate Texts in Mathematics 105, Springer-Verlag, New York, 1985.
\item \label{leplud} H.Leptin and J.Ludwig, Unitary representation theory of exponential Lie groups, De Gruyter Expositions in Mathematics 18, 1994.
\item \label{lin-ludwig} Y.-F.Lin and J.Ludwig, The $C^*$-algebras of $ax+b$-like groups, Journal of Functional Analysis 259, pp.104-130, 2010.
\item \label{Lipsman} R.L.Lipsman, Orbit theory and harmonic analysis on Lie groups with co-compact nilradical, Journal de Math\'ematiques Pures et Appliqu\'ees, tome 59, pp.337-374, 1980.
\item \label{ludwig-turowska} J.Ludwig and L.Turowska, The $C^*$-algebras of the Heisenberg Group and of thread-like Lie groups, Mathematische Zeitschrift 268, no.3-4, pp.897-930, 2011.
\item \label{ludwig-zahir} J.Ludwig and H.Zahir, On the Nilpotent $*$-Fourier Transform, Letters in Mathematical Physics 30, pp.23-24, 1994.
\item \label{Ma} G.W.Mackey, Unitary group representations in physics, probability and number theo\-ry, Benjamin-Cummings, 1978.
\item \label{puk} L.Pukanszky, Le\c{c}ons sur les repr\'{e}sentations des groupes, Dunod, Paris, 1967.
\item \label{hedidr} H.Regeiba, Les $C^*$-alg\`{e}bres des groupes de Lie nilpotents de dimension $\leq 6$, Ph.D. thesis at the Universit\'{e} de Lorraine, 2014.
\item \label{hedi} H.Regeiba and J.Ludwig, $C^*$-Algebras with Norm Controlled Dual Limits and Nilpotent Lie Groups, arXiv: 1309.6941, 2013.
\item \label{wal} N.Wallach, Real Reductive Groups I, Academic Press, Pure and Applied Mathematics, San Diego, 1988.
\item \label{walII} N.Wallach, Real Reductive Groups II, Academic Press, Pure and Applied Mathematics, San Diego, 1992.
\item \label{was} A.Wassermann, Une d\'emonstration de la conjecture de Connes-Kasparov pour les groupes de Lie lin\'eaires connexes r\'eductifs, Comptes Rendus de l'Acad\'emie des \mbox{Sciences}, Paris Series I Mathematics 304, no.18, pp.559-562, 1987.
\end{enumerate}

\end{document}